\documentclass[fleqn, 11pt]{elsarticle}

\PassOptionsToPackage{hyphens}{url}
\usepackage{xurl}
\usepackage{hyperref}

\journal{{\color{black}EJOR}}

\usepackage{pifont}

\newcommand{\cmark}{\ding{51}} 
\newcommand{\xmark}{\ding{55}} 

\makeatletter
\def\ps@pprintTitle{%
  \let\@oddhead\@empty
  \let\@evenhead\@empty
  \let\@oddfoot\@empty
  \let\@evenfoot\@oddfoot
}
\makeatother
\usepackage[usenames,dvipsnames]{color}

\let\underbar\relax
\usepackage[utf8]{inputenc}
\usepackage[english]{babel}
\usepackage{amsthm}
\usepackage{makecell}
\usepackage{thmtools}
\usepackage{thm-restate}
\usepackage{tabularx}
\usepackage{adjustbox}
\usepackage{blindtext}
\usepackage{amsmath,amssymb,amsfonts,amsthm}
\usepackage{cleveref}
\usepackage{fancyhdr}
\usepackage[table]{xcolor}
\usepackage{tikz}
\usetikzlibrary{shapes.geometric, positioning}
\usetikzlibrary{backgrounds}
\usetikzlibrary{arrows.meta, fit}
\usetikzlibrary{matrix}
\usepackage{pgfplots}
\usepackage{pgfplotstable}
\usepackage{pgfplots}\usetikzlibrary{patterns}
\usepgfplotslibrary{groupplots}
\usepackage{siunitx}
\usepackage[linesnumbered,ruled,vlined]{algorithm2e}

\usepackage{relsize}

\usepackage[short]{optidef}
\usepackage{float}
\graphicspath{{figures/}}
\usepackage{amssymb}
\usepackage{color}
\usepackage{mathtools}
\usepackage{subcaption}
\usepackage{changepage}
\usepackage{geometry}
\usepackage{lipsum}

\usepackage{multirow}
\usepackage{longtable}
\usepackage{supertabular}

\usetikzlibrary{positioning}
\usepgfplotslibrary{dateplot}
\pgfplotsset{compat=1.18}

\usepackage{graphics} 
\usepackage{setspace} 
\usepackage{graphicx} 
\usepackage{lscape}
\usepackage{makeidx} 
\usepackage{pdflscape} 
\usepackage{amsmath,amssymb,amsfonts} 
\usepackage{nccmath}
\usepackage{mathptmx} 
\usepackage[final]{pdfpages} 
\usepackage{calc} 
\usepackage{fancyhdr} 
\usepackage{multicol} 
\usepackage[T1]{fontenc}
\usepackage[utf8]{inputenc} 
\usepackage[english]{babel}
\usepackage{threeparttable, supertabular} 
\usepackage{endnotes} 
\usepackage{rotating, rotfloat} 
\usepackage{fullpage} 
\usepackage[usenames,dvipsnames]{color} 
\usepackage{colortbl} 
\usepackage{eurosym}
\usepackage[bottom, flushmargin]{footmisc} 
\usepackage{booktabs} 
\usepackage{tabularx}
\usepackage{mathrsfs}
\usepackage[format=plain, textfont=normal, justification=justified, singlelinecheck=false, font=normal]{caption}

\hypersetup{hidelinks}
\usepackage{enumitem}
\usepackage{diagbox}

\allowdisplaybreaks

\usepackage{etoolbox}

\makeatletter
\def\urlprefix{}
\def\DOIprefix{}
\def\newblock{\unskip\space}
\makeatother

\patchcmd{\thebibliography}{\newblock}{\space}{}{}

\pgfdeclarepatternformonly{sparsecrosshatch}
{\pgfpoint{0pt}{0pt}}{\pgfpoint{10pt}{10pt}}{\pgfpoint{10pt}{10pt}}
{
  \pgfsetlinewidth{0.32pt}
  \pgfpathmoveto{\pgfpoint{0pt}{0pt}}
  \pgfpathlineto{\pgfpoint{10pt}{10pt}}
  \pgfusepath{stroke}
  \pgfpathmoveto{\pgfpoint{0pt}{10pt}}
  \pgfpathlineto{\pgfpoint{10pt}{0pt}}
  \pgfusepath{stroke}
}

\pgfdeclarepatternformonly{sparsenorthwestlines}
{\pgfpoint{0pt}{0pt}}{\pgfpoint{10pt}{10pt}}{\pgfpoint{10pt}{10pt}}
{
  \pgfsetlinewidth{0.32pt}
  \pgfpathmoveto{\pgfpoint{0pt}{10pt}}
  \pgfpathlineto{\pgfpoint{10pt}{0pt}}
  \pgfusepath{stroke}
}

\pgfdeclarepatternformonly{sparsenortheastlines}
{\pgfpoint{0pt}{0pt}}{\pgfpoint{10pt}{10pt}}{\pgfpoint{10pt}{10pt}}
{
  \pgfsetlinewidth{0.32pt}
  \pgfpathmoveto{\pgfpoint{0pt}{0pt}}
  \pgfpathlineto{\pgfpoint{10pt}{10pt}}
  \pgfusepath{stroke}
}

\pgfdeclarepatternformonly{sparsehorizontallines}
  {\pgfpoint{0pt}{0pt}}{\pgfpoint{8pt}{8pt}}{\pgfpoint{8pt}{8pt}}{
    \pgfsetlinewidth{0.32pt}
    \pgfpathmoveto{\pgfpoint{0pt}{6pt}}
    \pgfpathlineto{\pgfpoint{8pt}{6pt}}
    \pgfusepath{stroke}
}

\makeatletter
\tikzset{
  dot distance/.store in=\pgf@pattern@distance, dot distance=4pt,
  dot radius/.store in=\pgf@pattern@radius,   dot radius=0.6pt
}
\pgfdeclarepatternformonly[\pgf@pattern@distance,\pgf@pattern@radius]{sparsedots}
  {\pgfpoint{0pt}{0pt}}
  {\pgfpoint{\pgf@pattern@distance}{\pgf@pattern@distance}}
  {\pgfpoint{\pgf@pattern@distance}{\pgf@pattern@distance}}
  {
    \pgfpathcircle{\pgfpoint{.5\pgf@pattern@distance}{.5\pgf@pattern@distance}}{\pgf@pattern@radius}%
    \pgfusepath{fill}
  }
\makeatother

\begin{document}

\begin{frontmatter}

\title{A Dynamic Strategic Plan for Transition to Campus-Scale Clean Electricity Using Multi-Stage Stochastic Programming}

\author[1]{Ahmet Emir \c{S}ener}

\ead{ahmet.sener@std.bogazici.edu.tr}

\address[1]{Department of Industrial Engineering, Bo\u{g}azi\c{c}i University, Istanbul, Turkey}

\author[2]{Burak Kocuk}
\ead{burak.kocuk@sabanciuniv.edu}

\author[2]{Tu\u{g}\c{c}e Y\"{u}ksel}
\ead{tugce.yuksel@sabanciuniv.edu}

\address[2]{Faculty of Engineering and Natural Sciences, Sabanc{\i} University, Istanbul, Turkey}

\begin{abstract} 
The decarbonization of energy systems at energy-intensive sites is an essential component of global climate mitigation, yet such transitions involve substantial capital requirements, ongoing technological progress, and the operational complexities of renewable integration. This study presents a dynamic strategic planning framework that applies multi-stage stochastic programming to guide clean electricity transitions at the campus level. The model jointly addresses technology investment, storage operation, and grid interaction decisions while explicitly incorporating uncertainties in future technology cost trajectories and efficiency improvements. By enabling adaptive, stage-wise decision-making, the framework provides a structured approach for large electricity consumers seeking to achieve self-sufficient and sustainable energy systems. The approach is demonstrated through a case study of Middle East Technical University (Ankara, Turkey), which has committed to achieving carbon-neutral electricity by 2040. Through the integration of solar photovoltaics, wind power, and lithium-ion batteries, the model links long-term investment planning with operational-level dynamics by incorporating high-resolution demand and meteorological data. 
Our findings from the case study, sensitivity analyses, and comparisons with simplified models indicate that accounting for uncertainty and temporal detail is crucial for both the economic viability and operational feasibility of campus-scale clean electricity transitions.
\end{abstract}

\begin{keyword}
\texttt{Clean energy transition, Strategic energy planning, Multi-stage stochastic programming, Mixed-integer linear programming}
\end{keyword}

\end{frontmatter}

\section{Introduction}

Global energy demand continues to increase, driven by population growth and sustained economic expansion. In 2023, primary energy consumption rose by approximately 2\% relative to the previous year~\cite{McGrath_Hunt_Taylor_Fitzgerald_2025}. While the global energy mix is gradually diversifying through greater adoption of renewables, nuclear power, and emerging options such as hydrogen and tidal energy, fossil fuels still supplied 81.5\% of total demand.
This represents only a modest 0.4\% decline from 2022, yet absolute fossil fuel use is still increasing, thereby aggravating global carbon emissions. 
In recognition of these trends, governments agreed to work together to transition away from fossil fuels and increase the penetration of renewable energy systems to achieve net zero emissions by 2050~\cite{arora2024cop28}. Within this effort, distributed energy systems, in which the electricity is generated where it is consumed, can play a  pivotal role  in fastening this transition.




Among these systems, large-scale facilities and campuses are particularly important due to their substantial electricity demand and reliance on fossil fuel–based grids. Such entities also offer some valuable opportunities for distributed energy generation. 
For instance, roofs on campus buildings can be used for small-scale solar panels, while  idle land spaces can be allocated for larger-scale solar and wind farms. 
However, transitioning campuses to on-site renewable technologies is both capital-intensive and operationally complex, making a one-step shift impractical. Instead, a phased strategic plan is required—one that identifies optimal portfolios of electricity generation and storage technologies, schedules their deployment while ensuring reliable operations.

There are two important considerations while making this strategic plan: i) techno-economic trends, and ii) operational details.
Ongoing techno-economical progress enhances renewable solutions by improving efficiency, lowering capital and operating costs, and increasing equipment reliability. However, these advances are subject to considerable uncertainty. On the one hand, we need strategic planning methodologies for long-term investment decisions that take into account this uncertainty.  On the other hand, since renewable generation is intermittent, operational feasibility should also be guaranteed at a high temporal resolution to meet campus electricity demand.
%

This study aims to obtain a dynamic strategic plan for transitioning to campus-scale clean electricity systems. 
We propose a multi-stage stochastic programming model that guides decisions on technology installations, salvage decisions, grid purchases, and storage operations, with the objective of minimizing the expected total transition cost. 
By explicitly representing uncertainties in technological advancement and linking long‑term strategic decisions with short‑term operational dynamics, our framework enables adaptive and reliable planning. 
We demonstrate the applicability of our approach through a detailed case study of Middle East Technical University (METU) in Ankara, Turkey, which has committed to generate its own electricity using renewable technologies  by 2040~\cite{METUPolicies}. 

\subsection{Literature Review}





Our work is closely related to two well-studied, strategic decision-making problems from the literature: Strategic Energy Planning (SEP) and Clean Energy Transition (CET). SEP provides a long-term, high-level framework that guides energy systems toward desired future configurations by optimizing the timing, scale, and composition of investments and operations. Beyond cost minimization, SEP often integrates secondary objectives such as decarbonization and energy security by incorporating technological trends, policy directives, and operational dynamics across scales ranging from individual households \cite{ling2024} to national power systems. Within this broader domain, CET has emerged as a growing subcategory aimed at achieving 100\% renewable electricity supply.


Table~\ref{tab:litreview} summarizes selected SEP and CET studies. While many adopt extended and phased planning horizons, they typically rely on deterministic formulations that do not explicitly capture uncertainty in technological progress.
Among the subset of studies that adopt stochastic methods, uncertainties in fuel prices, demand growth, or policy frameworks are more commonly emphasized \cite{Guevara2022,Lei2021}, whereas uncertainties in technological cost and efficiency trajectories are underrepresented. Those that do consider technological uncertainties generally omit efficiency improvements~\cite{Cano2016} and tend to rely on aggregated temporal resolutions or representative days~\cite{Ioannou2019,rathi2022}. This reliance can obscure critical operational mismatches between variable renewable generation and fluctuating electricity demand—a challenge that becomes increasingly severe as renewable penetration deepens~\cite{gao2023}. Moreover, reference \cite{Marcy2022} shows that such representative-hour approaches systematically underestimate peak-day ramping needs and storage requirements, reinforcing the necessity of high-resolution chronological data to accurately capture short-term operational dynamics.

A related body of research is capacity expansion planning (CEP), which addresses the cost-effective deployment of generation, storage, and transmission infrastructure to reliably meet future energy demand. MSSP is frequently employed in CEP studies to support sequential investment decisions under various uncertainties \cite{Zou2018, park2020}, often leveraging decomposition methods \cite{Singh2009, Rebennack2016}. While prior studies explore diverse uncertainty sources, they often rely on low temporal resolution and provide limited consideration to long-term technological change \cite{Hole2023, Liu2018}.

\begin{table}[H]
  \centering
  \normalsize
  \setlength{\tabcolsep}{7.7pt}
  \caption{An overview of SEP and CET literature.}
  \begin{adjustbox}{max width=0.9\textwidth}
    \begin{tabularx}{\textwidth}{
      c
      l
      l
      *{2}{>{\centering\arraybackslash}c}
      *{2}{>{\centering\arraybackslash}c}
      >{\centering\arraybackslash}c
    }
      \toprule
      &&& \multicolumn{2}{c}{\textbf{Technology}}
      & \multicolumn{2}{c}{\textbf{Temporal resolution}} &\\
      \cmidrule(lr){4-5}\cmidrule(lr){6-7}
      \textbf{Reference} & \textbf{Scope} & \textbf{Opt. type}
      & \textbf{Cost} & \textbf{Efficiency}
      & \textbf{Investment} & \textbf{Operational} & \textbf{Horizon} \\
      \midrule
      \cite{Husein2018} & SEP & Heuristic & D & – & One-phase & Hourly & 25 yrs \\
      \cite{Prina2019}& SEP & Heuristic & D & – & 5-year & Hourly & 35 yrs \\
      \cite{Gils2017}& SEP & D MILP & D & – & One-phase & Hourly & 1 yr \\
      \cite{Moret2020}& SEP & RO & D & – & One-phase & Monthly & 1 yr \\
      \cite{Cano2016}& SEP & MSSP & S & – & Annual & Hourly & 16 yrs \\
      \cite{Guevara2022}& SEP & MSSP & D & – & 10-year & Monthly & 30 yrs \\
      \cite{Heuberger2017}& SEP\textsuperscript{*} & D MILP & D & – & Annual & Rep. days & 35 yrs \\
      \cite{Powell2012}& SEP\textsuperscript{*} & ADP & D & – & Annual & Hourly & 20 yrs \\
      \cite{Ioannou2019}& SEP\textsuperscript{*} & MSSP & S & – & Annual & Annual & 15 yrs \\
      \cite{rathi2022}& SEP\textsuperscript{*} & MSSP & S & – & 5-year & Rep. day & 40 yrs \\
      \cite{Child2019}& CET & D MILP & D & – & 5-year & Hourly & 35 yrs \\
      \cite{Mavromatidis2021}& CET & D MILP & D & D & 5-year & Rep. days & 30 yrs \\
      \cite{Tian2022}& CET & D MINLP & D & – & Annual & Monthly & 15 yrs \\
      \cite{Zhao2021}& CET & MSARO & D & – & Annual & Annual & 20 yrs \\
      \cite{Lei2021}& CET & MSSP & D & – & 5-year & Rep. days & 30 yrs \\
      \textbf{This study} & CET & MSSP & S & S & Annual & Bihourly & 15 yrs \\
      \bottomrule
    \end{tabularx}
  \end{adjustbox}
  \label{tab:litreview}
\end{table}

\vspace{-0.4cm}
{\footnotesize\setlength{\baselineskip}{0.8\baselineskip}\noindent
\textbf{D MILP:} Deterministic Mixed-Integer Linear Programming, \textbf{D MINLP:} Deterministic Mixed-Integer Nonlinear Programming, \textbf{MSSP:} Multi-stage Stochastic Programming, \textbf{RO:} Robust Optimization, \textbf{MSARO:} Multi-stage Adaptive Robust Optimization, \textbf{ADP:} Approximate Dynamic Programming, \textbf{D:} Deterministic, \textbf{S:} Stochastic, {\textbf{Rep.:} Representative}, \textbf{Note:} SEP entries marked with \textbf{\textsuperscript{*}} refer to studies that support the clean energy transition through policy mechanisms—such as carbon taxes—but do not explicitly target a 100\% renewable energy transformation.\par}

At the campus scale, energy transition research has increasingly emphasized the integration of renewable generation and microgrid design \cite{Abulibdeh2024}. Deterministic MILP models are commonly used to schedule generation, storage, and demand response under varying tariff and outage conditions \cite{Li2022}, whereas MINLP formulations have been applied to optimize integrated electricity, heating, and cooling systems, as demonstrated in Cornell University’s energy redesign \cite{Tian2022}. Complementary approaches adopt simulation-based techno-economic analyses to size systems, evaluating both financial performance and environmental outcomes \cite{Paspatis2022,AlKassem2022}. Hybrid frameworks have also been proposed that link operational dispatch with financial evaluation, enabling sensitivity analyses under alternative policy and incentive settings \cite{Husein2018}. More recently, multi-objective metaheuristic methods have been introduced for microgrid design, jointly optimizing renewable generation and storage to balance cost, reliability, and sustainability \cite{Mishra2025}. While most of these campus-scale studies employ high temporal resolution, none account for uncertainty in the future cost and efficiency trajectories of emerging technologies.

In summary, studies that simultaneously incorporate technological uncertainty and high-resolution operational dynamics remain scarce. To the best of our knowledge, these issues are entirely overlooked in the context of campus-scale transition studies. 
Our contribution is to address this gap by proposing a multi-stage stochastic programming framework that models investment, operation, and decommissioning of solar PV, wind, and battery storage technologies under uncertainty in both cost and efficiency. In contrast to prior studies, our model adopts bihourly operational resolution and annual investment decisions over a 15-year horizon. This enables the formulation of more adaptive and realistic transition strategies for large-scale consumers, such as university campuses aiming to achieve carbon neutrality.


\subsection{Approach and Contributions}

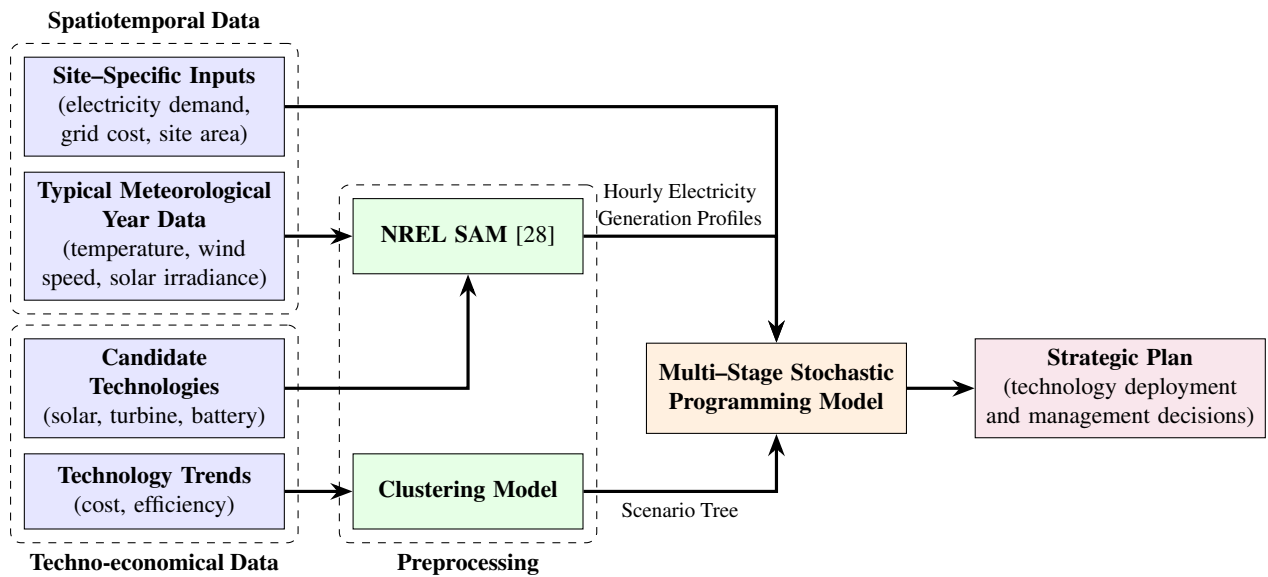
\begin{figure}[H]
  \centering
  \begin{tikzpicture}[
      font=\footnotesize,
      >=Stealth,
      data/.style  = {rectangle, draw, fill=blue!10,
                      text width=32mm, align=center,
                      minimum height=10mm},
      proc/.style  = {rectangle, draw, fill=green!10,
                      text width=28mm, align=center,
                      minimum height=10mm},
      model/.style = {rectangle, draw, fill=orange!12,
                      text width=32mm, align=center,
                      minimum height=12mm},
      plan/.style  = {rectangle, draw, fill=purple!10,
                      text width=36mm, align=center,
                      minimum height=10mm},
      flow/.style  = {->, very thick},
      node distance = 2mm and 16mm
    ]
        \node[data] (csi)   {\textbf{Site--Specific Inputs}\\ (electricity demand, grid cost, site area)};

    \node[data, below=of csi] (tmy)                  {\textbf{Typical Meteorological}\\\textbf{Year Data}\\(temperature, wind speed, solar irradiance)};
    \node[data, below=5mm of tmy] (ctiv)   {\textbf{Candidate Technologies}\\
    (solar, turbine, battery)};
    \node[data, below=of ctiv] (hiced)  {\textbf{Technology Trends}  \\ (cost, efficiency)};

    \node[proc, right=9mm of tmy]  (nrel)  {\textbf{NREL SAM} \cite{sam2025}};
    \node[proc, right=9mm of hiced] (clust){\textbf{Clustering Model}};

    \node[model, right=48mm of ctiv] (msp) {\textbf{Multi--Stage Stochastic}\\\textbf{Programming Model}};

        \node[plan, right=9mm of msp] (plan) {\textbf{Strategic Plan}\\ (technology deployment and management decisions)};

    \draw[flow] (tmy)  -- (nrel);     
    \draw[flow] (ctiv) -| (nrel);     
    \draw[flow] (hiced)-- (clust);
    \draw[flow] (csi)  -| (msp);
        \draw[flow] (msp)  -- (plan);

    \draw[flow] (nrel) -| node[above,pos=0.25]
          {\scriptsize \shortstack{Hourly Electricity \\ Generation Profiles}} (msp);

    \draw[flow] (clust) -| node[below,pos=0.25]
          {\scriptsize Scenario Tree} (msp);

    \begin{pgfonlayer}{background}
      \node[draw, dashed, rounded corners, inner sep=5pt,
            fit=(tmy)(csi),
            label=above:{\footnotesize\bfseries  Spatiotemporal Data}] {};
            
      \node[draw, dashed, rounded corners, inner sep=5pt,
            fit=(ctiv)(hiced),
            label=below:{\footnotesize\bfseries Techno-economical Data}] {};
            
      \node[draw, dashed, rounded corners, inner sep=5pt,
            fit=(nrel)(clust),
            label=below:{\footnotesize\bfseries Preprocessing}] {};
    \end{pgfonlayer}
  \end{tikzpicture}
  \caption{Overview of   the proposed planning framework, illustrating data sources, preprocessing modules, and their integration into the multi‑stage stochastic programming model.}
  \label{fig:contextdiagram}
\end{figure}

This study develops a multi-stage stochastic programming model to provide a long-term, adaptive strategy for campus-scale clean electricity transition. The framework simultaneously aims to ensure a self-sufficient electricity system, minimize expected transition costs under uncertainty, and determine the optimal timing and mix of technology investments. These objectives directly address practical challenges faced by institutions—achieving sustainability targets under emission-reduction commitments, managing long-term budgetary constraints, and coping with uncertainties in technological advancements and cost trajectories. The proposed approach guides decision-makers toward flexible and cost-effective strategies that remain responsive to evolving technological pathways while ensuring alignment with site-specific sustainability goals. Figure~\ref{fig:contextdiagram} provides a graphical overview of the methodological framework, with the corresponding steps described below.
\begin{itemize}[nosep,leftmargin=16.5pt]
\item \textbf{Data Acquisition:} Assembly   of high‑resolution meteorological records, candidate technology parameters, historical technology cost–efficiency series, and other site‑specific inputs. 

\item \textbf{High-Resolution Electricity Generation Simulation:} Derivation of hourly electricity generation profiles for each candidate renewable technology to capture short‑term variability using National Renewable Energy Laboratory’s System Advisor Model (NREL SAM) \cite{sam2025}.

\item \textbf{Technological Advancement Scenario Generation:} Construction of probabilistic scenario tree by applying clustering methods to historical technology cost and efficiency trends.

\item \textbf{Multi-Stage Stochastic Programming:} Incorporation of scenario trees into multi-stage stochastic programming model that jointly optimizes technology choice, investment timing, battery-storage operation, and grid‑electricity procurement at high temporal resolution.

\item \textbf{Case Study and Sensitivity Analysis:} 
{The framework is implemented for METU, which has adopted a climate action plan targeting a 100\% clean electricity generation  by 2040~\cite{METUPolicies}. The analysis explores investment and operational pathways for solar, wind, and battery technologies in line with METU’s sustainability goals, and includes sensitivity analyses under varying emission limits, financial constraints, and operational settings.}

\end{itemize}

The results of our case 
reveal interesting insights: 
i) Site‑specific conditions make solar power   
preferable over wind turbines for METU, even under conservative solar technology improvement scenarios. 
ii) Battery investments are typically deferred to later years, aligning with projected cost declines and rising renewable penetration. 
iii) Across all scenarios, excess generation capacity combined with large battery installations is needed to balance mismatches between electricity demand and renewable output. iv) As renewable penetration increases, the marginal cost of meeting residual demand with clean technologies rises, driven by the structural mismatch between electricity demand profiles and renewable generation patterns, along with the growing need for battery deployment. The \textit{Emission Allowance Case} sensitivity analysis demonstrates that allowing just 1\% of demand to be supplied by conventional grid resources reduces the total system cost by 20.5\%. 
v) Both the highly temporally aggregated and deterministic models fail to satisfy demand and emission constraints, with the deterministic model additionally violating budget constraints.

%
%

The primary contributions and innovations of the proposed approach are as follows:

\begin{itemize}[nosep,leftmargin=16.5pt]

\item 
Contrary to the majority of SEP and CET studies which consider a deterministic technology improvement trajectory, 
we develop
a multi-stage stochastic programming model 
that explicitly take into account the uncertainty in technology development in terms of cost and efficiency.

\item
Although cost and efficiency developments are highly correlated, the stochasticity
in efficiency improvements of clean electricity generation technologies is overlooked in the literature. To  forecast cost and efficiency developments together, we propose a new approach, in which scenario trees  are generated through clustering methods applied to historical data. These trees represent plausible trajectories of renewable electricity technology advancements, enabling uncertainty-aware planning and establishing a realistic foundation for strategic decision-making.

\item Unlike conventional methodologies that aggregate temporal data, the proposed framework employs high-resolution temporal profiles of electricity generation and consumption. Hourly renewable electricity generation profiles are simulated using NREL SAM \cite{sam2025}. This allows us to realistically model  detailed operational aspects in our multi-stage stochastic program, enhancing the reliability and feasibility of strategic-level decisions.

\end{itemize}



By combining site-specific data, technological uncertainty, and detailed operational dynamics, the framework delivers investment strategies that are both economically viable and practically implementable. While demonstrated on a campus-scale electricity transition at METU, the methodology is broadly applicable to other sectors—such as transportation \cite{karimi2024}, industrial manufacturing, and smart grid modernization—where long-term planning must account for technological uncertainty.

\subsection{Paper Organization}

The remainder of this paper is structured as follows. Section~\ref{sec:method} describes the problem formulation and proposes the   multi-stage stochastic programming model. 
Section~\ref{sec:data} provides the data sources, preprocessing methodologies, and scenario generation process, including the construction of probability distributions and scenario trees representing technological advancements. Section~\ref{sec:computational_experiments} presents comprehensive computational experiments based on the METU case study, including extensive sensitivity analyses evaluating the adaptability of the proposed approach under various cases. Finally, Section~\ref{sec:conclusion} summarizes the key findings and outlines potential directions for future research.

\section{Multi-Stage Stochastic Programming Model}
\label{sec:method}


This section describes the methodological framework developed for campus-scale clean electricity transition planning under technological uncertainty. A multi-stage stochastic programming model is formulated to capture both the ambiguity of future technological pathways and the high temporal resolution of system operations.
There are three time levels in our model: i) stages, ii) investment periods, iii) operational sub-periods.
The scenario tree (Section~\ref{subsec:Scenario_Tree_Generation}) represents plausible trajectories of technological evolution across multiple stages. Each stage comprises several time periods during which investment decisions are made. These periods are further subdivided into sub-periods to account for detailed operational decisions and system constraints. To retain computational tractability while preserving practical relevance, the model relies on a set of simplifying assumptions described below.

\noindent\textbf{Stage-wise assumptions and simplifications:}

\begin{itemize}[nosep,leftmargin=16.5pt]

\item Technological advancements occur exclusively at stage transitions. Consequently, each stage is characterized by distinct technological parameters and stochastic conditions.

\end{itemize}

\noindent\textbf{Period-wise assumptions and simplifications:}

\begin{itemize}[nosep,leftmargin=16.5pt]

\item For each renewable electricity generation technology, multiple versions—characterized by differing rated capacities, installation costs, and performance metrics—are considered for investment. All versions of a given technology follow an identical stochastic trajectory for cost and efficiency evolution.

\item Investment and retirement decisions for electricity generation technologies are modeled by integer variables to reflect discrete choice among technology versions. Storage technology deployment, by contrast, is modeled through continuous variables to reflect the unit‑based installation of a single battery version.

\item Newly installed technologies become operational in the same period in which they are commissioned. Conversely, technologies designated for retirement are decommissioned at the start of the specified retirement period.

\item The investment portfolio is restricted to renewable generation and battery storage technologies; conventional fossil fuel-based alternatives are deliberately excluded from consideration.



\item All costs incurred over the planning horizon are discounted to present value using a predefined real discount rate.

\end{itemize}

\noindent\textbf{Sub-period-wise assumptions and simplifications:}

\begin{itemize}[nosep,leftmargin=16.5pt]


\item Electricity demand, grid prices, and renewable generation are specified as deterministic and exogenous inputs in each sub-period. A supplementary sensitivity analysis introduces a safety demand margin to assess system reliability against potential variability in demand and generation levels.



\end{itemize}

Tables~\ref{tab:Sets}, \ref{tab:Variables}, and \ref{tab:Parameters} provide a comprehensive overview of the sets, decision variables, and parameters employed in the subsequent mathematical formulation.

\begin{table}[H]
    \centering
    \caption{List of index sets.}
    \label{tab:Sets}
    \small
    \begin{tabularx}{\textwidth}{|c|X|}
    \hline
    \textbf{Set} & \textbf{Description} \\
    \hline
    $\mathcal{N}$ & Set of nodes in the scenario tree, $\mathcal{N} = \{0, 1, \ldots, N\}$. \\
    $\mathcal{T}$ & Set of investment periods, $\mathcal{T} = \{0, 1, \ldots, T\}$. \\
    $\mathcal{T}_n$ & Subset of investment periods associated with node $n \in \mathcal{N}$. \\
    $\mathcal{Q}$ & Set of operational sub-periods within each investment period, $\mathcal{Q} = \{0, 1, \ldots, Q\}$. \\
    $\mathcal{J}$ & Set of renewable electricity generation technology types. \\
    $\mathscr{S}$ & Set of electricity storage technology types. \\
    $\mathcal{V}_u$ & Set of available technology versions for each technology type $u \in \mathcal{J} \cup \mathscr{S}$. \\
    $\mathcal{T}_{u,[t)}$ & Set of periods $t'$ during which technology type $u \in \mathcal{J} \cup \mathscr{S}$, installed in period $t \in \mathcal{T}$, remains operational, i.e., $t \leq t' < t + \tau_{(u,t)}$. \\
    $\mathcal{T}'_{u,[t')}$ & Set of possible installation periods $t$ of technology type $u \in \mathcal{J} \cup \mathscr{S}$ that is operational in period $t' \in \mathcal{T}$, i.e., $t \leq t' < t + \tau_{(u,t)}$. \\
    \hline
    \end{tabularx}
\end{table}

\vspace{-0.15cm}

\begin{table}[H]
    \centering
    \caption{List of decision variables.}
    \label{tab:Variables}
    \small
    \begin{tabularx}{\textwidth}{|c|X|}
    \hline
    \textbf{Variable} & \textbf{Description} \\
    \hline
    $v^+_{(u,v),t,n}$ & Number of installations of technology type $u \in \mathcal{J} \cup \mathscr{S}$, version $v \in \mathcal{V}_u$, in period $t \in \mathcal{T}_n$ at node $n \in \mathcal{N}$. \\
    $v_{(u,v),t,t',n}$ & Number of technology type $u \in \mathcal{J} \cup \mathscr{S}$, version $v \in \mathcal{V}_u$, that was installed in period $t \in \mathcal{T}'_{u,[t')}$ and remains operational during period $t' \in \mathcal{T}_n$ at node $n \in \mathcal{N}$. \\
    $v^-_{(u,v),t,t',n}$ & Number of technology type $u \in \mathcal{J} \cup \mathscr{S}$, version $v \in \mathcal{V}_u$, that was installed in period $t \in \mathcal{T}'_{u,[t')}$ and are salvaged in period $t' \in \mathcal{T}_n$ at node $n \in \mathcal{N}$. \\
    $c_{(q,t),n}$ & Electricity stored at the end of sub-period $q \in \mathcal{Q}$ in period $t \in \mathcal{T}_n$ at node $n \in \mathcal{N}$, carried forward to the next sub-period. \\
    $z_{(q,t),n}^{+}$ & Electricity charged into storage during sub-period $q \in \mathcal{Q}$ in period $t \in \mathcal{T}_n$ at node $n \in \mathcal{N}$. \\
    $z_{(q,t),n}^{-}$ & Electricity discharged from storage during sub-period $q \in \mathcal{Q}$ in period $t \in \mathcal{T}_n$ at node $n \in \mathcal{N}$. \\
    $g_{(q,t),n}$ & Electricity purchased from the grid in sub-period $q \in \mathcal{Q}$ in period $t \in \mathcal{T}_n$ at node $n \in \mathcal{N}$. \\
    \hline
    \end{tabularx}
\end{table}

\begin{table}[H]
    \centering
    \caption{List of parameters.}
    \label{tab:Parameters}
    \small
    \begin{tabularx}{\textwidth}{|c|X|}
    \hline
    \textbf{Parameter} & \textbf{Description} \\
    \hline
    $\pi_n$ & Probability of node $n \in \mathcal{N}$ in the scenario tree. \\
    $\mu_{(n,t)}$ & Ancestor node of $n \in \mathcal{N}$ in period $t \in \mathcal{T}$. \\
    $\kappa_{(q,t),n}$ & Sub-period, period, and node immediately preceding sub-period $q \in \mathcal{Q}$ of period $t \in \mathcal{T}_n$ at node $n \in \mathcal{N}$. \\
    $\beta_{\textnormal{nominal}}$ & Nominal discount rate applied to the periods. \\
    $\zeta$ & Inflation rate applied to the periods. \\
    $\beta_{\textnormal{real}}$ & Real discount rate, calculated as $\beta_{\textnormal{real}} = \frac{1 + \beta_{\textnormal{nominal}}}{1 + \zeta} - 1$. \\
    $\beta$ & Discount factor, calculated as $\beta = \frac{1}{1 + \beta_{\textnormal{real}}}$. \\
    $\alpha^{+}_{(u,v),t,n}$ & Installation cost of technology type $u \in \mathcal{J} \cup \mathscr{S}$, version $v \in \mathcal{V}_u$, in period $t \in \mathcal{T}_n$ at node $n \in \mathcal{N}$. \\
    $\alpha_{(u,v),t,t',n}$ & Operation and maintenance cost of technology type $u \in \mathcal{J} \cup \mathscr{S}$, version $v \in \mathcal{V}_u$, installed in period $t \in \mathcal{T}'_{u,[t')}$ and operated in period $t' \in \mathcal{T}_n$ at node $n \in \mathcal{N}$. \\
    $\alpha^-_{(u,v),t,t',n}$ & Salvage value of technology $u \in \mathcal{J} \cup \mathscr{S}$, version $v \in \mathcal{V}_u$, installed in period $t \in \mathcal{T}'_{u,[t')}$ and salvaged in period $t' \in \mathcal{T}_n$ at node $n \in \mathcal{N}$. \\
    $\lambda_{t}$ & Cost of purchasing grid electricity in period $t \in \mathcal{T}$. \\
    $\delta_{(q,t)}$ & Electricity demand in sub-period $q \in \mathcal{Q}$ of period $t \in \mathcal{T}$. \\
    $\Gamma_{(j,v),q,t,t',n}$ & Renewable electricity generation from technology type $j \in \mathcal{J}$, version $v \in \mathcal{V}_j$, in sub-period $q \in \mathcal{Q}$ of period $t' \in \mathcal{T}_n$ at node $n \in \mathcal{N}$, for installations made in period $t \in \mathcal{T}'_{j,[t')}$. \\
    $\upsilon_{(s,v),t,t',n}$ & Storage capacity from technology type $s \in \mathscr{S}$, version $v \in \mathcal{V}_s$, during period $t' \in \mathcal{T}_n$ at node $n \in \mathcal{N}$, for installations made in period $t \in \mathcal{T}'_{s,[t')}$. \\
    $\eta^{+}, \eta^{-}$ & Charging and discharging efficiencies of storage technologies, respectively. \\
    $\rho_{(u,v),t}$ & Spatial requirement of technology type $u \in \mathcal{J} \cup \mathscr{S}$, version $v \in \mathcal{V}_u$, in period $t \in \mathcal{T}$. \\
    $\tau_{(u,t)}$ & Economic lifetime of technology type $u \in \mathcal{J} \cup \mathscr{S}$ installed in period $t \in \mathcal{T}$. \\
    $\iota_{(u,v)}$ & Number of technology type $u \in \mathcal{J} \cup \mathscr{S}$, version $v \in \mathcal{V}_u$ existing at the start of the planning horizon. \\
    $\epsilon_{t}$ & Emission factor of grid electricity in period $t \in \mathcal{T}$. \\
    $\theta_{t}$ & Maximum permitted emissions in period $t \in \mathcal{T}$. \\
    $\phi_t$ & Available budget for installations in period $t \in \mathcal{T}$. \\
    $\psi_{t}$ & Cumulative maximum installation area available up to period $t \in \mathcal{T}$. \\
    \hline
    \end{tabularx}
\end{table}

\newcommand{\smallplus}{\mathrel{\scalebox{0.8}[0.8]{+}}}
\newcommand{\smallminus}{\mathrel{\scalebox{0.9}[0.9]{-}}}
\newcommand{\smalleq}{\mathrel{\scalebox{0.8}[0.8]{=}}}
\newcommand{\smallge}{\mathrel{\scalebox{0.9}[0.9]{$\ge$}}}
\newcommand{\smallle}{\mathrel{\scalebox{0.8}[0.8]{$\le$}}}
\newcommand{\ssetminussobj}{\mathbin{\tikz[baseline={(0,-0.14ex)}] \draw[line width=0.27pt] (0,0.18) -- (0.06,-0.03);}}
\newcommand{\ssetminusseq}{\mathbin{\tikz[baseline={(0,-0.14ex)}] \draw[line width=0.27pt] (0,0.24) -- (0.07,-0.03);}}
\newcommand{\zeroobj}{\tikz[baseline={(0,-0.14ex)}] \node[inner sep=0pt, anchor=base, scale=0.9] {\scriptsize\{0\}};}
\newcommand{\zeroeq}{\tikz[baseline={(0,-0.14ex)}] \node[inner sep=0pt, anchor=base, scale=1.18] {\scriptsize\{0\}};}

We propose the following multi-stage stochastic programming model, which is formulated as a mixed-integer linear program:

\begin{subequations}\label{eq:stochasticModel}
\begin{align}
    & \text{minimize} \ \sum_{n \in \mathcal{N}\ssetminussobj \zeroobj} \sum_{u \in \mathcal{J} \cup \mathscr{S}} \sum_{v \in \mathcal{V}_u} \sum_{t' \in \mathcal{T}_n} \sum_{t \in \mathcal{T}'_{u,[t')}}  \pi_n \beta^{(t'-1)} \left( \alpha_{(u,v),t,t',n} v_{(u,v),t,t',n} - \alpha^-_{(u,v),t,t',n} v^-_{(u,v),t,t',n} \right) \nonumber\\
    &  \hspace{11.95em} + \sum_{n \in \mathcal{N}\ssetminussobj \zeroobj} \sum_{t \in \mathcal{T}_n} \pi_n \beta^{(t-1)} \left( \sum_{u \in \mathcal{J} \cup \mathscr{S}} \sum_{v \in \mathcal{V}_u} \alpha^{+}_{(u,v),t,n} v^{+}_{(u,v),t,n} + \sum_{q \in \mathcal{Q}} \lambda_{t} g_{(q,t),n} \right) \label{eq:objFunc} 
\end{align}
\vspace{-0.038\textwidth}
\begin{align}
    & \text{subject to } \notag \\
    & g_{(q,t'),n} \smallminus z_{(q,t'),n}^{+} \smallplus z_{(q,t'),n}^{-} \smallplus \sum_{j \in \mathcal{J}} \sum_{v \in \mathcal{V}_j} \sum_{t \in \mathcal{T}'_{j,[t')}} \Gamma_{(j,v),q,t,t',n} v_{(j,v),t,t',n} \smallge \mathit{\delta}_{(q,t')}
    & n \in \mathcal{N} \ssetminusseq \zeroeq, t' \in \mathcal{T}_n, q \in \mathcal{Q} \label{eq:demandConstr}
\end{align}
\vspace{-0.0655\textwidth}
\begin{align}
    & v_{(u,v),t,t',n} = v^+_{(u,v),t,\mu_{(n,t)}} \smallminus \sum_{t'' = t}^{t'} v_{(u,v),t,t'',\mu_{(n,t'')}}^- \quad \quad \hspace{0.5em}
    & n \in \mathcal{N}, t' \in \mathcal{T}_n, u \in \mathcal{J} \cup \mathscr{S}, t \in \mathcal{T}'_{u,[t')}, v \in \mathcal{V}_u \label{eq:techbalanceConstr} \\
    & c_{(q,t'),n} = c_{{\kappa_{(q,t'),n}}} \smallplus \eta^{+} z_{(q,t'),n}^{+} \smallminus \frac{z_{(q,t'),n}^{-}}{\eta^{-}} \quad \quad \hspace{0.5em}
    & n \in \mathcal{N}, t' \in \mathcal{T}_n, q \in \mathcal{Q} \label{eq:storagebalanceConstr} \\
    & c_{(q,t'),n} \leq \sum_{s \in \mathscr{S}} \sum_{v \in \mathcal{V}_s} \sum_{t \in \mathcal{T}'_{s,[t')}} \upsilon_{(s,v),t,t',n} \ v_{(s,v),t,t',n} & n \in \mathcal{N}, t' \in \mathcal{T}_n, q \in \mathcal{Q} \label{eq:batterycapacityConstr} \\
    & \sum_{q \in \mathcal{Q}} \epsilon_{t} \ g_{(q,t),n} \leq \theta_{t}
    & n \in \mathcal{N}, t \in \mathcal{T}_n \label{eq:emissionConstr} \\
    & \sum_{u \in \mathcal{J} \cup \mathscr{S}} \sum_{v \in \mathcal{V}_u} \alpha^+_{(u,v),t,n} \ v^+_{(u,v),t,n} \leq \phi_{t}
    & n \in \mathcal{N}, t \in \mathcal{T}_n \label{eq:budgetConstr} \\
    & \sum_{u \in \mathcal{J} \cup \mathscr{S}} \sum_{v \in \mathcal{V}_u} \sum_{t \in \mathcal{T}'_{u,[t')}} \rho_{(u,v),t} \ v_{(u,v),t,t',n} \leq \psi_{t'}
    & n \in \mathcal{N}, t' \in \mathcal{T}_n \label{eq:spatialConstr} \\[-0.45em]
    & v^+_{(u,v),0,0} = \iota_{u,v}
    & u \in \mathcal{J} \cup \mathscr{S}, v \in \mathcal{V}_u \label{eq:initializationConstr} \\
    & v_{(j,v),t,n}^+ \leq M_{(j,v),t,n}, \quad v_{(j,v),t,n}^+ \in \mathbb{Z}_+
    & n \in \mathcal{N}, t \in \mathcal{T}_n, j \in \mathcal{J}, v \in \mathcal{V}_j \label{eq:VarDomain1}\\
    & v_{(j,v),t,t',n}^- \leq M_{(j,v),t,n}, \quad v_{(j,v),t,t',n}^- \in \mathbb{Z}_+
    & n \in \mathcal{N}, t' \in \mathcal{T}_n, j \in \mathcal{J}, t \in \mathcal{T}'_{j,[t')}, v \in \mathcal{V}_j \label{eq:VarDomain2}\\
    & v_{(s,v),t,n}^+ \in \mathbb{R}_+
    & n \in \mathcal{N}, t \in \mathcal{T}_n, s \in \mathscr{S}, v \in \mathcal{V}_s \label{eq:VarDomain3}\\
    & v_{(s,v),t,t',n}^- \in \mathbb{R}_+ 
    & n \in \mathcal{N}, t' \in \mathcal{T}_n, s \in \mathscr{S}, t \in \mathcal{T}'_{s,[t')}, v \in \mathcal{V}_s \label{eq:VarDomain4}\\
    & v_{(u,v),t,t',n} \in \mathbb{R}_+ 
    & n \in \mathcal{N}, t' \in \mathcal{T}_n, u \in \mathcal{J} \cup \mathscr{S}, t \in \mathcal{T}'_{u,[t')}, v \in \mathcal{V}_u  \label{eq:VarDomain5}\\
    & c_{(q,t),n}, \ g_{(q,t),n}, \ z_{(q,t),n}^{-}, \ z_{(q,t),n}^{+} \in \mathbb{R}_+ 
    & n \in \mathcal{N}, t \in \mathcal{T}_n, q \in \mathcal{Q}.\label{eq:VarDomain6}
\end{align}
\end{subequations}

The objective function~\eqref{eq:objFunc} minimizes the expected total cost over the planning horizon. This cost comprises installation expenditures, salvage values, operation and maintenance (O\&M) expenses, and grid electricity purchases. The stochasticity of the problem is explicitly incorporated through node probabilities.

Constraint set~\eqref{eq:demandConstr} ensures that electricity demand in each sub-period is met through a combination of on-site generation, grid purchases, and storage operations. Technology balance constraints~\eqref{eq:techbalanceConstr} capture operational continuity and aging dynamics of technologies, thereby ensuring intertemporal consistency. Storage-related constraints include (i) the balance equations~\eqref{eq:storagebalanceConstr}, which incorporate charging and discharging efficiencies, and (ii) the capacity bounds~\eqref{eq:batterycapacityConstr}, which restrict stored electricity to remain within the operational capacity limits of storage technologies. Emissions are restricted by~\eqref{eq:emissionConstr} to remain below predetermined thresholds in each period. Financial and spatial feasibility are imposed respectively by budget constraints~\eqref{eq:budgetConstr}, which cap annual installations according to allocated budgets, and spatial constraints~\eqref{eq:spatialConstr}, which limit cumulative installed footprints. Integration of pre-existing technologies at the start of the horizon is governed by~\eqref{eq:initializationConstr}.

Non-negative integer decision variables defined in~\eqref{eq:VarDomain1}–\eqref{eq:VarDomain2} govern installation and salvage decisions for generation technologies. Logical upper bounds, specified in~\eqref{eq:UpperBound}, are imposed to tighten the feasible region and reduce redundancy. These bounds correspond to the maximum number of installations that would be required if the entire residual demand over the horizon were met by a single technology version installed within a single period.

\begin{ceqn}
\begin{equation}
    M_{(j,v),t,n} = \left\lceil \max_{\substack{n' \in \mathcal{N},\, t' \in \mathcal{T}_{j,[t)} \cap \mathcal{T}_{n'},\, q \in \mathcal{Q}}} \left\{ \frac{\delta_{(q,t')}}{\Gamma_{(j,v),q,t,t',n'}} \right\} \right\rceil
    \quad 
    n \in \mathcal{N}, t \in \mathcal{T}_n , j \in \mathcal{J}  ,
    v \in \mathcal{V}_j  .
    \label{eq:UpperBound}
\end{equation}
\end{ceqn}

Continuous decision variables defined in~\eqref{eq:VarDomain3}–\eqref{eq:VarDomain5} capture installation, retirement, and operation of storage technologies \(s \in \mathscr{S}\). While operational decision variables for generation technologies \(u \in \mathcal{J}\), are modeled as continuous decision variables~\eqref{eq:VarDomain5}, their integrality is ensured through the technology balance constraints formulated in~\eqref{eq:techbalanceConstr}. Continuous decision variables defined in~\eqref{eq:VarDomain6} denote grid electricity purchases \(\bigl(g_{(q,t),n}\bigr)\), stored electricity amount \(\bigl(c_{(q,t),n}\bigr)\), and charging/discharging amounts \(\bigl(z_{(q,t),n}^{-}, z_{(q,t),n}^{+}\bigr)\) in each sub-period.


\section{Data Collection and Scenario Tree Generation}
\label{sec:data}

This section details the data utilized in the METU case study and elaborates on the methodology employed for constructing scenario trees. Section~\ref{subsec:Electricity_Technologies} introduces key parameters and technological advancement scenarios for each generation and storage technology, together with the simulation procedures used to obtain hourly electricity generation profiles. Section~\ref{subsec:CaseParameters} describes campus-specific parameters, encompassing operational data required to address short-term mismatches between generation and demand. Section~\ref{subsec:Scenario_Tree_Generation} then describes a systematic procedure for generating scenario trees that incorporate technological advancements into the multi-stage stochastic programming model. The complete dataset and optimization model developed in this research are publicly accessible in~\cite{campusapplicationrep}.

\subsection{Electricity Generation and Storage Technologies}
\label{subsec:Electricity_Technologies}

This subsection provides a comprehensive overview of the parameters associated with the evaluated electricity generation and storage technologies, namely solar photovoltaic, wind turbine, and lithium-ion battery technologies. 
We decide to use these technologies due to the availability of historical datasets and the access to recent cost information from domestic firms. 
The parameters we acquire for each technology include installation costs, O\&M expenditures, projected lifespans, annual performance degradation rates, and spatial land-use requirements. 

Technological advancements are explicitly modeled at stage boundaries, each corresponding to a five-year planning interval comprising discrete decision-making periods, facilitating a thorough evaluation of deployment and operational strategies. Technological advancement scenarios are derived from historical cost and efficiency datasets sourced from relevant literature and technical reports. Drawing inspiration from the exponential improvement trajectories described by Moore’s Law \cite{Moore1998}, logarithmic transformations are applied to these historical data, yielding normalized logarithmic five-year-forward ratios. Subsequently, a clustering model that minimizes the squared Euclidean distances between data points and cluster centroids is implemented to construct technological advancement scenarios. The probability associated with each scenario cluster is computed based on the proportion of historical data points assigned to it. Finally, exponential transformations of cluster centroids yield five-year technological advancement multipliers in cost and efficiency, thus establishing a scenario-based representation of technological uncertainties throughout the planning horizon. {We note that a similar approach is used in \cite{karimi2024} in the context of clean bus fleet transition problem for generating technological advancement scenarios.}

\subsubsection{Solar Photovoltaic Technology}
\label{subsubsec:tech_solar}

Technological progress in solar photovoltaic systems is evaluated using historical cost and efficiency time-series data from \cite{tracking_sun_2019}. Photovoltaic efficiency is defined as the proportion of incident solar energy converted into electricity. Historical data, shown in Figure~\ref{fig:solar_scatter_plot_comparison} (left), indicates a pronounced downward trend in inflation-adjusted installation costs accompanied by steady improvements in module efficiency.

In the advancement scenarios employing a two-cluster configuration (see Figure~\ref{fig:solar_scatter_plot_comparison}, right), the first cluster exhibits five-year logarithmic improvement rates of \(0.0542\) for efficiency and \(0.1552\) for cost, whereas the second cluster demonstrates higher rates of \(0.1247\) and \(0.5510\), respectively. Exponentiation of these rates yields multipliers of \(1.0557\) for efficiency and \(0.8562\) for cost in the first cluster, and \(1.1328\) for efficiency and \(0.5510\) for cost in the second. Scenario probabilities, derived from cluster populations, are \(4/12\) and \(8/12\). These clusters delineate two distinct technological pathways: a conservative trajectory with gradual efficiency improvements and modest cost reductions, and an accelerated trajectory characterized by significant efficiency gains and pronounced cost declines over successive five-year intervals.

To identify appropriate technology versions for on-site deployment at the METU campus, six photovoltaic configurations with varying installed capacities and cost profiles are assessed. Data are obtained from a domestic engineering, procurement, and construction (EPC) firm. Nominal capacities range from 0.018 to 12 \si{\mega\watt}, with the two smallest systems suitable for rooftop installation on campus buildings. The selected panels comprise 144 bifacial cells, achieve a module efficiency of 21.17\%, and have a design life of 25 years. Technical documentation specifies a linear degradation rate of 0.5\% per year, consistent with the established literature~\cite{Jordan2011}. Historical O\&M cost data from~~\cite{LBL2023} are fitted with an exponential function to project annual O\&M expenses per \si{\kilo\watt} of installed capacity, resulting in an annual multiplier of 0.918.

\begin{figure}[H]
\centering
\begin{subfigure}{0.60\textwidth}
\centering
\begin{tikzpicture}
    \pgfplotsset{
        width=\textwidth,
        height=6.2cm,
        xmin=2002, xmax=2018,
        y axis style/.style={yticklabel style=#1, ylabel style=#1, y axis line style=#1, ytick style=#1}
    }
    \pgfplotsset{every axis y label/.append style={rotate=0,yshift=0cm}}
    \begin{axis}[
      xticklabel style={/pgf/number format/1000 sep=},
      xtick pos=bottom,
      xtick={2002,2006,2010,2014,2018},
      axis y line*=left,
      y axis style=black!75!black,
      legend style={at={(0.8,0.8)},anchor=south,legend columns=1, nodes={scale=0.5, transform shape}, font=\LARGE},
      ymin=3, ymax=15,
      ytick={3,6,9,12,15},
      ylabel style={xshift=-0.02cm, yshift=-0.26cm, font=\fontsize{9pt}{10pt}\selectfont},
      xlabel style={yshift=0.1cm, font=\fontsize{9pt}{10pt}\selectfont},
      xticklabel style={yshift=-2pt},
      tick label style={font=\footnotesize},
      xlabel=Year,
      ylabel=Installation Cost (2024 USD/\unit{\W})
    ]
    \addplot[mark=*,black] 
      coordinates{
(2002, 13.9424443729401)
(2003, 13.5238860865593)
(2004, 12.3659295320511)
(2005, 12.0786227212906)
(2006, 12.3256687812805)
(2007, 12.213283805275)
(2008, 11.7068575906754)
(2009, 11.3113430204391)
(2010, 9.35903756890297)
(2011, 8.03157230482102)
(2012, 6.88202373743057)
(2013, 5.63355725450516)
(2014, 5.00414950094223)
(2015, 4.75094442675114)
(2016, 4.4341953799963)
(2017, 4.02545742025375)
(2018, 3.73654250013828)
    };
\addlegendentry{Cost Data}
\addplot[mark=o,black] 
      coordinates{
    (2002,0.12692)
    (2003,0.12987)
    (2004,0.13452)
    (2005,0.13452)
    (2006,0.13492)
    (2007,0.13492)
    (2008,0.14113)
    (2009,0.13805)
    (2010,0.14094)
    (2011,0.14417)
    (2012,0.14844)
    (2013,0.15409)
    (2014,0.15741)
    (2015,0.16352)
    (2016,0.16460)
    (2017,0.17378)
    (2018,0.18405)
    };
\addlegendentry{Efficiency Data}
\end{axis}
\pgfplotsset{every axis y label/.append style={rotate=180,yshift=0cm}}
    \begin{axis}[
      axis y line*=right,
      axis x line=none,
      ymin=0.12, ymax=0.20,
      ylabel=Module Efficiency (\%),
      y axis style=black!75!black,
      ylabel style={xshift=-0.02cm, yshift=-0.12cm, font=\fontsize{9pt}{10pt}\selectfont},
      tick label style={font=\footnotesize},
      yticklabel style={/pgf/number format/.cd, fixed, precision=2, fixed zerofill},
    ]
\addplot[mark=o,black] 
      coordinates{
    (2002,0.12692)
    (2003,0.12987)
    (2004,0.13452)
    (2005,0.13452)
    (2006,0.13492)
    (2007,0.13492)
    (2008,0.14113)
    (2009,0.13805)
    (2010,0.14094)
    (2011,0.14417)
    (2012,0.14844)
    (2013,0.15409)
    (2014,0.15741)
    (2015,0.16352)
    (2016,0.16460)
    (2017,0.17378)
    (2018,0.18405)
    };
\end{axis}
\end{tikzpicture}
\end{subfigure}
\hspace{0.05\textwidth}
\begin{subfigure}{0.28\textwidth}
\centering
\begin{tikzpicture}
    \begin{axis}[
        width=\textwidth,
        xtick pos=bottom,
        ytick pos=left,
        height=6.2cm,
        xlabel={Efficiency Improvement Rate},
        ylabel={Cost Improvement Rate},
        grid=major,
        xmin=0, xmax=0.20,
        ymin=0.05, ymax=0.85,
        xtick={0,0.10,0.20},
        ytick={0.05,0.25,0.45,0.65,0.85},
        yticklabel style={/pgf/number format/.cd, fixed, precision=2, fixed zerofill},
        ylabel style={xshift=0.01cm, yshift=-0.15cm, font=\fontsize{9pt}{10pt}\selectfont},
        xlabel style={yshift=0.1cm, font=\fontsize{9pt}{10pt}\selectfont},
        xticklabel style={yshift=-2pt},
        tick label style={font=\footnotesize},
        xticklabel style={
          /pgf/number format/.cd,
          fixed,
          precision=2,
          fixed zerofill
        },
        minor y tick num=2,
        scatter/classes={
            1={mark=triangle,blue},
            2={mark=diamond,green},
            3={mark=square,orange},
            representative1={mark=triangle*,draw=blue,fill=blue,scale=1.5},
            representative2={mark=diamond*,draw=green,fill=green,scale=1.5},
            representative3={mark=square*,draw=orange,fill=orange,scale=1.5}
        }
    ]
    \addplot[scatter,only marks,scatter src=explicit symbolic] 
    table[meta=label] {
        x y label
    0.0611055 0.1324135 1
    0.0831396 0.1442827 1
    0.0259040 0.0891390 1
    0.0466401 0.2550947 1
    };
    \addplot[scatter,only marks,scatter src=explicit symbolic] 
    table[meta=label] {
        x y label
    0.0663187 0.4283037 3
    0.0954772 0.5736114 3
    0.0878496 0.7314337 3
    0.1312421 0.8155386 3
    0.1486163 0.6779990 3
    0.1324902 0.5940341 3
    0.1576289 0.5362742 3
    0.1776783 0.4105804 3
    };
    \addplot[scatter,only marks,scatter src=explicit symbolic] 
    table[meta=label] {
        x y label 
        0.054197762	0.155231397 representative1
        0.1246674 0.595971548 representative3
    };
    \node[below] at (axis cs:0.054197762, 0.155231397) {\scriptsize \textbf{S}};
    \node[below] at (axis cs:0.1246674, 0.595971548) {\scriptsize \textbf{F}};
    \end{axis}
\end{tikzpicture}
\end{subfigure}
\vspace{-6pt}
\caption{Historical trends in solar technology installation cost and module efficiency (left), and the resulting two-cluster technological improvement scenarios derived from five-year improvement rates (right).}
\label{fig:solar_scatter_plot_comparison}
\end{figure}

Hourly electricity generation profiles are simulated for each technology version using a typical meteorological year (TMY) dataset for the METU campus site (39.84°N, 32.81°E), obtained from the National Solar Radiation Database (NSRDB) \cite{NSRDB2018}. Simulations employ the open-source NREL SAM \cite{sam2025}, combined with site-specific parameters and manufacturer specifications for panels and inverters. The resulting hourly generation data are incorporated into the optimization model, providing a high-resolution representation of the short-term variability inherent in solar photovoltaic output.

\subsubsection{Wind Turbine Technology} \label{subsubsec:tech_wind}

\begin{figure}[H]
\centering
\begin{subfigure}{0.45\textwidth}
\centering
\begin{tikzpicture}
  \pgfplotsset{
    width=\textwidth,
    height=6cm,
    xmin=2008, xmax=2040,
  }
  \begin{axis}[
      xticklabel style={/pgf/number format/1000 sep={}},
      xtick pos=bottom,
      axis y line*=left,
      axis line style={color=black!75!black},
      xlabel=Year,  
      ylabel={Cost (2024 USD/\unit{\kW})},
      ymin=0, ymax=2500,
      ytick={0,500,1000,1500,2000,2500},
      legend style={at={(0.705,0.7212)}, anchor=south, legend columns=1, nodes={scale=0.56, transform shape}, font=\large},
      legend cell align={left},
      ylabel style={xshift=-0.1cm, yshift=-0.12cm, font=\fontsize{9pt}{10pt}\selectfont},
      xlabel style={yshift=0.1cm, font=\fontsize{9pt}{10pt}\selectfont},
      xticklabel style={yshift=-2pt},
      tick label style={font=\footnotesize},
  ]
    \addplot[mark=*,black] coordinates {
(2008, 2431.120679)
(2009, 2474.007243)
(2010, 2337.495233)
(2011, 2311.367451)
(2012, 2379.951533)
(2013, 2196.583773)
(2014, 2146.296668)
(2015, 1978.751866)
(2016, 1967.177765)
(2017, 1963.520521)
(2018, 1872.455145)
(2019, 1778.173546)
(2020, 1609.563841)
(2021, 1525.479499)
(2022, 1370.595215)
    };
\addlegendentry{Historical Installation Cost}
\addplot[mark=o,color=red!60!white] 
        coordinates{
(2023, 1315.621107)
(2024, 1262.851991)
(2025, 1212.199427)
(2027, 1116.907781)
(2030, 987.8299272)
(2035, 804.9896274)
(2040, 655.9917678)
    };
\addlegendentry{Projected Installation Cost}
\addlegendimage{mark=diamond*,black}\addlegendentry{Historical O\&M Cost}
\addlegendimage{mark=diamond,color=red!80!white}\addlegendentry{Projected O\&M Cost}
\end{axis}
\pgfplotsset{every axis y label/.append style={rotate=180}}
  \begin{axis}[
      axis y line*=right,
      axis x line=none,
      ymin=0, ymax=125,
      ylabel={O\&M Cost (2024 USD/\unit{\kW}-year)},
      ytick= {0,25,50,75,100,125},
      axis line style={color=black!75!black},
      ylabel style={xshift=0.0cm, yshift=-0.1cm, font=\fontsize{9pt}{10pt}\selectfont},
      tick label style={font=\small},
  ]
    \addplot[mark=diamond*,black] coordinates {
      (2008,77.16877056)
      (2009,75.17617294)
      (2010,67.6970404)
      (2011,71.4602209)
      (2012,48.53819591)
      (2013,68.13888512)
      (2014,37.95703353)
      (2015,27.64036322)
      (2016,33.66907734)
      (2017,26.41085465)
      (2018,22.88243125)
      (2019,19.75564402)
      (2020,21.5135155)
      (2021,17.38309339)
      (2022,17.04699591)
    };
    \addplot[mark=diamond,color=red!80!white] 
          coordinates{
            (2023, 15.3040173)
            (2024, 13.73925042)
            (2025, 12.3344739)
            (2027, 9.941131361)
            (2030, 7.192975581)
            (2035, 4.194657845)
            (2040, 2.446157955)
          };
\end{axis}
\end{tikzpicture}
\caption{Historical and projected trends in wind turbine installation and annual O\&M costs.}
\label{fig:wind_cost_and_OM}
\end{subfigure}
\hfill
\begin{subfigure}{0.45\textwidth}
\centering
\begin{tikzpicture}
  \pgfplotsset{width=\textwidth, height=6cm}
  \begin{axis}[
      axis lines=box,
      xtick pos=bottom,
      ytick pos=left,
      xlabel=Commissioning Year,
      ylabel= {\color{white}dum}Wind Turbine Efficiency (\%){\color{white}dum},
      grid=major,
      xmin=1995, xmax=2025,
      ymin=0, ymax=0.6,
      xtick={1995,2005,2015,2025},
      xticklabel style={/pgf/number format/1000 sep={}},
      ytick={0,0.12,...,0.60},
      ylabel style={xshift=-0.01cm, yshift=-0.1cm, font=\fontsize{9pt}{10pt}\selectfont},
      xlabel style={yshift=0.1cm, font=\fontsize{9pt}{10pt}\selectfont},
      xticklabel style={yshift=-2pt},
      yticklabel style={/pgf/number format/.cd, fixed, precision=2, fixed zerofill},
      tick label style={font=\footnotesize},
      scatter/classes={
      1={mark=triangle,blue},
      representative1={mark=triangle*,draw=blue,fill=blue,scale=1.5}
      }
  ]
    \addplot[scatter,only marks,scatter src=explicit symbolic] 
      table[meta=label] {
        x y label
        2022 0.181 1
        2004 0.346 1
        2016 0.169 1
        1998 0.286 1
        1998 0.286 1
        2005 0.255 1
        2003 0.289 1
        2007 0.173 1
        2006 0.385 1
        2004 0.304 1
        2000 0.156 1
        2000 0.156 1
        2016 0.171 1
        2020 0.364 1
        2017 0.244 1
        2017 0.376 1
        2006 0.488 1
        2005 0.296 1
        2002 0.296 1
        2002 0.192 1
        2000 0.213 1
        2000 0.241 1
        2017 0.195 1
        2009 0.401 1
        2006 0.169 1
        2023 0.253 1
        2009 0.251 1
        2005 0.317 1
        2003 0.183 1
        2009 0.259 1
        2016 0.265 1
        2010 0.189 1
        2008 0.210 1
        2016 0.184 1
        2009 0.279 1
        2004 0.209 1
        2005 0.265 1
        2009 0.283 1
        2012 0.257 1
        2020 0.363 1
        2016 0.122 1
        2004 0.134 1
        2016 0.208 1
        2012 0.231 1
        2020 0.424 1
        2005 0.200 1
        2016 0.311 1
        2015 0.205 1
        2003 0.220 1
        2016 0.254 1
        2010 0.239 1
        2020 0.393 1
        2017 0.385 1
        2005 0.297 1
        2004 0.209 1
        2002 0.269 1
        2017 0.261 1
        2010 0.371 1
        2010 0.234 1
        2005 0.292 1
        2005 0.261 1
        2005 0.234 1
        2008 0.228 1
        2004 0.394 1
        2002 0.171 1
        2016 0.232 1
        2005 0.424 1
        2012 0.263 1
        2009 0.363 1
        2010 0.194 1
        2017 0.420 1
        2013 0.420 1
        2020 0.404 1
        2011 0.278 1
        2017 0.180 1
        2016 0.224 1
        2005 0.252 1
        2005 0.219 1
        2011 0.223 1
        2011 0.318 1
        2006 0.167 1
        2011 0.238 1
        2009 0.392 1
        2013 0.251 1
        2011 0.250 1
        2013 0.266 1
        2020 0.416 1
        2017 0.247 1
        2016 0.218 1
        2017 0.266 1
        2009 0.189 1
        2006 0.342 1
        2004 0.211 1
        2009 0.260 1
        2011 0.410 1
        2009 0.389 1
        2004 0.284 1
        2013 0.181 1
        2013 0.181 1
        2008 0.271 1
        2011 0.201 1
        2010 0.278 1
        2012 0.283 1
        2012 0.283 1
        2013 0.333 1
        2015 0.266 1
        2020 0.251 1
        2019 0.283 1
        2011 0.219 1
        2013 0.254 1
        2020 0.332 1
        2017 0.278 1
        2002 0.205 1
        2011 0.425 1
        2013 0.227 1
        2014 0.388 1
        2014 0.268 1
        2005 0.146 1
        2011 0.308 1
        2017 0.246 1
        2012 0.384 1
        2020 0.344 1
        2013 0.199 1
        2017 0.248 1
        2017 0.395 1
        2020 0.220 1
        2018 0.249 1
        2017 0.429 1
        2015 0.349 1
        2009 0.134 1
        2015 0.232 1
        2014 0.280 1
        2013 0.199 1
        2015 0.293 1
        2016 0.281 1
        2017 0.260 1
        2011 0.367 1
        2014 0.164 1
        2018 0.461 1
        2016 0.336 1
        2016 0.333 1
        2002 0.262 1
        2005 0.437 1
        2012 0.319 1
        2019 0.256 1
        2019 0.275 1
        2013 0.227 1
        2013 0.310 1
        2019 0.262 1
        2016 0.286 1
        2013 0.270 1
        2013 0.181 1
        2012 0.389 1
        2009 0.352 1
        2013 0.234 1
        2012 0.274 1
        2013 0.228 1
        2016 0.344 1
        2017 0.327 1
        2020 0.348 1
        2017 0.275 1
        2012 0.353 1
        2019 0.285 1
        2018 0.354 1
        2020 0.353 1
        2015 0.218 1
        2017 0.221 1
        2014 0.209 1
        2017 0.285 1
        2015 0.203 1
        2012 0.197 1
        2020 0.210 1
        2019 0.210 1
        2020 0.230 1
        2019 0.230 1
        2010 0.189 1
        2017 0.220 1
        2017 0.390 1
        2017 0.414 1
        2017 0.368 1
        2023 0.226 1
        2016 0.275 1
        2017 0.315 1
        2013 0.231 1
        2013 0.250 1
        2015 0.297 1
        2017 0.238 1
        2019 0.229 1
        2004 0.215 1
        2017 0.153 1
        2016 0.197 1
        2007 0.172 1
        2009 0.292 1
        2018 0.379 1
        2020 0.393 1
        2020 0.361 1
        2018 0.213 1
        2015 0.202 1
        2014 0.185 1
        2016 0.202 1
      };
  \end{axis}
\end{tikzpicture}
\caption{Efficiencies of wind turbine models by commissioning year.}
\label{fig:wind_efficiency}
\end{subfigure}
\caption{Techno-economical trends for wind turbines.}
\label{fig:wind_combined_plots}
\end{figure}
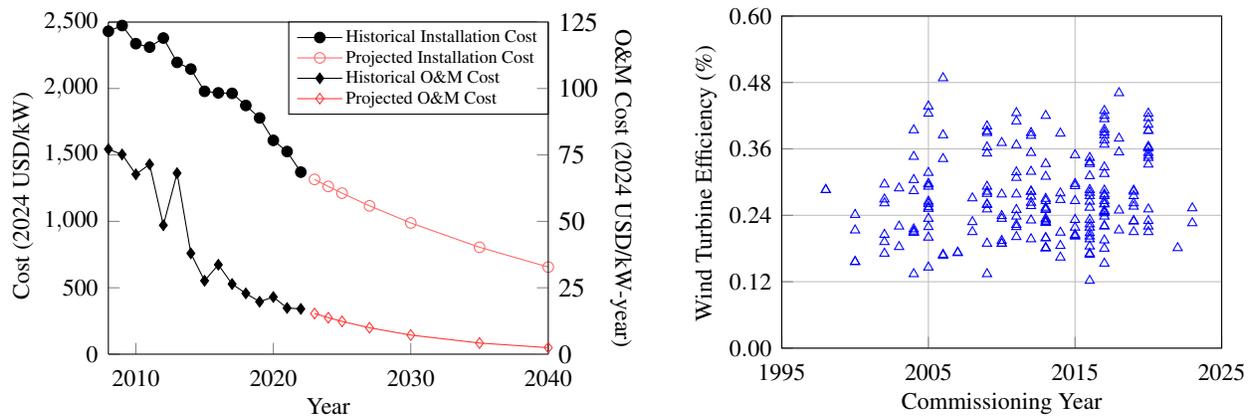

Wind turbine efficiency is the capability of a turbine to convert the kinetic energy of wind into electrical output. According to Betz’s law \cite{Bergey2012}, the theoretical maximum efficiency achievable by any wind turbine is 59.3\%; however, practical efficiencies remain lower due to aerodynamic and mechanical limitations \cite{Rosenberg2014}. Given the scarcity of comprehensive historical efficiency data, we estimate wind turbine efficiencies across different commissioning years using other technical specifications compiled from \cite{thewindpowerDB}. These specifications include key parameters such as rotor diameter, rated power, rated wind speed, cut‑in and cut‑out wind speeds, and commissioning year.
The power curve approximation function available in the NREL SAM \cite{sam2025} is used to derive the efficiency of each turbine. 
Analysis of the calculated efficiencies, plotted against commissioning years (Figure~\ref{fig:wind_efficiency}), revealed no discernible temporal trend. This result led us to assume a constant wind turbine efficiency over time, in line with the literature \cite{Rosenberg2014, Molteno2022Biomimetics}.

Following the efficiency assessment, technological advancement is quantified through the analysis of installation cost trends. Historical global weighted-average installed cost data for onshore wind power plants \cite{irena2023renewable} are presented in Figure~\ref{fig:wind_cost_and_OM} (left) as solid circles, alongside an exponential trend line projecting future expenditures. This decline, analogous to Moore’s Law, provides a deterministic technological advancement trajectory under the assumption of constant turbine efficiency. An exponential fit to cost data from 2009 to 2022 yields an annual cost multiplier of 0.9599, corresponding to a five-year multiplier of 0.8149. Installation costs for wind projects are influenced by turbine specifications, hub heights, and administrative and logistical requirements. For the case study, a representative project scheduled for 2024 was selected. Located approximately 222.5 kilometers from the METU campus, it consists of four Goldwind GW165-6.0 turbines—the sole candidate model considered—providing a total capacity of \SI{24}{\mega\watt} with an expected operational lifetime of 25 years.

Historical O\&M cost data, sourced from \cite{DOE2023} and normalized by installed capacity, are shown in Figure~\ref{fig:wind_cost_and_OM} (solid diamonds). These costs are extrapolated using an exponential trend to forecast future O\&M expenditures. Additionally, turbines are assumed to degrade at a constant annual rate of 1.6\% \cite{Staffell2014}, accounting for performance losses due to mechanical wear and environmental factors. Site area requirements are assumed to comply with national regulations \cite{windplantREGULATION}, which specify that the minimum required area must be at least equal to the turbine’s swept area for single-turbine installations.

%

To generate wind power simulation inputs, representative months and years from the TMY dataset \cite{NSRDB2018} are first identified. For these selected periods, hourly wind speed and direction data at 10 and 50 meters above ground level are obtained from NASA meteorological datasets \cite{NASAReference}. Since direct measurements at the 100-meter hub height are unavailable, wind speeds are extrapolated using the shear formula of \cite{Firtin2011}. A shear coefficient of 0.2215, derived from 42 years of NASA historical averages, yields an estimated annual mean wind speed of 4.90~m/s at 100 meters. This estimate aligns closely with the Global Wind Atlas \cite{globalwindatlas}, supporting the validity of the extrapolation. Finally, the extrapolated wind speeds are integrated with turbine specifications to produce hourly generation data using NREL’s SAM.

\subsubsection{Lithium-ion Battery Technology}
\label{subsubsec:tech_battery}

Lithium-ion battery storage systems are incorporated into this study. Time-series data on cost and energy density trends are extracted from \cite{walter2023} using the WebPlotDigitizer tool \cite{autometris}, ensuring accurate representation of historical and projected technological advancements. Advancement scenarios for battery technology are represented through the clusters proposed by \cite{karimi2024}. To establish a representative battery technology, technical specifications from various lithium-ion battery datasheets are reviewed, with key parameters from \cite{sungrowcatalogue} integrated into the model.

Within the modeling framework, batteries are assumed to reach end-of-life once their capacity declines to 80\% of the initial rated capacity, consistent with industry specifications. Battery systems are assumed to have an operational lifetime of 20 years, corresponding to an average annual degradation rate of 1\% \cite{sungrowcatalogue}. Charging and discharging efficiencies are fixed at 90\% \cite{ZHAO2024110398, Mavromatidis2021}. In the absence of project-specific cost data, installation costs are assumed to be \$350 per \si{\kilo\watt\hour}, following \cite{Cole2025NREL}. A depth-of-discharge limit of 90\% is imposed by adjusting the installation cost to \$388.89 per \si{\kilo\watt\hour} \cite{Mohamed2021, Tejero2024}. Consistent with the 0.9 divisor used to express costs per effective usable capacity, we likewise divide the baseline spatial requirement by 0.9; accordingly, the model adopts \SI[locale=US, group-digits=false]{0.033}{\metre\squared} per \si{\kilo\watt\hour} of installed capacity \cite{reber2023beyond}. Under these assumptions, a single representative battery version is modeled, with deployment decisions formulated as continuous variables representing installed storage capacity in \si{\kilo\watt\hour}.

Table~\ref{tab:VersionsTable} summarizes the key specifications of the candidate technology versions considered for deployment at METU. Reported values include nominal capacities, installation costs, spatial requirements, and annually aggregated electricity generation for the initial year of operation under the meteorological conditions of the METU campus in Ankara, along with the corresponding capacity factors. In addition, the salvage value of each technology is assumed to start at 20\% of its installation cost \cite{Paspatis2022} and to depreciate linearly to zero by the end of its operational lifetime. Salvage values are further adjusted across stage transitions in line with the respective installation cost multipliers.

\begin{table}[H]
\centering
\caption{Candidate electricity generation and storage technology versions for deployment.}
\label{tab:VersionsTable}
\begin{adjustbox}{max width=\textwidth}
\begin{tabular}{|c|c|c|c|c|c|c|}
    \hline
    \textbf{Technology} & \textbf{Capacity} & \textbf{Installation Cost} & \textbf{Spatial Requirement} & \textbf{Annual Electricity Generation}  & \textbf{Capacity Factor} \\ \hline
    Solar V$_{1}$    & \SI{18}{\kilo\watt}     & \$20,047 & \SI{88}{\meter\squared}  & \SI{0.029}{\giga\watt\hour} & \multirow{6}{*}{18.32\%}  \\ \cline{1-5}
    
    Solar V$_{2}$    & \SI{36}{\kilo\watt}     & \$32,570 & \SI{177}{\meter\squared} & \SI{0.059}{\giga\watt\hour} &  \\ \cline{1-5}
    
    Solar V$_{3}$    & \SI{0.6}{\mega\watt}    & \$631,742 & \SI{2837}{\meter\squared} & \SI{0.963}{\giga\watt\hour} &  \\ \cline{1-5}
    
    Solar V$_{4}$    & \SI{1.2}{\mega\watt}   & \$945,275 & \SI{5674}{\meter\squared} & \SI{1.926}{\giga\watt\hour} &  \\ \cline{1-5}
    
    Solar V$_{5}$    & \SI{6}{\mega\watt}   & \$3,190,170 & \SI{28370}{\meter\squared} & \SI{9.628}{\giga\watt\hour} &  \\ \cline{1-5}
    
    Solar V$_{6}$    & \SI{12}{\mega\watt}  & \$5,832,190 & \SI{56739}{\meter\squared} & \SI{19.256}{\giga\watt\hour} &  \\ \hline

    Wind V$_{1}$     & \SI{6}{\mega\watt}   & \$5,430,435 & \SI{86409}{\meter\squared} & \SI{5.466}{\giga\watt\hour} & 10.40\% \\ \hline
    
    Battery V$_{1}$  & \SI{1}{\kilo\watt\hour} & \$388.89 & \SI[locale=US, group-digits=false]{0.033}{\metre\squared} & - & - \\ \hline
\end{tabular}
\end{adjustbox}
\end{table}

\subsection{Case-Specific Parameters}
\label{subsec:CaseParameters}

Hourly electricity consumption data for the METU campus are obtained directly from university records for the years 2023 and 2024. The corresponding monthly consumption profiles are presented in Figure~\ref{fig:metuCostDemand} (left), which reveal a clear seasonal pattern. Owing to the limited availability of historical records, the average hourly consumption across 2023 and 2024 was adopted as the representative load profile for this study. To address discrepancies arising from a one-day misalignment between the two datasets, records from January 1, 2023, and December 30–31, 2024, are excluded. The finalized dataset therefore comprises 8,736 consecutive hourly observations, corresponding to a total annual demand of \SI{34.44}{\giga\watt\hour}.

Electricity purchasing cost data for METU are obtained from the Energy Market Regulatory Authority (EPDK), which governs the Turkish electricity market \cite{epdk}. Specifically, the cost data are drawn from the "Final Tariff Table" under the category "Distribution System Users, Single-Rate, Public and Private Services Sector and Others", which accurately reflects the pricing scheme applicable to public universities in Turkey. METU’s inflation-adjusted purchasing costs are depicted in Figure~\ref{fig:metuCostDemand} (right). Given the absence of a consistent long-term trend—due to volatility driven by fuel price fluctuations, regulatory revisions, and seasonal demand changes, a three-year average was computed and adopted as the representative cost. This value is calculated as \$0.144 per \si{\kilo\watt\hour}.

\begin{figure}[H]
\centering
\begin{subfigure}{0.47\textwidth}
\centering
\begin{tikzpicture}
    \pgfplotsset{width=1.02\textwidth, height=6cm}
    \begin{axis}[
      axis lines=box,
      xtick pos=bottom,
      ytick pos=left,
      ytick={0,1000,2000,3000,4000,5000},
      ymin=0, ymax=5000,
      xlabel={Month},
      ylabel={Electricity Demand (\si{\mega\watt\hour})},
      symbolic x coords={Jan, Feb, Mar, Apr, May, Jun, Jul, Aug, Sep, Oct, Nov, Dec},
      xtick={Jan, Mar, May, Jul, Sep, Nov},
      legend style={at={(0.98,0.05)}, anchor=south east, legend columns=1, nodes={scale=0.6, transform shape}, font=\large},
      y axis line style=black!75!black,
      ylabel style={xshift=0.01cm, yshift=-0.05cm, font=\fontsize{9pt}{10pt}\selectfont},
      xlabel style={yshift=0.1cm, font=\fontsize{9pt}{10pt}\selectfont},
      xticklabel style={yshift=-2pt},   
      tick label style={font=\footnotesize},
      ytick style=black!75!black,
      xtick style=black!75!black,
      yticklabel style={/pgf/number format/.cd, fixed},
      legend reversed,
    ]
    \addplot[mark=*,black,thick] coordinates {
        (Jan,3200)
        (Feb,2617)
        (Mar,2817)
        (Apr,2542)
        (May,2375)
        (Jun,2132)
        (Jul,2058)
        (Aug,2354)
        (Sep,2023)
        (Oct,2538)
        (Nov,2939)
        (Dec,3291)
    };
    \addlegendentry{2023}
    \addplot[mark=triangle*,red,thick] coordinates {
        (Jan,3301)
        (Feb,2807)
        (Mar,3063)
        (Apr,2690)
        (May,2613)
        (Jun,2485)
        (Jul,2579)
        (Aug,2656)
        (Sep,2422)
        (Oct,2950)
        (Nov,3330)
        (Dec,3694)
    };
    \addlegendentry{Average}
    \addplot[mark=square*,blue,thick] coordinates {
        (Jan,3402)
        (Feb,2997)
        (Mar,3308)
        (Apr,2838)
        (May,2850)
        (Jun,2838)
        (Jul,3100)
        (Aug,2958)
        (Sep,2820)
        (Oct,3362)
        (Nov,3721)
        (Dec,4096)
    };
    \addlegendentry{2024}
    \end{axis}
\end{tikzpicture}
\end{subfigure}
\hfill
\begin{subfigure}{0.47\textwidth}
\centering
\begin{tikzpicture}
    \pgfplotsset{width=1.02\textwidth, height=6cm,
    }
    \begin{axis}[
      axis lines=box,
      xtick pos=bottom,
      ytick pos=left,
      xmin={1}, xmax={168},
      ymin=0, ymax=0.25,
      xlabel={Year},
      ylabel={Electricity Tariff (2024 USD/\si{\kilo\watt\hour})},
      ytick={0.00,0.05,0.10,0.15,0.20,0.25},
      xtick={1, 49, 109, 168},
      xticklabels={2011,2015,2020,2025},
      y axis line style=black!75!black,
      ytick style=black!75!black,
      x axis line style=black!75!black,
      xtick style=black!75!black,
      xlabel style={yshift=0.1cm},
      xticklabel style={yshift=-2pt},
      tick label style={font=\footnotesize},
      yticklabel style={/pgf/number format/.cd, fixed, precision=2, fixed zerofill},
      ylabel style={xshift=0.04cm, yshift=-0.05cm, font=\fontsize{9pt}{10pt}\selectfont},
      xlabel style={yshift=0.1cm, font=\fontsize{9pt}{10pt}\selectfont},
    ]
    \addplot[const plot, thick, mark=none, black] coordinates {
    (1, 0.220146450607666)
    (4, 0.232216803628953)
    (7, 0.209169763723574)
    (10, 0.192491019198193)
    (13, 0.188301288471277)
    (16, 0.201360970111016)
    (22, 0.185017413024086)
    (28, 0.209357602677078)
    (31, 0.193964499560792)
    (34, 0.188065879759519)
    (37, 0.163203309152033)
    (40, 0.17556997443787)
    (43, 0.173224777450843)
    (46, 0.182520492685123)
    (49, 0.164295108784557)
    (52, 0.150389578108037)
    (55, 0.145030931407942)
    (58, 0.137859263580522)
    (61, 0.1358802586061)
    (64, 0.143649494363929)
    (67, 0.134378372443343)
    (70, 0.129790330110899)
    (73, 0.104540528257873)
    (76, 0.111070998929758)
    (79, 0.112090694028366)
    (82, 0.104027615088367)
    (85, 0.111617015741949)
    (88, 0.106161558302756)
    (92, 0.0753027508131329)
    (93, 0.0921644314120562)
    (94, 0.118115882888534)
    (97, 0.125899607933186)
    (100, 0.10919969530503)
    (103, 0.133861748701973)
    (106, 0.150508063278717)
    (109, 0.142002005667759)
    (112, 0.121580632989159)
    (115, 0.121633751800522)
    (118, 0.107350143286837)
    (121, 0.124662765372987)
    (124, 0.110151708539312)
    (127, 0.124114542705702)
    (133, 0.16768079847303)
    (134, 0.161177918536871)
    (135, 0.152106029829584)
    (136, 0.150363942048027)
    (138, 0.16688948184601)
    (139, 0.155539834256395)
    (141, 0.195533959178378)
    (142, 0.194566385800344)
    (145, 0.186572192189893)
    (148, 0.153606214913769)
    (151, 0.111194602160213)
    (154, 0.126669651793798)
    (157, 0.115965904780126)
    (163, 0.127709089592367)
    (168, 0.127709089592367)
    };
    \addplot[red, dashed, thick] coordinates {(1,0.144) (168,0.144)};
    \end{axis}
\end{tikzpicture}
\end{subfigure}
\caption{Monthly electricity demand profiles at METU for the years 2023 and 2024 (left), and grid electricity tariffs from 2011 to 2025 with the last 3-year average shown in red (right).}
\label{fig:metuCostDemand}
\end{figure}

The designated site for renewable electricity generation and storage deployment lies within the METU campus, near Eymir Lake (39.84°N, 32.81°E). The area remains largely undeveloped, covered by steppe vegetation, and offers ample space to accommodate the proposed installations without imposing spatial constraints on the model.

\subsection{Scenario Tree Generation}
\label{subsec:Scenario_Tree_Generation}

{A technology tree is a stage-wise probabilistic structure that captures the uncertain evolution of key parameters—namely installation cost and efficiency—for a single electricity technology across the planning horizon. Each node in the tree represents the state of the technology at a particular stage, while the outgoing branches denote mutually exclusive technological advancement scenarios.} The root node corresponds to the initial state, where the baseline installation cost and efficiency values of all candidate technology versions are initialized. For each subsequent node, parameter values are derived by applying technology-specific cost and efficiency multipliers to the values of their parent node. These multipliers, which encapsulate technological progress between stages, are obtained through clustering methods as outlined earlier. The probability associated with each branch reflects the likelihood of a given advancement scenario. By recursively applying this process, comprehensive technology trees are formed, wherein all nodes at the same depth correspond to a specific stage, thus capturing the stochastic dynamics of technological development over time.

Illustrative examples of the technology trees are presented in Figure~\ref{fig:tech_trees}. For clarity, each example depicts a single candidate version of the respective technology. Specifically, the wind power technology tree features a single branch, reflecting the assumption of exponential cost reduction coupled with deterministic and constant efficiency over time. In contrast, the solar power technology tree includes two branches, each corresponding to distinct advancement scenarios identified through clustering. The associated branch probabilities, denoted as \(\theta^{s}_{solar}\) and \(\theta^{f}_{solar}\), reflect slow and fast technological advancement scenarios, respectively. Similarly, the battery technology tree comprises two branches, each corresponding to distinct cost and efficiency improvement scenarios, also derived from clustering. While the illustrative examples adopt two clusters for solar and battery technologies, the model is extendable to incorporate additional clusters, which would increase the number of branches emerging from each node proportionally.

{\small
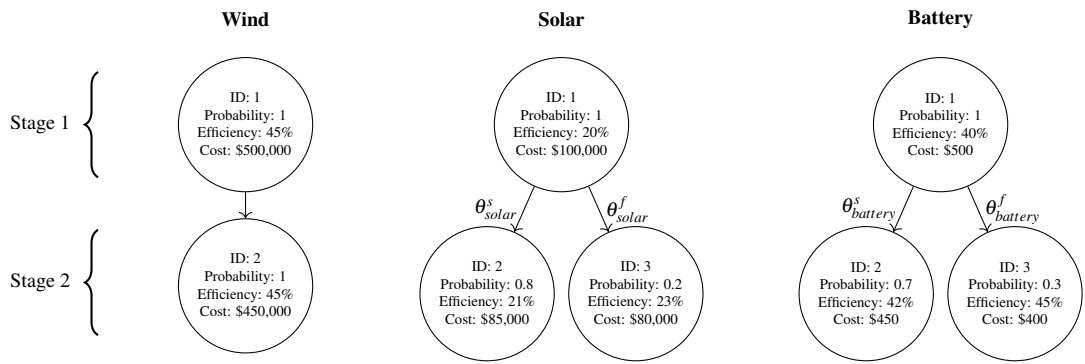
\begin{figure}[H]
\centering
\begin{tikzpicture}[scale=0.7, transform shape, ->,
circleNode/.style={circle, draw, align=center, minimum size=2.5cm, font=\scriptsize}
]
\def\dx{6cm}
\begin{scope}[shift={(0.4*\dx,0cm)}]
\node at (0,2cm) {\normalsize\textbf{Wind}};
\node (wind1) [circleNode] {ID: 1 \\ Probability: 1\\Efficiency: 45\%\\Cost: \$500,000 };
\node (wind2) [circleNode, below=0.5cm of wind1] {ID: 2 \\ Probability: 1\\Efficiency: 45\%\\Cost: \$450,000 };
\draw (wind1) -- (wind2);
\end{scope}
\begin{scope}[shift={(1.4*\dx,0cm)}]
\node at (0,2cm) {\normalsize\textbf{Solar}};
\node (sol1) [circleNode] {ID: 1 \\ Probability: 1\\Efficiency: 20\%\\Cost: \$100,000 };
\node (sol2) [circleNode, below left=1.4cm and -0.4cm of sol1] {ID: 2\\Probability: 0.3\\Efficiency: 21\%\\Cost: \$85,000};
\node (sol3) [circleNode, below right=1.4cm and -0.4cm of sol1] {ID: 3\\Probability: 0.7\\Efficiency: 23\%\\Cost: \$80,000};
\draw (sol1) -- node[midway, left] {$\theta^{s}_{solar}$} (sol2);
\draw (sol1) -- node[midway, right] {$\theta^{f}_{solar}$} (sol3);
\end{scope}
\begin{scope}[shift={(2.6*\dx,0cm)}]
\node at (0,2cm) {\normalsize\textbf{Battery}};
\node (battery1) [circleNode] {ID: 1 \\ Probability: 1\\Efficiency: 40\%\\Cost: \$350};
\node (battery2) [circleNode, below left=1.4cm and -0.4cm of battery1] {ID: 2\\Probability: 0.4\\Efficiency: 42\%\\Cost: \$300};
\node (battery3) [circleNode, below right=1.4cm and -0.4cm of battery1] {ID: 3\\Probability: 0.6\\Efficiency: 45\%\\Cost: \$250};
\draw (battery1) -- node[midway, left] {$\theta^{s}_{battery}$} (battery2);
\draw (battery1) -- node[midway, right] {$\theta^{f}_{battery}$} (battery3);

\draw[-, decorate, decoration={brace, amplitude=5pt}, thick] 
    (-16, -1) -- (-16, 1) node[midway, xshift=-1.1cm] {Stage 1};
\draw[-, decorate, decoration={brace, amplitude=5pt}, thick] 
    (-16, -4) -- (-16, -2) node[midway, xshift=-1.1cm] {Stage 2};

\end{scope}
\end{tikzpicture}
\caption{Examples of technology trees with two stages.}
\label{fig:tech_trees}
\end{figure}
}

Since each technology advances independently, it is essential to consider all possible combinations of their respective advancement trajectories. This is accomplished by constructing a scenario tree \cite{shapiro2021}, wherein each stage consists of nodes formed by taking the Cartesian product of the nodes from the individual technology trees at that stage. Consequently, each scenario node represents a unique joint realization of the technological states across all considered technologies. The probability associated with each scenario node is computed as the product of the corresponding branch probabilities from the individual technology trees. The root node at stage 0 serves as the initialization point for the optimization model and adopts parameter values consistent with those of its immediate successors. For example, combining the three-stage versions of the three distinct technology trees depicted in Figure~\ref{fig:tech_trees} yields a fully expanded three-stage scenario tree as illustrated in Figure~\ref{fig:scenario_tree}.

Combining the three-stage versions of the three technology trees illustrated in Figure~\ref{fig:tech_trees} results in a fully expanded scenario tree, as shown in Figure~\ref{fig:scenario_tree}. Since two branches are defined for two of the technologies, there are four branches extending from each scenario node. Specifically, the probability associated with scenario node ID 3 at the second stage is calculated as \(\pi_{3} = \theta^{s}_{solar} \times \theta^{f}_{battery}\), and this node corresponds to technology nodes with IDs 2, 2, and 3 in the wind, solar, and battery technology trees, respectively.

Each path from the root to a leaf node defines a distinct scenario path, denoted by $S_{i}$, representing a trajectory of technological advancements over the planning horizon. For instance, scenario path $S_{3}$, characterized by the sequence \(ss \times fs\), indicates a trajectory in which the solar technology experiences slow advancement followed by fast advancement, while the battery technology undergoes two consecutive stages of slow improvement.

\newcommand{\smalltimes}{\mathbin{\mathsmaller{\times}}}
\newcommand{\scriptcompact}[1]{
  {\scriptsize
   \spaceskip=.12em
   \xspaceskip=.2em
   #1}
}

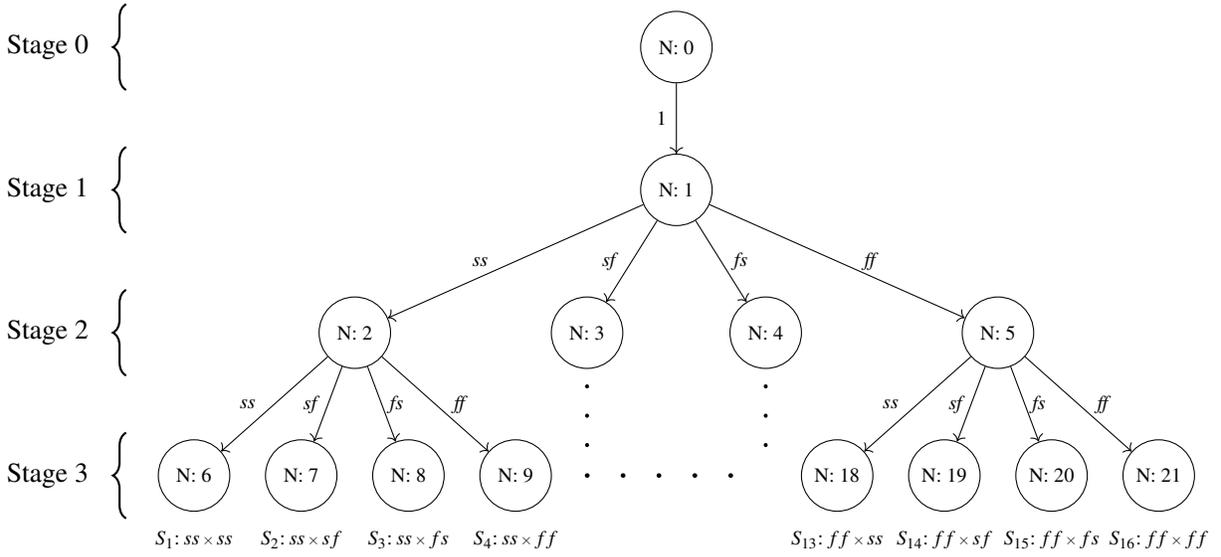
\begin{figure}[H]
\centering
\begin{tikzpicture}[scale=0.95, transform shape, ->, 
    node style/.style={draw, font=\scriptsize, align=center, text centered, minimum width=1cm, minimum height=1cm},
    edge label/.style={font=\normalsize}]

\node (0) [node style, circle] at (0, 0) {N: 0};
\node (1) [node style, circle] at (0, -2) {N: 1};
\node (2) [node style, circle] at (-4.5, -4) {N: 2};
\node (3) [node style, circle] at (-1.25, -4) {N: 3};
\node (4) [node style, circle] at (1.25, -4) {N: 4};
\node (5) [node style, circle] at (4.5, -4) {N: 5};
\node (6) [node style, circle, 
label={[yshift=-4pt]below:{\scriptcompact {$S_{1}$: \(ss \smalltimes ss\)}}}] at (-6.75, -6) {N: 6};
\node (7) [node style, circle,
label={[yshift=-4pt]below:{\scriptcompact {$S_{2}$: \(ss \smalltimes sf\)}}}] at (-5.25, -6) {N: 7};
\node (8) [node style, circle,
label={[yshift=-4pt]below:{\scriptcompact {$S_{3}$: \(ss \smalltimes fs\)}}}] at (-3.75, -6) {N: 8};
\node (9) [node style, circle,
label={[yshift=-4pt]below:{\scriptcompact {$S_{4}$: \(ss \smalltimes ff\)}}}] at (-2.25, -6) {N: 9};
\node (18) [node style, circle,
label={[yshift=-4pt]below:{\scriptcompact {$S_{13}$: \(ff \smalltimes ss\)}}}] at (2.25, -6) {N: 18};
\node (19) [node style, circle,
label={[yshift=-4pt]below:{\scriptcompact {$S_{14}$: \(ff \smalltimes sf\)}}}] at (3.75, -6) {N: 19};
\node (20) [node style, circle,
label={[yshift=-4pt]below:{\scriptcompact {$S_{15}$: \(ff \smalltimes fs\)}}}] at (5.25, -6) {N: 20};
\node (21) [node style, circle,
label={[yshift=-4pt]below:{\scriptcompact {$S_{16}$: \(ff \smalltimes ff\)}}}] at (6.75, -6) {N: 21};

\draw (0) -- (1) node [midway, left, edge label] {\scriptsize 1};
\draw (1) -- (2) node [midway, left, shift={(-.25,0)}, edge label] {\scriptsize \textit{ss}};
\draw (1) -- (3) node [midway, left, shift={(-.05,0)}, edge label] {\scriptsize \textit{sf}};
\draw (1) -- (4) node [midway, right, shift={(.05,0)}, edge label] {\scriptsize \textit{fs}};
\draw (1) -- (5) node [midway, right, shift={(.25,0)}, edge label] {\scriptsize \textit{ff}};
\draw (2) -- (6) node [midway, left, shift={(-.13,0)}, edge label] {\scriptsize \textit{ss}};
\draw (2) -- (7) node [midway, left, shift={(.02,0)}, edge label] {\scriptsize \textit{sf}};
\draw (2) -- (8) node [midway, right, shift={(-.02,0)}, edge label] {\scriptsize \textit{fs}};
\draw (2) -- (9) node [midway, right, shift={(.13,0)}, edge label] {\scriptsize \textit{ff}};
\draw (5) -- (18) node [midway, left, shift={(-.13,0)}, edge label] {\scriptsize \textit{ss}};
\draw (5) -- (19) node [midway, left, shift={(.05,0)}, edge label] {\scriptsize \textit{sf}};
\draw (5) -- (20) node [midway, right, shift={(-.03,0)}, edge label] {\scriptsize \textit{fs}};
\draw (5) -- (21) node [midway, right, shift={(.13,0)}, edge label] {\scriptsize \textit{ff}};

\draw[-, decorate, decoration={brace, amplitude=5pt}, thick] 
    (-7.7, -0.6) -- (-7.7, 0.6) node[midway, xshift=-1.1cm] {Stage 0};
\draw[-, decorate, decoration={brace, amplitude=5pt}, thick] 
    (-7.7, -2.6) -- (-7.7, -1.4) node[midway, xshift=-1.1cm] {Stage 1};
\draw[-, decorate, decoration={brace, amplitude=5pt}, thick] 
    (-7.7, -4.6) -- (-7.7, -3.4) node[midway, xshift=-1.1cm] {Stage 2};
\draw[-, decorate, decoration={brace, amplitude=5pt}, thick] 
    (-7.7, -6.6) -- (-7.7, -5.4) node[midway, xshift=-1.1cm] {Stage 3};

\draw[-, dash pattern=on 0pt off 13.47pt, line width=1.75pt, line cap=round] (-1.25,-6) -- (1.25,-6);

\draw[-, dash pattern=on 0pt off 11pt, line width=1.5pt, line cap=round] (-1.25,-4.75) -- (-1.25,-5.75);

\draw[-, dash pattern=on 0pt off 11pt, line width=1.5pt, line cap=round] (1.25,-4.75) -- (1.25,-5.75);

\end{tikzpicture}
\vspace{-0.012\textwidth}
\caption{Scenario tree instance with three stages.} 
\label{fig:scenario_tree}
\end{figure}

By embedding the constructed scenario trees into the multi-stage stochastic programming framework, decision-making concerning the transition and deployment of renewable electricity technologies is optimized while explicitly accounting for uncertainties in technological progress. This integration facilitates the determination of the optimal timing, scale, and mix of technology investments across a diverse set of technological advancement trajectories, thereby enabling adaptive and economically efficient strategic planning.

\section{Computational Experiments}
\label{sec:computational_experiments}

This section details the computational experiments conducted using the multi-stage stochastic programming framework, specifically tailored to the METU case study, accompanied by a comprehensive sensitivity analysis. All experiments are executed using the Gurobi solver integrated through Python, on a 64-bit workstation equipped with two Intel® Xeon® Gold 6248R CPUs operating at 3.00 GHz, and 256 GB of RAM. 
All runs use a 24-hour time limit and a 1\% optimality-gap tolerance. Here, the optimality gap is defined as $100\times(1-\frac{\text{LB}}{\text{UB}})$, where $\text{UB}$ is the objective function value of the incumbent solution reported by the solver and $\text{LB}$ is the proven lower bound given by the solver obtained via the branch-and-bound procedure.


\subsection{Base Case}
\label{sec:compexBase}

In the \texttt{Base Case}, we obtain a dynamic strategic plan of METU’s transition to a clean campus, assuming a planning horizon of 15 years (2026-2040). This is consistent with METU's  target of having a  100\% renewable electricity supply by 2040. The planning horizon is divided into three stages, each spanning a five-year interval. Each year comprises 4,368 bihourly sub-periods, corresponding to consecutive intervals during which electricity load, generation, and storage activities are aggregated. Investment decisions are taken yearly whereas the operational decisions are taken bihourly.
We utilize the input data described in Sections~\ref{subsec:Electricity_Technologies} and~\ref{subsec:CaseParameters}.
We set the investment budget as 5 million USD for the first five-year period, followed by an increased budget of 10 million USD for the subsequent ten years. {All costs are discounted using an annual discount factor of $\beta=0.97$}. In alignment with METU’s sustainability objectives, the maximum allowable emissions level is set to zero at the end of the 15-year planning horizon. In this context, achieving zero emissions for electricity implies complete independence from grid electricity, as the model considers only renewable technology installations and treats grid electricity as a carbon-emitting source due to its predominantly conventional generation mix. 
Technological advancements of both solar photovoltaic and lithium-ion battery systems are modeled using two distinct scenarios, categorized as slow and fast technological advancements for each stage. These scenarios are structured within a scenario tree comprising 16 leaf nodes and 16 corresponding scenario paths, as shown in Figure~\ref{fig:scenario_tree}.
%
Our optimization model comprises 1,400,967 constraints, 1,864,027 continuous variables, and 9,989 integer variables.
The base optimization model yields an objective value of 69.3 million USD, representing the total expected cost, including installation (\$47.1 M), grid electricity procurement (\$20.0 M) and  O\&M costs (\$2.2 M) across all scenario paths. 

The decision tree presented in Figure~\ref{fig:decision_tree_base_model} clearly illustrates annual technology installations for solar PV, wind power, and usable battery (i.e., the effective storage capacity available after accounting for depth-of-discharge limits) at each year and scenario node throughout the planning horizon, which can serve as the dynamic strategic plan for METU. We can roughly categorize the installation decisions as \textit{economically-driven} and \textit{environmentally-driven}. Economically-driven decisions, which are typically made in the first years of each stage, take advantage of the reduced costs and improved efficiencies. For example, we   observe technology installations in the first two years of the third, fourth, and fifth nodes (i.e., years 6-7), and the allowed budget is not utilized at all in year 8. This suggests that at those technological levels, the model opts for \textit{cleaning} only some portion of the demand, and waits for future technological advancements to make more economical decisions. In contrast, environmentally-driven decisions are typically made in the last years of the planning horizon just to meet the emission target. Note that such decisions might be quite costly: Even in scenario $S_{16}$, in which the technology improvements are always fast,  12-\si{\mega\watt} wind turbine, 12-\si{\mega\watt} solar, and 21.5-\si{\mega\watt\hour} battery investments are made in year 15, which have a cost of more than 6 million USD. Environmentally-driven decisions are observed not only in the final stage but also in node 2, where deployment occurs in every year of the interval (i.e., years 6–10). These deployments are driven by the insufficiency of the total budget in years 11–15 to meet the emission target in scenario path $S_1$. Consequently, the model undertakes additional installations in the preceding stage. This behavior signifies the intrinsic complications of simultaneously accounting for emission targets, budget limitations and economical considerations at the same time.

\colorlet{solaryellow}{yellow!8}
\colorlet{windcyan}{cyan!8}
\colorlet{batterypurple}{red!7}
\colorlet{gridgray}{gray!16}
\colorlet{renewablecolor}{green!8}

\clearpage
\begin{landscape}
  \vspace*{\fill}
    \begin{figure}[H]
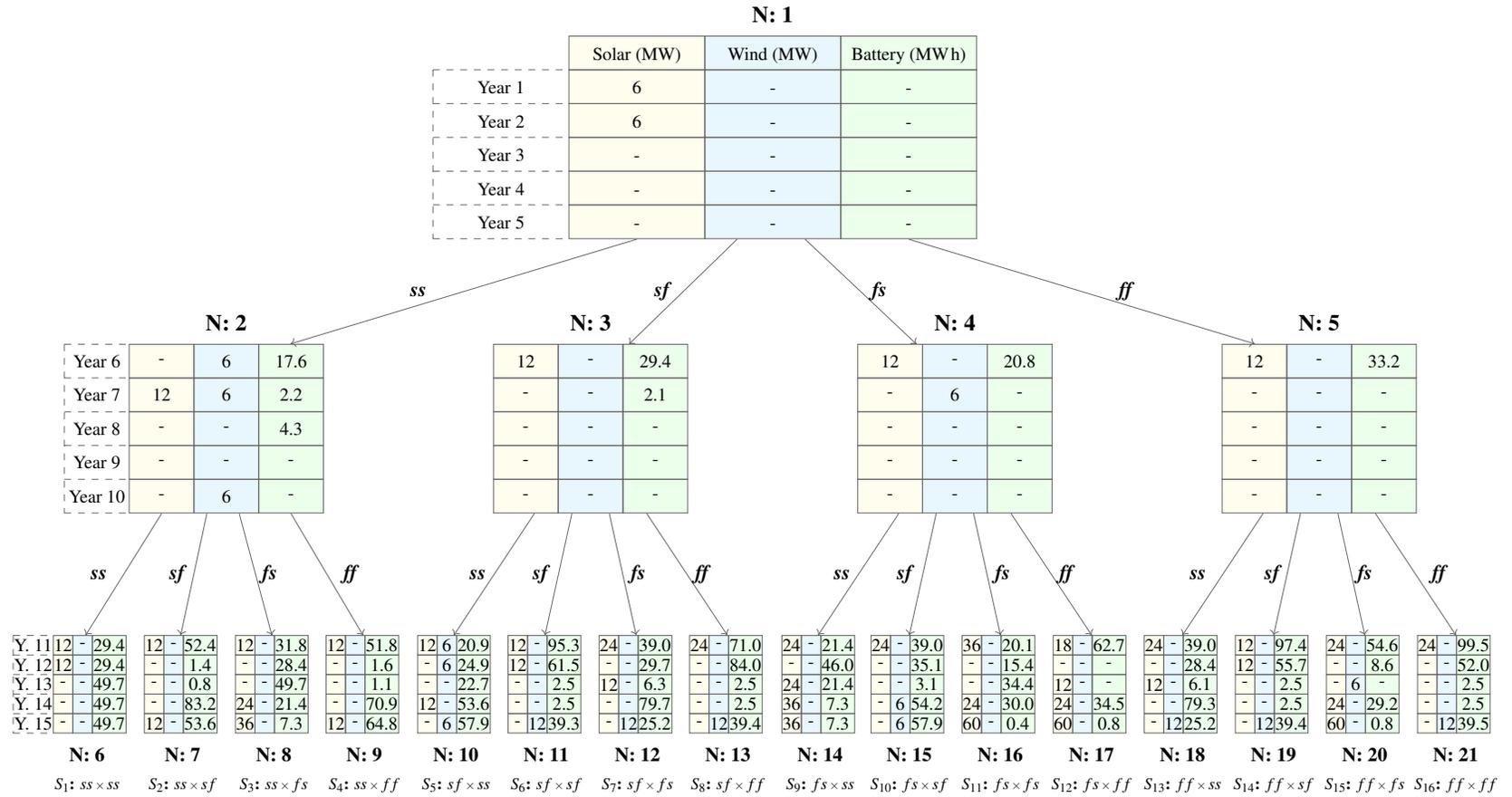

    \centering

    };
    \draw ([xshift=-1.05cm]1.south) -- ([xshift=1.5cm]2.north) node [midway, left, shift={(-.45,0)}, edge label] {\footnotesize \textbf{\textit{ss}}};
    \draw ([xshift=0.5cm]1.south) -- ([xshift=0.6cm]3.north) node [midway, left, shift={(-.05,0)}, edge label] {\footnotesize \textbf{\textit{sf}}};
    \draw ([xshift=1.55cm]1.south) -- ([xshift=-0.6cm]4.north) node [midway, right, shift={(.05,0)}, edge label] {\footnotesize \textbf{\textit{fs}}};
    \draw ([xshift=3.15cm]1.south) -- ([xshift=-1cm]5.north) node [midway, right, shift={(.45,0)}, edge label] {\footnotesize \textbf{\textit{ff}}};
    \draw ([xshift=-0.5cm]2.south) -- ([xshift=0.25cm]6.north) node [midway, left, shift={(-.13,0)}, edge label] {\footnotesize \textbf{\textit{ss}}};
    \draw ([xshift=0.2cm]2.south) -- (7.north) node [midway, left, shift={(.02,0)}, edge label] {\footnotesize \textbf{\textit{sf}}};
    \draw ([xshift=0.7cm]2.south) -- (8.north) node [midway, right, shift={(-.02,0)}, edge label] {\footnotesize \textbf{\textit{fs}}};
    \draw ([xshift=1.5cm]2.south) -- (9.north) node [midway, right, shift={(.13,0)}, edge label] {\footnotesize \textbf{\textit{ff}}};
    \draw (3) -- (10.north) node [midway, left, edge label] {\footnotesize \textbf{\textit{ss}}};
    \draw (3) -- (11.north) node [midway, left, edge label] {\footnotesize \textbf{\textit{sf}}};
    \draw (3) -- (12.north) node [midway, right, edge label] {\footnotesize \textbf{\textit{fs}}};
    \draw (3) -- (13.north) node [midway, right, edge label] {\footnotesize \textbf{\textit{ff}}};
    \draw (4) -- (14.north) node [midway, left, edge label] {\footnotesize \textbf{\textit{ss}}};
    \draw (4) -- (15.north) node [midway, left, edge label] {\footnotesize \textbf{\textit{sf}}};
    \draw (4) -- (16.north) node [midway, right, edge label] {\footnotesize \textbf{\textit{fs}}};
    \draw (4) -- (17.north) node [midway, right, edge label] {\footnotesize \textbf{\textit{ff}}};
    \draw (5) -- (18.north) node [midway, left, shift={(-.13,0)}, edge label] {\footnotesize \textbf{\textit{ss}}};
    \draw (5) -- (19.north) node [midway, left, shift={(.05,0)}, edge label] {\footnotesize \textbf{\textit{sf}}};
    \draw (5) -- (20.north) node [midway, right, shift={(-.03,0)}, edge label] {\footnotesize \textbf{\textit{fs}}};
    \draw (5) -- (21.north) node [midway, right, shift={(.13,0)}, edge label] {\footnotesize \textbf{\textit{ff}}};
    \node[above=2pt of 1, xshift=1.05cm] {\small \textbf{N: 1}};
    \node[above=2pt of 2, xshift=0.5cm] {\small \textbf{N: 2}};
    \node[above=2pt of 3] {\small \textbf{N: 3}};
    \node[above=2pt of 4] {\small \textbf{N: 4}};
    \node[above=2pt of 5] {\small \textbf{N: 5}};
    \end{tikzpicture}
    \vspace{-0.012\textwidth}
    \caption{Decision tree for the \texttt{Base Case} (see Figure~\ref{fig:scenario_tree}  for the corresponding scenario tree). Node indices and stage transitions are written in boldface. For each node in the tree, we provide a $5\times3$ matrix, in which a cell gives the installation decision for  solar (\si{\mega\watt}), wind (\si{\mega\watt}), and usable battery storage (\si{\mega\watt\hour}) technologies (given by the columns) for each   investment period (indexed by the rows). For example, 6 \si{\mega\watt} solar PV, 6 \si{\mega\watt} wind turbine and 21.0 \si{\mega\watt\hour} battery storage are installed in year 6 of node 3.} 
    \label{fig:decision_tree_base_model}
    \end{figure}
  \vspace*{\fill}
\end{landscape}
\clearpage

Over the 15-year planning horizon, procuring electricity exclusively from the grid without any technology installations results in a total cost of 58.8 million USD. By contrast, the optimization model yields a cost of 69.3 million USD, corresponding to a 17.8\% increase. The additional 10.5 million USD represents the premium required to facilitate the transition to clean electricity technologies while ensuring full demand satisfaction. From a cost–benefit perspective, this premium internalizes environmental benefits, reduces long-term dependence on grid electricity, and enhances energy self-sufficiency. The model prescribes a phased transition toward 100\% renewable technologies, such that no early retirements or salvage decisions are economically justified, since all installed assets are retained to meet the final year’s demand entirely with renewable sources. These findings underscore the economic viability of renewable deployment even within a 15-year horizon. Moreover, the extended lifetimes of renewable assets beyond the modeled period indicate that benefits will continue to accrue after the planning horizon, further strengthening the long-term cost-effectiveness of the clean energy transition.

The median usable battery storage capacity across all scenario paths is 338.2-\si{\mega\watt\hour}. The full effective capacity is sufficient to cover 42.9 consecutive average sub-period demand, underscoring both the critical role and heavy reliance on storage technologies in enabling system flexibility. However, the prominence of battery storage in the optimization outcomes also highlights a methodological limitation: relying on representative days for modeling is inadequate to capture the full extent of storage dynamics. Since battery performance and utilization are inherently tied to consecutive sub-periods variations, simplifications through representative days risk underestimating or misrepresenting the system’s actual dependence on storage. This result therefore emphasizes the need for modeling approaches that preserve chronological detail to more accurately reflect the operational and economic significance of clean energy transition.

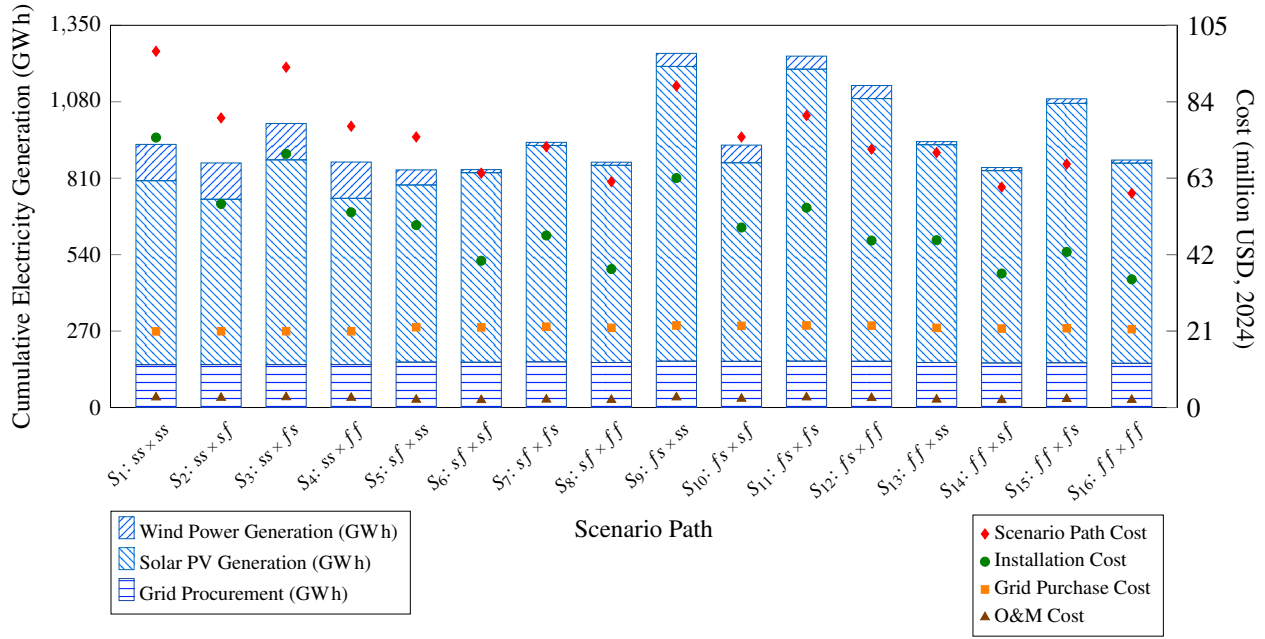
\begin{figure}[H]
\centering
\begin{tikzpicture}[scale=0.9]
  \begin{axis}[
      width=\textwidth,
      height=0.42\textwidth,
      ymin=0,
      ymax=1900,
      xmin=0.3,
      xmax=16.7,
      bar width=16pt,
      ybar stacked,
      ylabel={Cumulative Electricity Generation (\si{\giga\watt\hour})},
      xlabel={Scenario Path},
      xtick={1,2,3,4,5,6,7,8,9,10,11,12,13,14,15,16},
      xticklabels={
      $S_{1}$: \(ss \smalltimes ss\),
      $S_{2}$: \(ss \smalltimes sf\),
      $S_{3}$: \(ss \smalltimes fs\),
      $S_{4}$: \(ss \smalltimes ff\),
      $S_{5}$: \(sf \smalltimes ss\),
      $S_{6}$: \(sf \smalltimes sf\),
      $S_{7}$: \(sf \smalltimes fs\),
      $S_{8}$: \(sf \smalltimes ff\),
      $S_{9}$: \(fs \smalltimes ss\),
      $S_{10}$: \(fs \smalltimes sf\),
      $S_{11}$: \(fs \smalltimes fs\),
      $S_{12}$: \(fs \smalltimes ff\),
      $S_{13}$: \(ff \smalltimes ss\),
      $S_{14}$: \(ff \smalltimes sf\),
      $S_{15}$: \(ff \smalltimes fs\),
      $S_{16}$: \(ff \smalltimes ff\)},
      xticklabel style={yshift=-5pt, xshift=7pt, font=\scriptsize, rotate=45, anchor=east},
      xlabel style={yshift=0.1cm, font=\fontsize{10pt}{10pt}\selectfont},
      tick label style={font=\footnotesize},
      xtick style={draw=none},
      ytick={0, 380, 760, 1140, 1520, 1900},
      scaled y ticks=false,
      ylabel style={xshift=0cm, yshift=-0.06cm, font=\fontsize{10pt}{10pt}\selectfont},
      legend style={legend columns=2, font=\fontsize{7pt}{7pt}\selectfont, at={(-0.005,-0.28)}, anchor=north west, legend cell align={left}, /tikz/every even column/.append style={column sep=0.4em}},
      legend reversed,
  ]
      \addplot+ [
        ybar,
        fill=gridgray,
        draw=black,
        postaction={pattern=sparsehorizontallines, pattern color=black}
    ] coordinates {
(1, 144.73497743090863)
(2, 146.30454775209373)
(3, 145.33170583536983)
(4, 146.30454775209373)
(5, 167.8144281168286)
(6, 169.73959542770032)
(7, 166.44413729301928)
(8, 166.76623144846565)
(9, 152.96256000239953)
(10, 154.9837454160114)
(11, 156.00715033253306)
(12, 156.04541640270494)
(13, 163.01184734147608)
(14, 163.942144206494)
(15, 163.80945421777824)
(16, 163.26833394511144)
    };
    \addlegendentry{Grid Procurement}
    \addplot+ [
        ybar,
        fill=solaryellow,
        draw=black,
        postaction={pattern=sparsenorthwestlines, pattern color=black}
    ] coordinates {
(1, 851.0371767825849)
(2, 829.5905273338342)
(3, 990.4479863079846)
(4, 829.5905273338342)
(5, 624.2143435129063)
(6, 581.8572108516237)
(7, 665.4081991190809)
(8, 642.6258122257462)
(9, 1450.105822067629)
(10, 810.0128004428173)
(11, 1338.5751630693549)
(12, 1086.2774895342964)
(13, 824.8775238874508)
(14, 699.4593327096221)
(15, 999.3691265519404)
(16, 801.093853778979)
    };
    \addlegendentry{Solar PV Generation}
    \addplot+ [
        ybar,
        fill=windcyan,
        draw=black,
        postaction={pattern=sparsenortheastlines, pattern color=black}
    ] coordinates {
(1, 128.27527514604)
(2, 128.27527514604)
(3, 128.27527514604)
(4, 128.27527514604)
(5, 61.58261421720001)
(6, 61.58261421720001)
(7, 66.95471460636001)
(8, 61.58261421720001)
(9, 50.663710987200005)
(10, 61.58261421720001)
(11, 50.663710987200005)
(12, 50.663710987200005)
(13, 10.91890323)
(14, 27.03520439748)
(15, 16.11630116748)
(16, 10.91890323)
};
\addlegendentry{Wind Power Generation}
  \end{axis}
  \begin{axis}[
      width=\textwidth,
      height=0.42\textwidth,
      ymin=0,
      ymax=105,
      xmin=0.3,
      xmax=16.7,
      axis y line*=right,
      axis x line=none,
      ytick={0, 21, 42, 63, 84, 105},
      ylabel={Cost (million USD, 2024)},
      ylabel near ticks,
      ylabel style={rotate=180, xshift=0.04cm, yshift=-0.20cm, font=\fontsize{10pt}{10pt}\selectfont},
      legend style={legend columns=2, font=\fontsize{7pt}{7pt}\selectfont, at={(0.702,-0.28)}, anchor=north west, legend cell align={left}, /tikz/every even column/.append style={column sep=0.4em}},
  ]
      \addplot[mark=diamond*, red, mark size=2pt, only marks] coordinates {
(1, 94.937313188977)
(2, 79.03300603094891)
(3, 91.2962017694786)
(4, 79.03300603094891)
(5, 70.89202791919222)
(6, 62.29143316580303)
(7, 67.23576116125942)
(8, 58.84017442593553)
(9, 81.04142914499555)
(10, 70.7897443376794)
(11, 73.39819656478254)
(12, 66.00228098472806)
(13, 65.04764504088017)
(14, 57.637088342533254)
(15, 59.35824774812731)
(16, 54.47209347815759)
      };
      \addlegendentry{Scenario Path Cost}
      \addplot[mark=*, green!50!black, mark size=1.8pt, only marks] coordinates {
(1, 73.74708567247859)
(2, 57.70343526241119)
(3, 69.93193670336649)
(4, 57.70343526241119)
(5, 47.99958153242678)
(6, 39.26745543817863)
(7, 44.43379060411933)
(8, 36.0513927245024)
(9, 58.794096640528835)
(10, 48.967054127175246)
(11, 51.45484067442219)
(12, 43.89017072391701)
(13, 42.729844052786035)
(14, 35.328511918788)
(15, 36.84875781462435)
(16, 32.192080961566504)
      };
      \addlegendentry{Installation Cost}
      \addplot[mark=square*, orange, mark size=1.5pt, only marks] coordinates {
(1, 18.568877789803167)
(2, 18.72588197533729)
(3, 18.629471028525153)
(4, 18.72588197533729)
(5, 21.03916982376097)
(6, 21.22720098386603)
(7, 20.89901266714612)
(8, 20.924676972425335)
(9, 19.45492777233136)
(10, 19.658070896352573)
(11, 19.762363006388075)
(12, 19.76279300647973)
(13, 20.47566376096524)
(14, 20.56678298059759)
(15, 20.555281094594452)
(16, 20.49437510870829)
      };
      \addlegendentry{Grid Purchase Cost}
      \addplot[mark=triangle*, orange!50!black, mark size=2pt, only marks] coordinates {
(1, 2.621349726695231)
(2, 2.603688793200422)
(3, 2.7347940375869575)
(4, 2.603688793200422)
(5, 1.853276563004464)
(6, 1.7967767437583664)
(7, 1.902957889993966)
(8, 1.8641047290077917)
(9, 2.792404732135356)
(10, 2.164619314151587)
(11, 2.2641366897747144)
(12, 2.349317254331324)
(13, 1.842137227128898)
(14, 1.7417934431476565)
(15, 1.9542088389085144)
(16, 1.7856374078828003)
      };
      \addlegendentry{O\&M Cost}
  \end{axis}
\end{tikzpicture}
\caption{Cumulative electricity generation and grid electricity procurement, along with scenario path cost components, for each scenario path under the \texttt{Base Case}. Diagonally patterned bars represent cumulative renewable electricity generation from wind power and solar PV technologies, while horizontally lined bars denote cumulative grid electricity procurement summed across the entire planning horizon}.
\label{fig:scenario_paths_base}
\end{figure}

Figure~\ref{fig:scenario_paths_base} depicts 
electricity generation from solar PV and wind power technologies,  and  
the total electricity demand satisfied by  grid procurement in each scenario path. 
This figure also
provides a detailed breakdown of costs for each scenario path, distinctly categorizing installation expenditures, grid procurement costs, and O\&M expenses. 
%
Installation costs constitute the largest portion of total scenario path expenditures, underscoring the capital-intensive nature of renewable energy technology deployment, while O\&M costs account for a comparatively minor share. 

The total installed capacities of renewable and storage technologies, along with scenario path costs, demonstrate significant sensitivity to technological advancement trajectories whereas the total grid procurement is insensitive to uncertainty realizations. 
{Across all scenario paths, 
renewable generation capacities exceed electricity load requirements significantly, particularly in the later years of the planning horizon. The need for substantial electricity generation capacity arises from notable temporal mismatches between renewable electricity generation and campus demand profiles. 
Consequently, electricity storage technologies are critical for managing excess generation and ensuring a stable balance between intermittent renewable supply and campus electricity demand.} Furthermore, wind power installations remain relatively limited due to lower simulated capacity factors at the campus location. Despite wind energy's comparatively limited economic competitiveness relative to solar PV systems, wind power is strategically implemented in smaller capacities to complement the generation profile provided by solar PV. 
As expected, the electricity generation from wind turbines is the highest in scenario paths $S_{1}$–$S_{4}$, in which both solar PV and storage technology improvements are slow after the first state transition.

Considering scenario paths characterized by slow first stage transition (i.e., the transition from Stage~1 to Stage~2) in solar technological advancement, the expected cost for the reduced scenario tree (i.e., $S_{1}$–$S_{8}$) is \$76.2 million, whereas for scenario paths involving fast first stage transition in solar advancements (i.e., $S_{9}$–$S_{16}$), the expected cost amounts to \$65.8 million.
In comparison, scenario paths regarding battery technology advancements in the first stage transition result  in expected costs of \$76.8 million for slow advancement (i.e., $S_{1}$–$S_{4}$, $S_{9}$–$S_{12}$) and \$60.6 million for fast advancement (i.e., $S_{5}$–$S_{8}$, $S_{13}$–$S_{16}$). These outcomes indicate that total scenario path costs are considerably more sensitive to advancements in battery technology, despite the relatively similar overall contributions from battery and solar installations to the total expected costs, amounting to \$20.6 million and \$20.5 million, respectively.

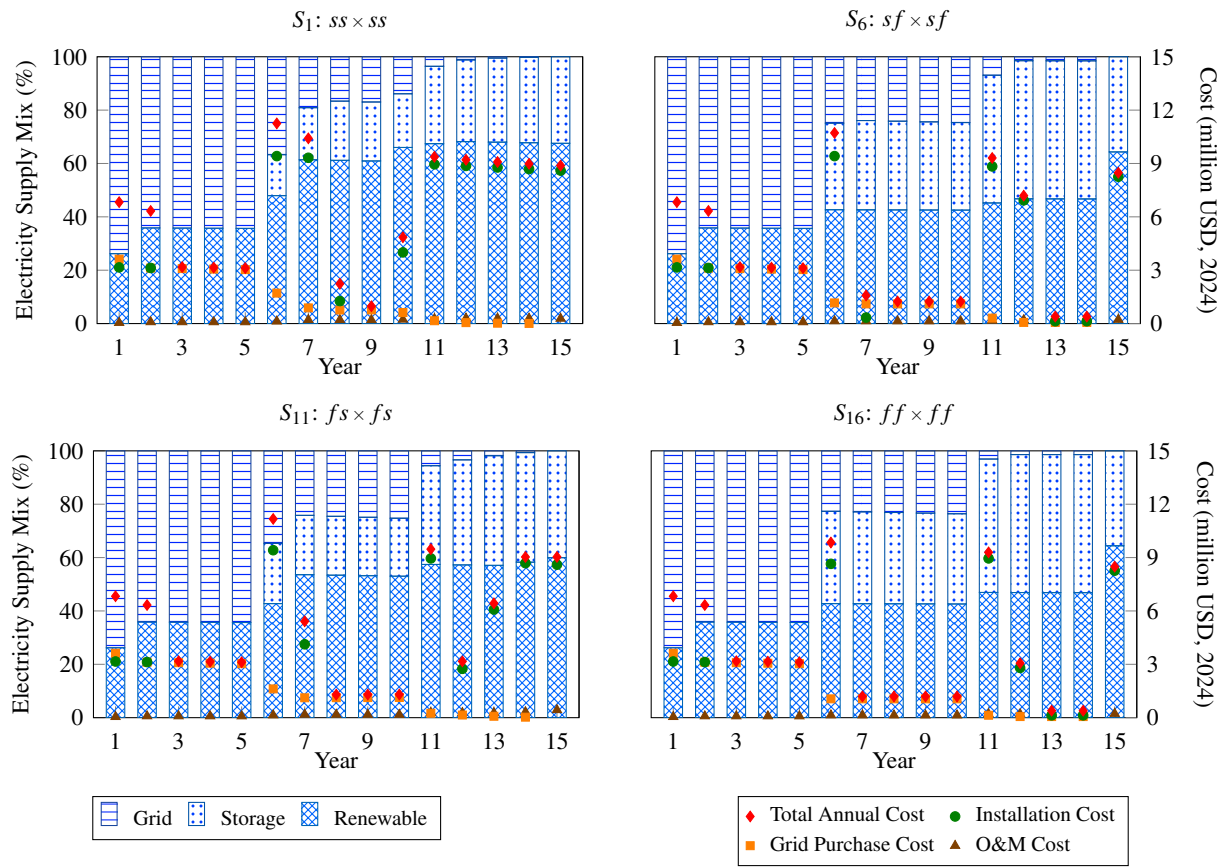
\begin{figure}[H]
\centering
\begin{subfigure}[t]{0.46\textwidth}
\centering
\begin{tikzpicture}[remember picture]
\pgfplotsset{every axis y label/.append style={rotate=0}}
\begin{axis}[
    xtick pos=bottom,
    width=0.85\linewidth,
    height=3.55cm,
    scale only axis,
    title={\fontsize{9pt}{10pt}\selectfont $S_{1}$: \(ss \smalltimes ss\)},
    xmin=0.3, xmax=15.7, xtick={1, 3, 5, 7, 9, 11, 13, 15},
    ymin=0,   ymax=100, title style={yshift=-4.5pt},
    ytick pos=left,
    ybar stacked,
    legend to name=sharedpercentage1,
    bar width=7pt, tick label style={font=\footnotesize},
    xlabel={\footnotesize Year}, xticklabel style={yshift=-2pt},
    legend columns=3,
    legend style={column sep=4pt, font=\fontsize{7pt}{7pt}\selectfont, at={(0,-0.23)}, anchor=north west, legend cell align={left}, legend reversed,},
    ylabel style={
        at={(axis description cs:-0.102,0.5)},
        anchor=south
    },
    ylabel={\footnotesize Electricity Supply Mix (\%)},
    xlabel style={yshift=0.22cm},
]

\addplot+ [ybar, 
          fill=renewablecolor,
          draw=black,
          postaction={pattern=sparsecrosshatch, pattern color=black}
]
coordinates {
(1, 26.203814530085367)
(2, 35.84267164665875)
(3, 35.78423824670897)
(4, 35.72573050326233)
(5, 35.666868469458684)
(6, 47.9502199031178)
(7, 53.32706181248954)
(8, 61.228883137993506)
(9, 67.72223246771422)
(10, 69.11992941699803)
(11, 68.91534736672075)
(12, 68.7058935627242)
(13, 68.49405853431871)
(14, 68.2798895002548)
(15, 68.63402101405107)
};
\addlegendentry{Renewable}
\addplot+ [ybar, 
            fill=batterypurple,
            draw=black,
            postaction={pattern=sparsedots, pattern color=black}
            ] 
coordinates {
(1, 0)
(2, 0)
(3, 0)
(4, 0)
(5, 0)
(6, 14.01906322748058)
(7, 27.634004186084926)
(8, 24.655041136980543)
(9, 23.480680475618755)
(10, 23.960180129797628)
(11, 29.597566129507616)
(12, 30.661255318102768)
(13, 31.187801104978842)
(14, 31.574347661667307)
(15, 31.365978985948928)
};
\addlegendentry{Storage}
\addplot+ [ybar,         
            fill=gridgray,
            draw=black,
            postaction={pattern=sparsehorizontallines, pattern color=black}
            ]
coordinates {
(1, 73.7961854699147)
(2, 64.15732835334131)
(3, 64.21576175329112)
(4, 64.2742694967376)
(5, 64.33313153054128)
(6, 38.03071686940162)
(7, 19.038934001425535)
(8, 14.116075725025953)
(9, 8.797087056667023)
(10, 6.9198904532043395)
(11, 1.487086503771632)
(12, 0.632851119173034)
(13, 0.3181403607024466)
(14, 0.1457628380778961)
(15, 0)
};
\addlegendentry{Grid}
\end{axis}
\pgfplotsset{every axis y label/.append style={rotate=180}}
\begin{axis}[
    width=0.85\linewidth,
    height=3.55cm,
    scale only axis,
    xmin=0.3, xmax=15.7,
    ymin=0,   ymax=10,
    ytick=\empty,
    legend reversed,
    legend to name=sharedcostcomp1,
    axis y line*=right,
    axis x line=none,
    legend columns=2,
    legend style={column sep=4pt, font=\fontsize{7pt}{7pt}\selectfont, at={(0,-0.23)}, anchor=north west, legend cell align={left}},
    ylabel style={
        at={(axis description cs:1.095,0.5)},
        anchor=south
    },
]

\addplot[mark=triangle*, orange!50!black, mark size=2pt, only marks] coordinates {
(1, 0.044638220590268216)
(2, 0.08305857542731857)
(3, 0.08056681816449901)
(4, 0.07814981361956404)
(5, 0.0758053192109771)
(6, 0.10953924817967578)
(7, 0.15083318797663722)
(8, 0.1736291683342705)
(9, 0.2275864062639121)
(10, 0.2837635753386835)
(11, 0.27525066807852305)
(12, 0.2669931480361673)
(13, 0.2589833535950823)
(14, 0.2512138529872298)
(15, 0.26133837089242223)
};
\addlegendentry{O\&M Cost}
\addplot[mark=square*, orange, mark size=1.5pt, only marks] coordinates {
(1, 3.549933785293645)
(2, 2.9936728851777707)
(3, 2.9065074872198258)
(4, 2.8218809715004376)
(5, 2.739731278018097)
(6, 1.5710120709512345)
(7, 0.762885521718837)
(8, 0.548658937496376)
(9, 0.3316645869323729)
(10, 0.2530644651093597)
(11, 0.05275211950209852)
(12, 0.021775942857091547)
(13, 0.010618567182936224)
(14, 0.004719170843094101)
};
\addlegendentry{Grid Purchase Cost}
\addplot[mark=*, green!50!black, mark size=1.8pt, only marks] coordinates {
(1, 3.0944653267999995)
(2, 3.0016313669959995)
(6, 8.329720049289996)
(7, 7.644997566394739)
(8, 3.4682995443629987)
(9, 7.160556017942948)
(10, 7.364832592227026)
(11, 7.153014030880803)
(12, 6.938423609954377)
(13, 6.7302709016557465)
(14, 6.528362774606073)
(15, 6.332511891367892)
};
\addlegendentry{Installation Cost}
\addplot[mark=diamond*, red, mark size=2pt, only marks] coordinates {
(1, 6.689037332683912)
(2, 6.078362827601089)
(3, 2.9870743053843247)
(4, 2.9000307851200016)
(5, 2.815536597229074)
(6, 10.010271368420906)
(7, 8.558716276090214)
(8, 4.1905876501936445)
(9, 7.719807011139233)
(10, 7.90166063267507)
(11, 7.481016818461424)
(12, 7.227192700847636)
(13, 6.999872822433765)
(14, 6.784295798436397)
(15, 6.593850262260314)
};
\addlegendentry{Total Annual Cost}
\end{axis}
\end{tikzpicture}
\label{fig:base_annual_subfig1}
\end{subfigure}
\hspace{0.05\textwidth}
\begin{subfigure}[t]{0.46\textwidth}
\centering
\begin{tikzpicture}[remember picture]
\pgfplotsset{every axis y label/.append style={rotate=0}}
\begin{axis}[
    xtick pos=bottom,
    width=0.85\linewidth,
    height=3.55cm,
    scale only axis,
    title={\fontsize{9pt}{10pt}\selectfont $S_{6}$: \(sf \smalltimes sf\)},
    xmin=0.3, xmax=15.7, xtick={1, 3, 5, 7, 9, 11, 13, 15},
    ymin=0,   ymax=100, title style={yshift=-4.5pt},
    ytick pos=left,
    ybar stacked,
    bar width=7pt, tick label style={font=\footnotesize},
    xlabel={\footnotesize Year}, xticklabel style={yshift=-2pt},
    ytick=\empty,
    xlabel style={yshift=0.22cm},
]

\addplot+ [ybar, 
          fill=renewablecolor,
          draw=black,
          postaction={pattern=sparsecrosshatch, pattern color=black}
]
coordinates {
(1, 26.203814530085367)
(2, 35.84267164665875)
(3, 35.78423824670897)
(4, 35.72573050326233)
(5, 35.666868469458684)
(6, 51.60858479469847)
(7, 51.40642171435924)
(8, 51.20236808953091)
(9, 50.99590700620292)
(10, 50.787910093640456)
(11, 53.89995155477956)
(12, 53.71453063535151)
(13, 55.20017338218165)
(14, 56.151409141393394)
(15, 68.76016432626527)
};
\addplot+ [ybar, 
            fill=batterypurple,
            draw=black,
            postaction={pattern=sparsedots, pattern color=black}
            ] 
coordinates {
(1, 0)
(2, 0)
(3, 0)
(4, 0)
(5, 0)
(6, 16.175806670834746)
(7, 22.55553683882999)
(8, 22.420266978201955)
(9, 22.28213030965742)
(10, 22.1336913931333)
(11, 36.639325717981315)
(12, 36.85220377413222)
(13, 40.93864097832197)
(14, 42.94333154272346)
(15, 31.23983567373473)
};
\addplot+ [ybar,         
            fill=gridgray,
            draw=black,
            postaction={pattern=sparsehorizontallines, pattern color=black}
            ]
coordinates {
(1, 73.7961854699147)
(2, 64.15732835334131)
(3, 64.21576175329112)
(4, 64.2742694967376)
(5, 64.33313153054128)
(6, 32.21560853446678)
(7, 26.03804144681077)
(8, 26.377364932267138)
(9, 26.721962684139662)
(10, 27.078398513226244)
(11, 9.460722727239121)
(12, 9.433265590516273)
(13, 3.8611856394963797)
(14, 0.9052593158831408)
(15, 0)
};
\end{axis}
\pgfplotsset{every axis y label/.append style={rotate=180}}
\begin{axis}[
    width=0.85\linewidth,
    height=3.55cm,
    scale only axis,
    xmin=0.3, xmax=15.7,
    ymin=0,   ymax=10,
    ytick={0, 2, 4, 6, 8, 10},
    axis y line*=right,
    axis x line=none,
    ylabel={\footnotesize Cost (million USD, 2024)}, tick label style={font=\footnotesize},
    ylabel style={
        at={(axis description cs:1.095,0.5)},
        anchor=south
    },
    legend style={at={(0.9,-0.1)}, anchor=north, font=\tiny, scale=0.8},
]

\addplot[mark=triangle*, orange!50!black, mark size=2pt, only marks] coordinates {
(1, 0.044638220590268216)
(2, 0.08305857542731857)
(3, 0.08056681816449901)
(4, 0.07814981361956404)
(5, 0.0758053192109771)
(6, 0.13456442498764548)
(7, 0.1305274922380161)
(8, 0.12661166747087563)
(9, 0.12281331744674936)
(10, 0.11912891792334689)
(11, 0.14361438904157303)
(12, 0.13930595737032583)
(13, 0.15738782341155033)
(14, 0.17249422015643928)
(15, 0.18810978669921732)
};
\addplot[mark=square*, orange, mark size=1.5pt, only marks] coordinates {
(1, 3.549933785293645)
(2, 2.9936728851777707)
(3, 2.9065074872198258)
(4, 2.8218809715004376)
(5, 2.739731278018097)
(6, 1.3307955791232269)
(7, 1.043338079337824)
(8, 1.0252266493608118)
(9, 1.0074617493913158)
(10, 0.9902729648840753)
(11, 0.33560467035223596)
(12, 0.32459175030499543)
(13, 0.12887474895751438)
(14, 0.029308384944261457)
};
\addplot[mark=*, green!50!black, mark size=1.8pt, only marks] coordinates {
(1, 3.0944653267999995)
(2, 3.0016313669959995)
(6, 8.32972004929)
(7, 1.2916583505105506)
(11, 5.769213571941215)
(12, 0.0714654736784643)
(13, 4.848426632988442)
(14, 6.528362774606073)
(15, 6.332511891367891)
};
\addplot[mark=diamond*, red, mark size=2pt, only marks] coordinates {
(1, 6.689037332683912)
(2, 6.078362827601089)
(3, 2.9870743053843247)
(4, 2.9000307851200016)
(5, 2.815536597229074)
(6, 9.795080053400872)
(7, 2.4655239220863905)
(8, 1.1518383168316875)
(9, 1.1302750668380652)
(10, 1.1094018828074221)
(11, 6.248432631335024)
(12, 0.5353631813537856)
(13, 5.134689205357507)
(14, 6.730165379706774)
(15, 6.5206216780671085)
};
\end{axis}
\end{tikzpicture}
\label{fig:base_annual_subfig2}
\end{subfigure}
\\[-0.05cm]
\vspace{-0.02\textwidth}
\begin{subfigure}[t]{0.46\textwidth}
\centering
\begin{tikzpicture}[remember picture]
\pgfplotsset{every axis y label/.append style={rotate=0}}
\begin{axis}[
    xtick pos=bottom,
    width=0.85\linewidth,
    height=3.55cm,
    scale only axis,
    title={\fontsize{9pt}{10pt}\selectfont $S_{11}$: \(fs \smalltimes fs\)},
    xmin=0.3, xmax=15.7, xtick={1, 3, 5, 7, 9, 11, 13, 15},
    ymin=0,   ymax=100, title style={yshift=-4.5pt},
    ytick pos=left,
    ybar stacked,
    bar width=7pt, tick label style={font=\footnotesize},
    xlabel={\footnotesize Year}, xticklabel style={yshift=-2pt},
    ylabel={\footnotesize Electricity Supply Mix (\%)},
    ylabel style={
        at={(axis description cs:-0.102,0.5)},
        anchor=south
    },
    xlabel style={yshift=0.22cm},
]

\addplot+ [ybar, 
          fill=renewablecolor,
          draw=black,
          postaction={pattern=sparsecrosshatch, pattern color=black}
]
coordinates {
(1, 26.203814530085367)
(2, 35.84267164665875)
(3, 35.78423824670897)
(4, 35.72573050326233)
(5, 35.666868469458684)
(6, 47.9502199031178)
(7, 53.525397757734645)
(8, 53.33893865634851)
(9, 55.2357069724915)
(10, 56.379661951692114)
(11, 55.71056154299549)
(12, 56.72950094441467)
(13, 56.556928531763525)
(14, 59.019257612737974)
(15, 60.21181723044592)
};
\addplot+ [ybar, 
            fill=batterypurple,
            draw=black,
            postaction={pattern=sparsedots, pattern color=black}
            ] 
coordinates {
(1, 0)
(2, 0)
(3, 0)
(4, 0)
(5, 0)
(6, 14.019063227480586)
(7, 27.2894861191726)
(8, 27.163219532065565)
(9, 29.11344246031007)
(10, 29.815406104169725)
(11, 35.40830252791709)
(12, 38.92324662290326)
(13, 41.22828053335351)
(14, 40.37611630495088)
(15, 39.78818276955408)
};
\addplot+ [ybar,         
            fill=gridgray,
            draw=black,
            postaction={pattern=sparsehorizontallines, pattern color=black}
            ]
coordinates {
(1, 73.7961854699147)
(2, 64.15732835334131)
(3, 64.21576175329112)
(4, 64.2742694967376)
(5, 64.33313153054128)
(6, 38.030716869401616)
(7, 19.185116123092754)
(8, 19.497841811585925)
(9, 15.650850567198434)
(10, 13.804931944138161)
(11, 8.881135929087424)
(12, 4.347252432682077)
(13, 2.2147909348829655)
(14, 0.60462608231115)
(15, 0)
};
\end{axis}
\pgfplotsset{every axis y label/.append style={rotate=180}}
\begin{axis}[
    width=0.85\linewidth,
    height=3.55cm,
    scale only axis,
    xmin=0.3, xmax=15.7,
    ymin=0,   ymax=10,
    ytick=\empty,
    axis y line*=right,
    axis x line=none,
    legend style={at={(0.91,-0.3)}, anchor=north, font=\tiny, scale=0.8},
    ylabel style={
        at={(axis description cs:1.095,0.5)},
        anchor=south
    },
]

\addplot[mark=triangle*, orange!50!black, mark size=2pt, only marks] coordinates {
(1, 0.044638220590268216)
(2, 0.08305857542731857)
(3, 0.08056681816449901)
(4, 0.07814981361956404)
(5, 0.0758053192109771)
(6, 0.10953924817967578)
(7, 0.15083318797663722)
(8, 0.1463081923373381)
(9, 0.1772868501949997)
(10, 0.20347062532049415)
(11, 0.16444920571595756)
(12, 0.18450833232652372)
(13, 0.17897308235672799)
(14, 0.252916015674968)
(15, 0.33363320267876545)
};
\addplot[mark=square*, orange, mark size=1.5pt, only marks] coordinates {
(1, 3.549933785293645)
(2, 2.9936728851777707)
(3, 2.9065074872198258)
(4, 2.8218809715004376)
(5, 2.739731278018097)
(6, 1.5710120709512339)
(7, 0.7687430042935288)
(8, 0.7578356322396026)
(9, 0.5900626940569154)
(10, 0.5048544831655116)
(11, 0.315044715056836)
(12, 0.14958576779186156)
(13, 0.07392305171932889)
(14, 0.01957517990348241)
};
\addplot[mark=*, green!50!black, mark size=1.8pt, only marks] coordinates {
(1, 3.0944653267999995)
(2, 3.0016313669959995)
(6, 8.329720049289996)
(7, 5.872691166063905)
(9, 2.4431500116314933)
(10, 2.3698555112825486)
(11, 2.5901179587274914)
(12, 4.162063716001034)
(13, 6.7302709016557465)
(14, 6.528362774606073)
(15, 6.332511891367892)
};
\addplot[mark=diamond*, red, mark size=2pt, only marks] coordinates {
(1, 6.689037332683912)
(2, 6.078362827601089)
(3, 2.9870743053843247)
(4, 2.9000307851200016)
(5, 2.815536597229074)
(6, 10.010271368420906)
(7, 6.792267358334072)
(8, 0.9041438245769406)
(9, 3.2104995558834086)
(10, 3.0781806197685544)
(11, 2.9864680736978557)
(12, 4.496157816119419)
(13, 6.983167035731803)
(14, 6.800853970184523)
(15, 6.666145094046657)
};
\end{axis}
\end{tikzpicture}
\label{fig:base_annual_subfig3}
\end{subfigure}
\hspace{0.05\textwidth}
\begin{subfigure}[t]{0.46\textwidth}
\centering
\begin{tikzpicture}[remember picture]
\pgfplotsset{every axis y label/.append style={rotate=0}}
\begin{axis}[
    xtick pos=bottom,
    width=0.85\linewidth,
    height=3.55cm,
    scale only axis,
    title={\fontsize{9pt}{10pt}\selectfont $S_{16}$: \(ff \smalltimes ff\)},
    xmin=0.3, xmax=15.7, xtick={1, 3, 5, 7, 9, 11, 13, 15},
    ymin=0,   ymax=100, title style={yshift=-4.5pt},
    ytick pos=left,
    ybar stacked,
    bar width=7pt, tick label style={font=\footnotesize},
    xlabel={\footnotesize Year}, xticklabel style={yshift=-2pt},
    ytick=\empty,
    xlabel style={yshift=0.22cm},
]

\addplot+ [ybar, 
          fill=renewablecolor,
          draw=black,
          postaction={pattern=sparsecrosshatch, pattern color=black}
]
coordinates {
(1, 26.203814530085367)
(2, 35.84267164665875)
(3, 35.78423824670897)
(4, 35.72573050326233)
(5, 35.666868469458684)
(6, 42.8548112527982)
(7, 42.821747957030496)
(8, 42.7884671863179)
(9, 42.75475893106295)
(10, 42.72028053617896)
(11, 47.110119382694954)
(12, 47.08659987003519)
(13, 47.06294592079107)
(14, 48.167527434406566)
(15, 65.94295326110955)
};
\addplot+ [ybar, 
            fill=batterypurple,
            draw=black,
            postaction={pattern=sparsedots, pattern color=black}
            ] 
coordinates {
(1, 0)
(2, 0)
(3, 0)
(4, 0)
(5, 0)
(6, 32.8344839280402)
(7, 32.643958545760725)
(8, 32.43751407214084)
(9, 32.22709398427928)
(10, 32.01089224601772)
(11, 46.823270793676556)
(12, 46.87754479281417)
(13, 46.86676407931229)
(14, 50.60984478828479)
(15, 34.05704673889045)
};
\addplot+ [ybar,         
            fill=gridgray,
            draw=black,
            postaction={pattern=sparsehorizontallines, pattern color=black}
            ]
coordinates {
(1, 73.7961854699147)
(2, 64.15732835334131)
(3, 64.21576175329112)
(4, 64.2742694967376)
(5, 64.33313153054128)
(6, 24.310704819161604)
(7, 24.53429349720878)
(8, 24.774018741541255)
(9, 25.018147084657773)
(10, 25.26882721780332)
(11, 6.0666098236284896)
(12, 6.035855337150645)
(13, 6.070289999896641)
(14, 1.2226277773086482)
(15, 0)
};
\end{axis}
\pgfplotsset{every axis y label/.append style={rotate=180}}
\begin{axis}[
    width=0.85\linewidth,
    height=3.55cm,
    scale only axis,
    xmin=0.3, xmax=15.7,
    ymin=0,   ymax=10,
    ytick={0, 2, 4, 6, 8, 10},
    axis y line*=right,
    axis x line=none,
    ylabel={\footnotesize Cost (million USD, 2024)}, tick label style={font=\footnotesize},
    ylabel style={
        at={(axis description cs:1.095,0.5)},
        anchor=south
    },
    legend style={at={(0.9,-0.1)}, anchor=north, font=\tiny, scale=0.8},
]

\addplot[mark=triangle*, orange!50!black, mark size=2pt, only marks] coordinates {
(1, 0.044638220590268216)
(2, 0.08305857542731857)
(3, 0.08056681816449901)
(4, 0.07814981361956404)
(5, 0.0758053192109771)
(6, 0.12358151325058721)
(7, 0.11987406785306959)
(8, 0.11627784581747751)
(9, 0.11278951044295318)
(10, 0.10940582512966458)
(11, 0.16224232768762778)
(12, 0.15737505785699893)
(13, 0.15265380612128898)
(14, 0.16790222338488575)
(15, 0.2013164833256197)
};
\addplot[mark=square*, orange, mark size=1.5pt, only marks] coordinates {
(1, 3.549933785293645)
(2, 2.9936728851777707)
(3, 2.9065074872198258)
(4, 2.8218809715004376)
(5, 2.739731278018097)
(6, 1.0042516646580337)
(7, 0.9830832594524322)
(8, 0.9629083227537271)
(9, 0.9432251113575008)
(10, 0.9240958779705976)
(11, 0.21520370575416156)
(12, 0.20768935525817742)
(13, 0.20260799994532955)
(14, 0.03958340434856603)
};
\addplot[mark=*, green!50!black, mark size=1.8pt, only marks] coordinates {
(1, 3.0944653267999995)
(2, 3.0016313669959995)
(6, 7.818932662720286)
(11, 5.6219962894506725)
(12, 0.04524137747993088)
(13, 0.05522709926179919)
(14, 6.222074947489925)
(15, 6.332511891367891)
};
\addplot[mark=diamond*, red, mark size=2pt, only marks] coordinates {
(1, 6.689037332683912)
(2, 6.078362827601089)
(3, 2.9870743053843247)
(4, 2.9000307851200016)
(5, 2.815536597229074)
(6, 8.946765840628906)
(7, 1.1029573273055018)
(8, 1.0791861685712045)
(9, 1.056014621800454)
(10, 1.0335017031002622)
(11, 5.999442322892462)
(12, 0.4103057905951072)
(13, 0.4104889053284177)
(14, 6.429560575223376)
(15, 6.533828374693511)
};
\end{axis}
\end{tikzpicture}
\label{fig:base_annual_subfig4}
\end{subfigure}
\vspace{0.05\textwidth}
\begin{tikzpicture}[remember picture, overlay]
  \node[anchor=north west] at (-15.25, -0.1)
        {\pgfplotslegendfromname{sharedpercentage1}};
  \node[anchor=north east] at (-1.2, -0.1)
        {\pgfplotslegendfromname{sharedcostcomp1}};
\end{tikzpicture}
\vspace{-5pt}
\caption{Annual electricity supply mix percentages and cost breakdowns for extreme scenario paths in the \texttt{Base Case}. The stacked bars show the annual proportions of electricity demand satisfied by renewable generation (solar PV and wind power), battery storage utilization, and grid electricity procurement. Cost components—including installation, grid procurement, and O\&M expenses—along with the total annual total costs, are depicted on the secondary (right) axis.}
\label{fig:extreme_scenario_paths_base}
\end{figure}

Figure~\ref{fig:extreme_scenario_paths_base} further presents a detailed annual breakdown of electricity source contributions and corresponding cost components across four representative scenario paths. 
These selected paths represent extreme technological advancement scenarios characterized by either two consecutive periods of slow or fast technological improvements for solar PV and storage technologies.
{Grid electricity procurement is similar across all scenario paths and is primarily concentrated in the initial years, common to all paths.} 
{As mentioned earlier in relation to Figure~\ref{fig:decision_tree_base_model}, annual installation patterns distinctly favor technology deployments at the onset of each stage. This trend aligns with anticipated technological advancements occurring between stages, incentivizing early investments to capitalize on projected cost reductions sooner as much as the budget allows, and thereby reduce reliance on grid electricity procurement.} We observe that the proportion of renewables in the annual supply mix increases steadily, starting from the initial years. On the other hand, the battery storage deployments typically occur later in the planning horizon. One of the underlying reasons for this strategy is to increase renewable generation first, and then, when surplus occurs, invest in battery storage to manage intermittency.
%



Notably, even though we only have a single emission target at the end of the planning horizon, the deployment of renewable electricity generation technologies begins immediately.
Across all scenario paths, emissions decrease to less than 26\% of the total supply mix in year 7 and further reduce to below 10\% in year 12. This outcome highlights a critical trade-off: early installations, despite incurring higher initial costs, extend the duration during which renewable generation meets electricity demand, whereas delayed installations, although initially less expensive, result in increased reliance on grid electricity procurement and consequently higher grid purchasing costs. Thus, these investment decisions demonstrate the economic competitiveness of renewable technologies relative to conventional grid electricity procurement, emphasizing their viability in both environmental and economic terms.
However, we also note that our stringent emission constraint enforced for the last year entails a large investment amount at the end of the planning horizon. For example, in $S_{16}$ where the technological advancements are always fast, the grid procurement amounts to 1.22\% of the total demand in year 14. As previously mentioned
in relation to Figure~\ref{fig:extreme_scenario_paths_base}, cleaning up this portion of the demand requires an investment of more than 6 million USD in the last year.

\subsection{Sensitivity Analysis}
\label{sec:sensitivityanalysis}

To evaluate the adaptability of the proposed multi-stage stochastic programming framework, we perform a set of sensitivity analyses under four distinct cases: i) \texttt{Relaxed Budget}, ii) \texttt{Relaxed Emission}, iii) \texttt{Safety Margin}, and iv) \texttt{Low Turbine Price}. These cases are designed to evaluate the effects of increased financial flexibility, modified emission policies, enhanced system reliability constraints, and policy-induced cost reductions on optimal investment and operational decisions. Summary statistics of expected cost components and solution times are reported in Table~\ref{tab:optimizationresults-sens}. Comprehensive results—including cumulative electricity generation and cost components for each node, as well as annual electricity supply mixes and cost breakdowns for extreme scenario paths-are provided in \ref{sec:appendixA}.

\begin{table}[h]
    \centering
    \caption{Expected cost components, optimality gap (\%), and solution time (seconds) for each sensitivity analysis case. }
    \begin{adjustbox}{max width=\textwidth}
    \label{tab:optimizationresults-sens}
    \begin{tabular}{|c|c|c|c|c|c|c|}
    \hline
    \textbf{Case} & \textbf{Total Cost} & \textbf{Installation Cost} & \textbf{Grid Purchase Cost} & \textbf{O\&M Cost} & \textbf{Optimality Gap} & \textbf{Time} \\ \hline
    \texttt{Base} & \$69.3 M & \$47.1 M & \$20.0 M & \$2.2 M & 0.92\% & 6706 s \\ \hline
    \texttt{Relaxed Budget}   & \$65.4 M & \$45.2 M & \$17.5 M & \$2.7 M & 0.72\% & 39785 s \\ \hline
    \texttt{Relaxed Emission}  & \$55.0 M & \$29.2 M & \$24.3 M & \$1.6 M & 0.97\% & 9384 s \\ \hline
    \texttt{Safety Margin} & \$76.6 M & \$53.8 M & \$20.3 M & \$2.6 M & 0.99\% & 39809 s \\ \hline
    \texttt{Low Turbine Price} & \$67.0 M & \$46.0 M & \$18.6 M & \$2.5 M & 0.81\% & 41565 s \\ \hline
    \end{tabular}
    \end{adjustbox}
\end{table}

In the \texttt{Relaxed Budget Case}, we increase the installation budget to \$15 million for all periods. This adjustment reduces the objective value to \$65.4 million, primarily due to the model’s ability to allocate investments earlier in the planning horizon, thereby extending the utilization of renewable electricity sources and decreasing reliance on grid procurement. The expanded budget consistently lowers grid usage across all scenario paths, with the expected cost of grid purchases declining by \$2.5 million.  Moreover, the additional budget enables the installation of 12-\si{\mega\watt} PV and 6 \si{\mega\watt} wind capacity in the first year, (see Figure~\ref{fig:decision_tree_relaxed_budget}), compared to the \texttt{Base Case}, in which 12 \si{\mega\watt} of PV capacity is distributed across the first two years. 
In the second-stage nodes, technology deployment occurs exclusively in the sixth year, reflecting economically-driven decisions. This outcome arises because the less stringent budget constraints provide sufficient flexibility for environmentally-driven investments to be deferred to the third-stage nodes. Although the improvements are moderate, this case highlights the advantages of enhanced financial flexibility in accelerating the transition toward clean energy.

In the \texttt{Relaxed Emission Case}, we allow up to 1\% of university’s electricity demand to be met from grid sources during the final period to examine the trade-off between strict decarbonization targets and economic viability.
The resulting objective function value decreases to \$55.0 million, underscoring the cost of stringent clean energy targets. 
In the \texttt{Base Case}, the requirement to fully align generation and demand profiles leads to substantial investments in technologies whose capacities remain underutilized. By permitting limited grid usage, the model reduces this investment burden, particularly in later periods.
In particular, by increasing the expected grid purchase cost by \$4.3 million, the model is able to reduce the expected installation cost by \$14.3 million compared to the \texttt{Base Case}. 
Notably, most reductions in deployment occur in battery technologies across scenario paths (ranging between 28.8-58.5\%), reflecting a shift away from overcapacity in storage technologies (see Figure~\ref{fig:decision_tree_relaxed_emission}).

Figure~\ref{fig:emission_allowance_cost_breakdown} further illustrates the cost breakdown under varying emission allowance levels. The results indicate that the total expected costs, installation costs, and O\&M costs decline as the emission allowance increases, whereas the expected grid purchase costs rise and eventually surpasses the expected installation costs. Furthermore, the marginal reduction in expected total cost diminishes with higher allowance levels (e.g., a 20.5\% reduction from 0\% to 1\% allowance, compared to a 5.6\% reduction from 1\% to 2\%), thereby underscoring the cost implications of pursuing more stringent clean energy targets. Notably, the total expected cost results with different emission allowance levels, 1\% to 5\% allowance results in \$55.0, \$52.0, \$50.1, \$48.9, \$48.1 million, is less than cost of fullfilling demand via grid electricity procurements (\$58.8 million). This result also shows the economic viability of transitioning to clean energy technologies.

\vspace{-10pt}

\begin{figure}[H]
\centering
\begin{tikzpicture}[scale=1]
    \pgfplotsset{width=0.6\textwidth, height=5cm}
    \begin{axis}[
      axis lines=box,
      xtick pos=bottom,
      enlarge x limits=0.05,
      ytick pos=left,
      ytick={0,14,28,42,56,70},
      ymin=0, ymax=70,
      xlabel={Emission Allowance (\%)},
      xmin = 0,
      xmax = 5,
      xtick={0,...,5},
      ylabel={Cost (million USD, 2024)},
      legend style={at={(1,1)}, anchor=north east, legend columns=2, nodes={scale=0.6, transform shape}, font=\large},
      y axis line style=black!75!black,
      ylabel style={xshift=0.01cm, font=\fontsize{9pt}{10pt}\selectfont},
      xlabel style={yshift=0.1cm, font=\fontsize{9pt}{10pt}\selectfont},
      xticklabel style={yshift=-2pt},   
      tick label style={font=\footnotesize},
      ytick style=black!75!black,
      xtick style=black!75!black,
      yticklabel style={/pgf/number format/.cd, fixed},
      legend reversed,
    ]
\addplot[mark=triangle*, orange!50!black, mark size=2pt] coordinates {
        (0,2.181177969)
        (1,1.567574673)
        (2,1.55023887)
        (3,1.521145897)
        (4,1.480827477)
        (5,1.455782851)
    };
    \addlegendentry{O\&M Cost}
\addplot[mark=square*, orange, mark size=1.5pt] coordinates {
        (0,19.96541383)
        (1,24.31075922)
        (2,24.39657599)
        (3,24.0167952)
        (4,24.90004453)
        (5,24.53821147)
    };
    \addlegendentry{Grid Purchase Cost}
\addplot[mark=*, green!50!black, mark size=1.8pt] coordinates {
        (0,47.12831367)
        (1,29.16233972)
        (2,26.03764907)
        (3,24.54635533)
        (4,22.56890686)
        (5,22.08374921)
    };
    \addlegendentry{Installation Cost}
\addplot[mark=diamond*, red, mark size=2pt] coordinates {
        (0,69.2643004)
        (1,55.04067361)
        (2,51.98446392)
        (3,50.08429642)
        (4,48.94977886)
        (5,48.07774353)
    };
    \addlegendentry{Total Cost}
    \end{axis}
\end{tikzpicture}
\caption{Cost breakdown under the \texttt{Relaxed Emission Case} for different emission allowance levels.}
\label{fig:emission_allowance_cost_breakdown}
\end{figure}

In the Safety Margin Case, a buffer is introduced for each sub-period to explicitly account for uncertainties arising from high-resolution stochastic inputs related to electricity demand and renewable generation. Accordingly, the demand allocated to each sub-period is increased by 8\%. Beyond addressing uncertainty, this adjustment also mitigates potential mismatches between supply and demand within individual sub-periods. Under this case, the objective value increases to \$76.6 million, representing a moderate additional cost for enhanced system reliability. Across all scenario paths, battery installations rise relative to the base case. In 13 of 16 scenario paths, solar installations also increase, ranging from 6.7\% to 20\%. In the two scenario paths where solar capacity remains unchanged relative to the \texttt{Base Case} ($S_{1}$ and $S_{14}$), significant increments in battery installations are observed, amounting to 11.3\% and 17.9\%, respectively. In the only scenario path where solar capacity decreases relative to the \texttt{Base Case} ($S_{13}$), an additional 6-\si{\mega\watt} wind turbine is installed in year 11, alongside a 9.4\% increase in battery capacity. Moreover, in the tenth year of node 4, a 36-\si{\mega\watt} solar installation is observed, indicating that higher demand partially shifts environmentally driven decisions from descendant nodes back to the parent node. In the first stage, an 8.1-\si{\mega\watt\hour} battery installation is introduced, which is not present in the base case. These findings underscore the importance of incorporating a safety margin to address potential mismatches and renewable generation variability. However, they also highlight that such measures impose significant costs, thereby influencing long-term strategic investment decisions.


In the \texttt{Low Turbine Price Case}, wind turbine installation costs are assumed to be 80\% of their baseline level, reflecting a potential policy-driven incentive. Relative to the \texttt{Base Case}, this reduction promotes earlier deployment of wind capacity, with a 6-\si{\mega\watt} turbine installed in year 2 (see Figure~\ref{fig:decision_tree_low_turbine_price}), as wind energy becomes more cost-competitive. The increase in wind-based electricity generation ranges from 17.6\% to 819.2\%, with a median of 78.9\%, particularly along scenario paths featuring rapid improvements in solar PV and battery technologies during the first-stage transition. Across all scenario paths, total system costs, installation costs, and grid purchase costs decline, whereas O\&M costs rise relative to the \texttt{Base Case}. Moreover, the share of emissions in the total supply mix decreases to below 8\% by year 11, compared to 12\% in the \texttt{Base Case}. These findings highlight wind power’s complementary role to solar PV in enhancing generation portfolio diversity. Nevertheless, even under favorable cost assumptions, solar technologies remain the more economically attractive option for the campus site, where wind potential is constrained by relatively low capacity factors.

%


\subsection{Temporally Aggregated Models}
\label{sec:aggregatemodelcase}

To investigate the effects of temporal resolution on strategic-level clean energy planning, we develop alternative models with varying degrees of temporal aggregation—specifically 24-hour, 12-hour, 8-hour, and 4-hour intervals. The solutions of these models are then evaluated against our bihourly framework to assess the impacts of temporal aggregation on model outcomes. Table~\ref{tab:aggregatedmodelresults} summarizes the expected cost components, average electricity demand violations, and the average cost of fulfilling unmet demand through grid procurement for each level of aggregation. Notably, the 24-hour, 12-hour, 8-hour, and 4-hour models experience demand shortfalls corresponding to 24.2\%, 20.3\%, 0.9\%, and 0.8\% of total demand, respectively, when averaged across scenario paths.

As temporal aggregation increases, the model’s objective value decreases, while demand constraint violations become more pronounced. This inverse relationship arises from the reduced visibility of short-term fluctuations in electricity demand and renewable generation profiles, which limits the model’s ability to accurately align supply with demand. High-resolution models capture variability more effectively, thereby motivating substantial investments in both generation and storage technologies. In contrast, aggregated models obscure these dynamics, reducing the perceived need for storage capacity and additional generation. In particular, battery storage investments decrease substantially as temporal resolution is reduced, underscoring the critical role of storage technologies in mitigating renewable intermittency and ensuring demand–supply alignment. Furthermore, an increase in aggregation leads to a marked reduction in wind turbine investments. This outcome reflects the relatively weak competitiveness of wind power compared to solar in the selected installation site, as the diminished temporal resolution conceals the complementary role of wind observed in the \texttt{Base Case}.

Temporal aggregation not only compromises the fulfillment of electricity demand but also leads to violations of emission restriction constraints. Even when accounting for grid procurements and incorporating the associated average costs of unmet demand, the total system cost reported by aggregated models remains artificially lower than the actual cost of transitioning to a clean electricity system. This discrepancy arises because supplying even the final, relatively small portion of electricity demand exclusively with clean technologies entails disproportionately high costs, as demonstrated in the \textit{Emission Allowance Case} sensitivity analysis. These findings highlight the importance of employing high temporal resolution models, despite their higher computational burden, to ensure realistic and reliable strategic-level clean energy planning decisions.



 \begin{table}[h]
     \centering
     \caption{Objective values, demand constraint violations, cost of fulfilling unmet demand via grid procurements, optimality gap (\%), and solution time (seconds) for temporally aggregated models simulated under the base model with bihourly resolution.}
     \label{tab:aggregatedmodelresults}
     \scalebox{0.93}{
     \begin{tabular}{|c|c|c|c|c|c|}
         \hline
         \makecell[c]{Operational\\Resolution}
         & \makecell[c]{Objective\\Value}
         & \makecell[c]{Average Demand\\Violation}
         & \makecell[c]{Average Cost of \\Demand Violation} & Optimality Gap & Time \\ \hline
         {24-hour} & \$52.5 M & \SI{124.92}{\giga\watt\hour}& \$15.10 M & 0.86\% & 284 s \\ \hline
         {12-hour} & \$54.7 M & \SI{105.00}{\giga\watt\hour} & \$12.76 M & 1.00\% & 86407 s \\ \hline
         {8-hour} & \$67.5 M & \SI{4.48}{\giga\watt\hour} & \$0.54 M & 0.75\% & 2635 s \\ \hline
         {4-hour} & \$68.0 M & \SI{3.98}{\giga\watt\hour} & \$0.48 M & 0.79\% & 6325 s \\ \hline
     \end{tabular}
    }
\end{table}


\subsection{Deterministic Model}

%

To assess the importance of incorporating uncertainties in technological advancements, we construct a deterministic model by assigning a single representative cluster center to each technology. This deterministic model yields an objective value of \$67.1 million, with the solution obtained in 483 seconds. By utilizing a single deterministic trajectory to represent technological advancements, this approach models average improvements in cost and efficiency. However, in practice, realized scenarios may involve slower cost reductions and more limited efficiency gains. Consequently, when applied to such realizations, the deterministic solution may fail to satisfy budget, demand, and emission constraints, thereby rendering it infeasible within the multi-stage stochastic programming framework.

Table~\ref{tab:deterministicmodelinstallations} presents the annual installation decisions for solar PV, wind power, and battery storage derived from the deterministic model. Unlike the stochastic model, only a 6-\si{\mega\watt} solar PV installation occurs in the first stage. Battery investments are postponed, with most capacity added in the final five years. Moderate amounts of wind power are installed in years six and seven, whereas the stochastic model exhibits delayed wind power deployment across most scenario paths. This divergence arises from the deterministic model’s exclusion of extreme technological advancement trajectories that would otherwise strongly favor solar or battery technologies, leading instead to more balanced investments across technologies. Notably, the deterministic model generates economically-driven decisions in the second stage, resulting in unrealistic outcomes under slow technological improvement scenarios, where last stage deployments alone are insufficient to achieve a full transition to clean electricity. These findings highlight the deterministic model’s limited ability to adapt to actual realizations of technological progress.

\newcolumntype{C}[1]{>{\centering\arraybackslash}m{#1}}

\begin{table}[htbp]
\centering
\caption{Installation decisions in the deterministic model.}
\label{tab:deterministicmodelinstallations}
\setlength{\tabcolsep}{2pt}
\renewcommand{\arraystretch}{1.1}
\scalebox{0.84}{
\begin{tabular}{|C{4cm}|*{15}{C{0.85cm}|}}
\hline
\diagbox[width=4cm, height=1cm]{Technology}{Year} & 1 & 2 & 3 & 4 & 5 & 6 & 7 & 8 & 9 & 10 & 11 & 12 & 13 & 14 & 15 \\ \hline
\rowcolor{solaryellow}
Solar (\si{\mega\watt}) & 6 & - & - & - & - & 12 & - & - & - & - & 24 & - & - & 12 & 36 \\ \hline
\rowcolor{windcyan}
Wind (\si{\mega\watt}) & - & - & - & - & - & 6 & 6 & - & - & - & - & - & - & - & - \\ \hline
\rowcolor{batterypurple}
Battery (\si{\mega\watt\hour}) & - & - & - & - & - & 9.9 & 9.2 & - & - & - & 55.4 & 0.8 & 103.0 & 85.0 & 32.0 \\ \hline
\end{tabular}}
\end{table}

Tables~\ref{tab:deterministic_budget_violation} and~\ref{tab:deterministic_demand_violation} summarize the extent of budget and demand violations that arise when applying the deterministic solution across different scenario paths, thereby underscoring the limitations of the deterministic approach compared with the more flexible stochastic programming framework. For example, in scenario $S_1$, characterized by slow technological progress throughout the planning horizon, the deterministic solution would require more than an additional 20.7 million USD to implement its prescribed strategic plan.

Slower-than-anticipated improvements in technological efficiency lead to reduced electricity generation, resulting in demand constraint violations in 12 of the 16 scenario paths. In these cases, minor violations occur in the final year of the planning horizon due to insufficient renewable electricity supply. This shortfall must be compensated through conventional grid-based generation, thereby preventing a complete transition to clean electricity. The results are further misleading because, as discussed in Section~\ref{sec:compexBase}, addressing this residual demand with renewable technologies in the final year would be prohibitively expensive.

The deterministic model not only risks infeasibility under scenarios of slow technological advancement but also tends to induce overinvestment and excessive grid electricity procurement in cases of rapid technological progress. These inefficiencies highlight the deterministic model’s inability to dynamically adapt to diverse and evolving technological conditions.

\section{Conclusion}
\label{sec:conclusion}

This study develops a comprehensive multi-stage stochastic programming framework for long-term clean energy transition planning, with particular emphasis on technological uncertainty and high-resolution temporal dynamics. By explicitly modeling uncertainties in both the costs and efficiencies of renewable technologies and embedding detailed operational constraints, the framework effectively links strategic investment planning with short-term system operations. Over the planning horizon, the model jointly optimizes investment and operational decisions, thereby enabling large-scale consumers to pursue ambitious decarbonization targets in a cost-effective and operationally feasible manner. A central methodological contribution lies in aligning strategic decision-making with both evolving technological developments and the temporal variability of electricity demand and generation.

The case study of METU demonstrates the applicability of the proposed framework in a university campus setting. Our results show that solar photovoltaic technology is the most economically advantageous option, reflecting METU’s strong solar potential and expected cost reductions in PV systems. By contrast, the relatively low capacity factor of wind turbines at the study site constrains their role to a complementary resource. These findings highlight the strong sensitivity of transition strategies to site-specific conditions, emphasizing the importance of localized data in energy planning. At the same time, the analysis yields generalizable insights: battery storage is indispensable for addressing intermittency between renewable generation and demand. Substantial lithium-ion battery investments are therefore required to ensure system reliability and to cover demand during periods of low renewable output.

Sensitivity analyses further enrich these insights. A clear trade-off emerges between economically driven and environmentally driven decisions. In trajectories of rapid technological progress, the model prioritizes early, cost-effective deployments, whereas under slower progress it defers investments and makes more environmentally motivated choices toward the end of the horizon. The relaxed emission case illustrates this trade-off: allowing 1\% of demand to be met by grid electricity reduces total expected cost by 20.5\%, underscoring the economic burden associated with stringent zero-emission targets.



The results also demonstrate the importance of incorporating both stochastic technological advancements and high-resolution temporal detail. Comparison with a deterministic model reveals that ignoring technological uncertainty may lead to demand shortfall and severe budget overruns under slow technological advancement trajectories. Similarly, coarser temporal resolutions underestimate storage requirements, leading to significant demand violations. These findings underscore that high temporal resolution and explicit stochastic modeling are critical for producing transition strategies that are both economically viable and practically feasible.

The framework is designed to be adaptable to diverse energy systems and consumer categories, including university campuses, industrial facilities, and institutional complexes. Its flexible structure supports customization across different geographic regions, system sizes, and technology portfolios. In contexts where planners must balance long-term sustainability goals with short-term operational reliability under uncertain technological pathways, the proposed model offers a systematic and adaptive decision-support tool. Although this study focuses on electricity transitions at the campus scale, the framework can be extended to integrate multiple energy vectors and sustainable technologies.
Particularly, there are two specific future work directions that are promising. Firstly, we plan to develop a model that jointly optimizes both electricity and thermal energy transition for the METU campus. Secondly, we aim to explore how our model can be used for policy-making in order to design mechanisms for renewable generation incentives, especially in the case of large energy consumers. As these future directions are pursued, the complexity of the problem is expected to increase, and we will therefore focus on developing and applying more advanced modeling and solution techniques (such as decomposition algorithms) tailored to our context.

\subsection*{Acknowledgements}
This work was supported by the  Scientific and Technological Research Council of Turkey [grant number 222M243]. The authors thank METU for providing detailed datasets, Ilg{\i}n Savur for reporting them, and Neman Karimi  for his help related to implementation.

\bibliography{mybibfile}

\begin{thebibliography}{10}
\expandafter\ifx\csname url\endcsname\relax
  \def\url#1{\texttt{#1}}\fi
\expandafter\ifx\csname urlprefix\endcsname\relax\def\urlprefix{URL }\fi
\expandafter\ifx\csname href\endcsname\relax
  \def\href#1#2{#2} \def\path#1{#1}\fi

\bibitem{McGrath_Hunt_Taylor_Fitzgerald_2025}
D.~McGrath, H.~Hunt, J.~R. Taylor, S.~Fitzgerald, Meeting climate goals through mitigation and intervention: developments in emissions reduction, greenhouse gas removal, and solar radiation modification, Global Sustainability 8 (2025) e20.

\bibitem{arora2024cop28}
P.~Arora, Cop28: ambitions, realities, and future, Environmental Sustainability 7 (2024) 107--113.

\bibitem{METUPolicies}
{Middle East Technical University}, \href{https://sustainablecampus2.metu.edu.tr/en/metu-climate-action-and-sustainability-plan}{Metu climate action and sustainability plan}, accessed Sep 20, 2025 (2022).
\newline\urlprefix\url{https://sustainablecampus2.metu.edu.tr/en/metu-climate-action-and-sustainability-plan}

\bibitem{ling2024}
J.-M. Ling, F.~I. Mulani, Planning and optimization of a residential microgrid utilizing renewable resources and integrated energy storage, Journal of Energy Storage 97 (2024) 112933.

\bibitem{Guevara2022}
E.~Guevara, F.~Babonneau, T.~H. de~Mello, Uncertainty dynamics in energy planning models: An autoregressive and markov chain modeling approach, Computers \& Industrial Engineering 191 (2024) 110084.

\bibitem{Lei2021}
Y.~Lei, D.~Wang, H.~Jia, J.~Li, J.~Chen, J.~Li, Z.~Yang, Multi-stage stochastic planning of regional integrated energy system based on scenario tree path optimization under long-term multiple uncertainties, Applied Energy 300 (2021) 117224.

\bibitem{Cano2016}
E.~L. Cano, J.~M. Moguerza, A.~Alonso-Ayuso, A multi-stage stochastic optimization model for energy systems planning and risk management, Energy and Buildings 110 (2016) 49--56.

\bibitem{Ioannou2019}
A.~Ioannou, G.~Fuzuli, F.~Brennan, S.~W. Yudha, A.~Angus, Multi-stage stochastic optimization framework for power generation system planning integrating hybrid uncertainty modelling, Energy Economics 80 (2019) 760--776.

\bibitem{rathi2022}
T.~Rathi, Q.~Zhang, Capacity planning with uncertain endogenous technology learning, Computers \& Chemical Engineering 164 (2022) 107868.

\bibitem{gao2023}
N.~Gao, D.~Gao, X.~Fang, Manage real-time power imbalance with renewable energy: Fast generation dispatch or adaptive frequency regulation?, Power Systems, IEEE Transactions on (12 2022).

\bibitem{Marcy2022}
C.~Marcy, T.~Goforth, D.~Nock, M.~Brown, Comparison of temporal resolution selection approaches in energy systems models, Energy 251 (2022) 123969.

\bibitem{Zou2018}
J.~Zou, S.~Ahmed, X.~Sun, Partially adaptive stochastic optimization for electric power generation expansion planning, INFORMS Journal on Computing 30 (2018) 388--401.

\bibitem{park2020}
H.~Park, Generation capacity expansion planning considering hourly dynamics of renewable resources, Energies 13~(21) (2020).

\bibitem{Singh2009}
K.~J. Singh, A.~B. Philpott, R.~K. Wood, Dantzig-wolfe decomposition for solving multistage stochastic capacity-planning problems, Operations Research 57~(5) (2009) 1271--1286.

\bibitem{Rebennack2016}
S.~Rebennack, Combining sampling-based and scenario-based nested benders decomposition methods: application to stochastic dual dynamic programming, Math. Program. 156~(1–2) (2016) 343–389.

\bibitem{Hole2023}
J.~Hole, A.~Philpott, O.~Dowson, Capacity planning of renewable energy systems using stochastic dual dynamic programming, European Journal of Operational Research 322~(2) (2025) 573--588.

\bibitem{Liu2018}
Y.~Liu, R.~Sioshansi, A.~J. Conejo, Multistage stochastic investment planning with multiscale representation of uncertainties and decisions, IEEE Transactions on Power Systems 33~(1) (2018) 781--791.

\bibitem{Husein2018}
M.~Husein, I.-Y. Chung, Optimal design and financial feasibility of a university campus microgrid considering renewable energy incentives, Applied Energy 225 (2018) 273--289.

\bibitem{Prina2019}
M.~G. Prina, M.~Lionetti, G.~Manzolini, W.~Sparber, D.~Moser, Transition pathways optimization methodology through energyplan software for long-term energy planning, Applied Energy 235 (2019) 356--368.

\bibitem{Gils2017}
H.~C. Gils, Y.~Scholz, T.~Pregger, D.~{Luca de Tena}, D.~Heide, Integrated modelling of variable renewable energy-based power supply in europe, Energy 123 (2017) 173--188.

\bibitem{Moret2020}
S.~Moret, F.~Babonneau, M.~Bierlaire, F.~Mar{\'e}chal, Decision support for strategic energy planning: A robust optimization framework, European Journal of Operational Research 280~(2) (2020) 539--554.

\bibitem{Heuberger2017}
C.~F. Heuberger, E.~S. Rubin, I.~Staffell, N.~Shah, N.~{Mac Dowell}, Power capacity expansion planning considering endogenous technology cost learning, Applied Energy 204 (2017) 831--845.

\bibitem{Powell2012}
W.~Powell, A.~George, H.~Simão, W.~Scott, A.~Lamont, J.~Stewart, Smart: A stochastic multiscale model for the analysis of energy resources, technology, and policy, INFORMS Journal on Computing 24 (2012) 665--682.

\bibitem{Child2019}
M.~Child, C.~Kemfert, D.~Bogdanov, C.~Breyer, Flexible electricity generation, grid exchange and storage for the transition to a 100\% renewable energy system in europe, Renewable Energy 139 (2019) 80--101.

\bibitem{Mavromatidis2021}
G.~Mavromatidis, I.~Petkov, Mango: A novel optimization model for the long-term, multi-stage planning of decentralized multi-energy systems, Applied Energy 288 (2021) 116585.

\bibitem{Tian2022}
X.~Tian, Y.~Zhou, B.~Morris, F.~You, Sustainable design of cornell university campus energy systems toward climate neutrality and 100\% renewables, Renewable and Sustainable Energy Reviews 161 (2022) 112383.

\bibitem{Zhao2021}
N.~Zhao, F.~You, New york state's 100\% renewable electricity transition planning under uncertainty using a data-driven multistage adaptive robust optimization approach with machine-learning, Advances in Applied Energy 2 (2021) 100019.

\bibitem{Abulibdeh2024}
A.~Abulibdeh, Towards zero-carbon, resilient, and community-integrated smart schools and campuses: A review, World Development Sustainability 5 (2024) 100193.

\bibitem{Li2022}
L.~Bin, M.~Shahzad, H.~Javed, H.~A. Muqeet, M.~N. Akhter, R.~Liaqat, M.~M. Hussain, Scheduling and sizing of campus microgrid considering demand response and economic analysis, Sensors 22~(16) (2022).

\bibitem{Paspatis2022}
A.~Paspatis, K.~Fiorentzis, Y.~Katsigiannis, E.~Karapidakis, Smart campus microgrids towards a sustainable energy transition—the case study of the hellenic mediterranean university in crete, Mathematics 10 (2022) 1065.

\bibitem{AlKassem2022}
A.~AlKassem, A.~Draou, A.~Alamri, H.~Alharbi, Design analysis of an optimal microgrid system for the integration of renewable energy sources at a university campus, Sustainability 14~(7) (2022).

\bibitem{Mishra2025}
J.~Mishra, A.~Shankar, Optimizing wind-pv-battery microgrids for sustainable and resilient residential communities, Scientific Reports 15 (07 2025).

\bibitem{sam2025}
{\relax National Renewable Energy Laboratory}, System advisor model™, version 2025.4.16 (2025).

\bibitem{karimi2024}
N.~Karimi, B.~Kocuk, T.~Yuksel, \href{https://arxiv.org/abs/2410.23387}{{A Dynamic Strategic Plan for the Transition to a Clean Bus Fleet using Multi-Stage Stochastic Programming with a Case Study in Istanbul}} (2024).
\newblock \href {http://arxiv.org/abs/2410.23387} {\path{arXiv:2410.23387}}.
\newline\urlprefix\url{https://arxiv.org/abs/2410.23387}

\bibitem{campusapplicationrep}
A.~E. Sener, T.~Yuksel, B.~Kocuk, \href{https://github.com/emirseners/Campus-Application.git}{Campus application}, accessed Sep 20, 2025 (2025).
\newline\urlprefix\url{https://github.com/emirseners/Campus-Application.git}

\bibitem{Moore1998}
G.~Moore, Cramming more components onto integrated circuits, Proceedings of the IEEE 86~(1) (1998) 82--85.

\bibitem{tracking_sun_2019}
G.~L. Barbose, N.~R. Darghouth, Tracking the sun: Pricing and design trends for distributed photovoltaic systems in the united states - 2019 edition, Tech. rep., Lawrence Berkeley National Laboratory (2019).

\bibitem{Jordan2011}
D.~Jordan, S.~Kurtz, Photovoltaic degradation rates—an analytical review, Progress in Photovoltaics: Research and Applications 21 (01 2013).

\bibitem{LBL2023}
M.~Bolinger, J.~Seel, J.~M. Kemp, C.~Warner, A.~Katta, D.~Robson, Utility-scale solar, 2023 edition: Empirical trends in deployment, technology, cost, performance, ppa pricing, and value in the united states, Tech. rep., Lawrence Berkeley National Laboratory (2023).

\bibitem{NSRDB2018}
M.~Sengupta, Y.~Xie, A.~Lopez, A.~Habte, G.~Maclaurin, J.~Shelby, The national solar radiation data base (nsrdb), Renewable and Sustainable Energy Reviews 89 (2018) 51--60.

\bibitem{Bergey2012}
K.~H. Bergey, The lanchester-betz limit (energy conversion efficiency factor for windmills), Journal of Energy 3~(6) (1979) 382--384.

\bibitem{Rosenberg2014}
A.~Rosenberg, S.~Selvaraj, A.~Sharma, A novel dual-rotor turbine for increased wind energy capture, Journal of Physics: Conference Series 524 (06 2014).

\bibitem{thewindpowerDB}
M.~P. EI, \href{https://www.thewindpower.net/index.php}{The wind power}, accessed Sep 20, 2025 (2025).
\newline\urlprefix\url{https://www.thewindpower.net/index.php}

\bibitem{Molteno2022Biomimetics}
T.~C.~A. Molteno, Nature’s wind turbines: The measured aerodynamic efficiency of spinning seeds approaches theoretical limits, Biomimetics 7~(4) (2022).

\bibitem{irena2023renewable}
IRENA, Renewable power generation costs in 2022, Tech. rep., International Renewable Energy Agency (2023).

\bibitem{DOE2023}
R.~Wiser, M.~Bolinger, B.~Hoen, D.~Millstein, J.~Rand, G.~Barbose, N.~Darghouth, W.~Gorman, S.~Jeong, E.~O'Shaughnessy, , B.~Paulos, Land-based wind market report: 2023 edition, Tech. rep., U.S. Department of Energy (2023).

\bibitem{Staffell2014}
I.~Staffell, R.~Green, How does wind farm performance decline with age?, Renewable Energy 66 (2014) 775--786.

\bibitem{windplantREGULATION}
{T.C. Resmî Gazete [Official Gazette of the Republic of Turkey]}, {\relax Rüzgar Kaynağına Dayalı Elektrik Üretimi Başvurularının Teknik Değerlendirmesi Hakkında Yönetmelik [Regulation on the Technical Evaluation of Applications for Electricity Generation Based on Wind Resource]}, \url{https://www.mevzuat.gov.tr/mevzuat?MevzuatNo=21189&MevzuatTur=7&MevzuatTertip=5}, accessed Sep 20, 2025 (2015).

\bibitem{NASAReference}
\href{https://power.larc.nasa.gov/}{National aeronautics and space administration (nasa) langley research center's prediction of worldwide energy resources (power)}, accessed Sep 20, 2025 (2021).
\newline\urlprefix\url{https://power.larc.nasa.gov/}

\bibitem{Firtin2011}
E.~Fırtın, Önder Güler, S.~A. Akdağ, Investigation of wind shear coefficients and their effect on electrical energy generation, Applied Energy 88~(11) (2011) 4097--4105.

\bibitem{globalwindatlas}
\href{https://globalwindatlas.info/en/}{Global wind atlas}, accessed Sep 20, 2025 (2015).
\newline\urlprefix\url{https://globalwindatlas.info/en/}

\bibitem{walter2023}
D.~Walter, K.~Bond, S.~Butler-Sloss, L.~Speelman, Y.~Numata, W.~Atkinson, X-change: Batteries — the battery domino effect, Tech. rep., Rocky Mountain Institute (RMI) (2023).

\bibitem{autometris}
A.~Rohatgi, \href{https://automeris.io}{Webplotdigitizer}, version 5.2 (2024).
\newline\urlprefix\url{https://automeris.io}

\bibitem{sungrowcatalogue}
Sungrow Power Supply Co., Ltd., Hefei, China, \href{https://vas-co.com/wp-content/uploads/2024/05/EN-BR-Sungrow-Energy-Storage-System-Products-Catalogue.pdf}{Energy Storage System Products Catalogue}, europe edition; Version 1.1 / 2021--2022 (2021).
\newline\urlprefix\url{https://vas-co.com/wp-content/uploads/2024/05/EN-BR-Sungrow-Energy-Storage-System-Products-Catalogue.pdf}

\bibitem{ZHAO2024110398}
Y.~Zhao, Y.~Zhang, Y.~Li, Y.~Chen, W.~Huo, H.~Zhao, Optimal configuration of energy storage for alleviating transmission congestion in renewable energy enrichment region, Journal of Energy Storage 82 (2024) 110398.

\bibitem{Cole2025NREL}
W.~Cole, V.~Ramasamy, M.~Turan, Cost projections for utility-scale battery storage: 2025 update, Tech. rep., National Renewable Energy Laboratory (2025).

\bibitem{Mohamed2021}
A.~A.~R. Mohamed, R.~J. Best, X.~Liu, D.~J. Morrow, Residential battery energy storage sizing and profitability in the presence of pv and ev, in: 2021 IEEE Madrid PowerTech, IEEE, 2021, pp. 1--6.

\bibitem{Tejero2024}
J.~A. Tejero-G{\'o}mez, {\'A}.~A. Bayod-R{\'u}jula, Sizing of battery energy storage systems for firming pv power including aging analysis, Energies 17~(6) (2024).

\bibitem{reber2023beyond}
D.~Reber, S.~R. Jarvis, M.~P. Marshak, Beyond energy density: flow battery design driven by safety and location, Energy Advances 2~(9) (2023) 1357--1365.

\bibitem{epdk}
{\relax Enerji Piyasası Düzenleme Kurumu [Energy Market Regulatory Authority]}, Elektrik faturalar{\i}na esas tarife tablolar{\i} [tariff tables underpinning electricity bills], \url{https://www.epdk.gov.tr/Detay/Icerik/3-1327/elektrik-faturalarina-esas-tarife-tablolari}, accessed Sep 23, 2025 (2025).

\bibitem{shapiro2021}
A.~Shapiro, D.~Dentcheva, A.~Ruszczynski, Lectures on stochastic programming: modeling and theory, SIAM, 2021.

\end{thebibliography}

\appendix

\clearpage
\section{Deterministic Model Constraint Violations}

 \begin{table}[H]
  \centering
  \caption{Budget constraint violations (thousand USD, 2024) of the solution given by the deterministic model under different scenarios.}
  \label{tab:deterministic_budget_violation}
  \begin{adjustbox}{width=0.98\textwidth}
\begin{tabular}{lrrrrrrrrrrr}
    \toprule
    & Y6 & Y7 & Y8 & Y9 & Y10 & Y11 & Y12 & Y13 & Y14 & Y15 & Total \\
    \midrule
    $S_{1}$: \(ss \smalltimes ss\) & 1554 & - & - & - & - & 5169 & - & 4121 & 4947 & 4966 & 20757 \\
    $S_{2}$: \(ss \smalltimes sf\) & 1554 & - & - & - & - & 2371 & - & - & 1032 & 3534 & 8492 \\
    $S_{3}$: \(ss \smalltimes fs\) & 1554 & - & - & - & - & 2989 & - & 4121 & 3952 & 2071 & 14687 \\
    $S_{4}$: \(ss \smalltimes ff\) & 1554 & - & - & - & - & 191 & - & - & 37 & 639 & 2422 \\
    $S_{5}$: \(sf \smalltimes ss\) & 635 & - & - & - & - & 2371 & - & - & 1032 & 3534 & 7573 \\
    $S_{6}$: \(sf \smalltimes sf\) & 635 & - & - & - & - & 835 & - & - & - & 2748 & 4219 \\
    $S_{7}$: \(sf \smalltimes fs\) & 635 & - & - & - & - & 191 & - & - & 37 & 639 & 1503 \\
    $S_{8}$: \(sf \smalltimes ff\) & 635 & - & - & - & - & - & - & - & - & - & 635 \\
    $S_{9}$: \(fs \smalltimes ss\) & 72 & - & - & - & - & 2989 & - & 4121 & 3952 & 2071 & 13204 \\
    $S_{10}$: \(fs \smalltimes sf\) & 72 & - & - & - & - & 191 & - & - & 37 & 639 & 939 \\
    $S_{11}$: \(fs \smalltimes fs\) & 72 & - & - & - & - & 1586 & - & 4121 & 3312 & 207 & 9298 \\
    $S_{12}$: \(fs \smalltimes ff\) & 72 & - & - & - & - & - & - & - & - & - & 72 \\
    $S_{13}$: \(ff \smalltimes ss\) & - & - & - & - & - & 191 & - & - & 37 & 639 & 867 \\
    \bottomrule
  \end{tabular}
\end{adjustbox}
\end{table}

\vspace{5pt}

 \begin{table}[H]
  \centering
  \caption{Demand constraint violations (\si{\mega\watt\hour}) of the solution given by the deterministic model under different scenarios.}
  \label{tab:deterministic_demand_violation}
  \begin{adjustbox}{width=0.98\textwidth}
  \begin{tabular}{lrrrrrrrrrrr}
    \toprule
    & Y6 & Y7 & Y8 & Y9 & Y10 & Y11 & Y12 & Y13 & Y14 & Y15 & Total \\
    \midrule
    $S_{1}$: \(ss \smalltimes ss\) & 97.5 & 87.7 & 88.0 & 88.6 & 89.5 & 138.0 & 140.8 & 95.2 & 21.7 & 22.6 & 869.5 \\
    $S_{2}$: \(ss \smalltimes sf\) & 97.5 & 87.7 & 88.0 & 88.6 & 89.5 & 138.0 & 140.8 & 95.2 & 21.7 & 22.6 & 869.5 \\
    $S_{3}$: \(ss \smalltimes fs\) & 97.5 & 87.7 & 88.0 & 88.6 & 89.5 & 50.9 & 52.5 & 26.6 & 5.7 & 6.9 & 593.8 \\
    $S_{4}$: \(ss \smalltimes ff\) & 97.5 & 87.7 & 88.0 & 88.6 & 89.5 & 50.9 & 52.5 & 26.6 & 5.7 & 6.9 & 593.8 \\
    $S_{5}$: \(sf \smalltimes ss\) & 97.5 & 87.7 & 88.0 & 88.6 & 89.5 & 138.0 & 140.8 & 95.2 & 21.7 & 22.6 & 869.5 \\
    $S_{6}$: \(sf \smalltimes sf\) & 97.5 & 87.7 & 88.0 & 88.6 & 89.5 & 138.0 & 140.8 & 95.2 & 21.7 & 22.6 & 869.5 \\
    $S_{7}$: \(sf \smalltimes fs\) & 97.5 & 87.7 & 88.0 & 88.6 & 89.5 & 50.9 & 52.5 & 26.6 & 5.7 & 6.9 & 593.8 \\
    $S_{8}$: \(sf \smalltimes ff\) & 97.5 & 87.7 & 88.0 & 88.6 & 89.5 & 50.9 & 52.5 & 26.6 & 5.7 & 6.9 & 593.8 \\
    $S_{9}$: \(fs \smalltimes ss\) & - & - & - & - & - & 14.4 & 14.7 & 7.5 & 2.5 & 4.5 & 43.6 \\
    $S_{10}$: \(fs \smalltimes sf\) & - & - & - & - & - & 14.4 & 14.7 & 7.5 & 2.5 & 4.5 & 43.6 \\
    $S_{13}$: \(ff \smalltimes ss\) & - & - & - & - & - & 14.4 & 14.7 & 7.5 & 2.5 & 4.5 & 43.6 \\
    $S_{14}$: \(ff \smalltimes sf\) & - & - & - & - & - & 14.4 & 14.7 & 7.5 & 2.5 & 4.5 & 43.6 \\
    \bottomrule
  \end{tabular}
\end{adjustbox}
\end{table}

\begin{landscape}

\section{Detailed Results of Sensitivity Analyses}
\label{sec:appendixA}

\vspace*{\fill}
\begin{figure}[H]
\centering
%
};
\draw ([xshift=-1.05cm]1.south) -- ([xshift=1.5cm]2.north) node [midway, left, shift={(-.45,0)}, edge label] {\footnotesize \textbf{\textit{ss}}};
\draw ([xshift=0.5cm]1.south) -- ([xshift=0.6cm]3.north) node [midway, left, shift={(-.05,0)}, edge label] {\footnotesize \textbf{\textit{sf}}};
\draw ([xshift=1.55cm]1.south) -- ([xshift=-0.6cm]4.north) node [midway, right, shift={(.05,0)}, edge label] {\footnotesize \textbf{\textit{fs}}};
\draw ([xshift=3.15cm]1.south) -- ([xshift=-1cm]5.north) node [midway, right, shift={(.45,0)}, edge label] {\footnotesize \textbf{\textit{ff}}};
\draw ([xshift=-0.5cm]2.south) -- ([xshift=0.25cm]6.north) node [midway, left, shift={(-.13,0)}, edge label] {\footnotesize \textbf{\textit{ss}}};
\draw ([xshift=0.2cm]2.south) -- (7.north) node [midway, left, shift={(.02,0)}, edge label] {\footnotesize \textbf{\textit{sf}}};
\draw ([xshift=0.7cm]2.south) -- (8.north) node [midway, right, shift={(-.02,0)}, edge label] {\footnotesize \textbf{\textit{fs}}};
\draw ([xshift=1.5cm]2.south) -- (9.north) node [midway, right, shift={(.13,0)}, edge label] {\footnotesize \textbf{\textit{ff}}};
\draw (3) -- (10.north) node [midway, left, edge label] {\footnotesize \textbf{\textit{ss}}};
\draw (3) -- (11.north) node [midway, left, edge label] {\footnotesize \textbf{\textit{sf}}};
\draw (3) -- (12.north) node [midway, right, edge label] {\footnotesize \textbf{\textit{fs}}};
\draw (3) -- (13.north) node [midway, right, edge label] {\footnotesize \textbf{\textit{ff}}};
\draw (4) -- (14.north) node [midway, left, edge label] {\footnotesize \textbf{\textit{ss}}};
\draw (4) -- (15.north) node [midway, left, edge label] {\footnotesize \textbf{\textit{sf}}};
\draw (4) -- (16.north) node [midway, right, edge label] {\footnotesize \textbf{\textit{fs}}};
\draw (4) -- (17.north) node [midway, right, edge label] {\footnotesize \textbf{\textit{ff}}};
\draw (5) -- (18.north) node [midway, left, shift={(-.13,0)}, edge label] {\footnotesize \textbf{\textit{ss}}};
\draw (5) -- (19.north) node [midway, left, shift={(.05,0)}, edge label] {\footnotesize \textbf{\textit{sf}}};
\draw (5) -- (20.north) node [midway, right, shift={(-.03,0)}, edge label] {\footnotesize \textbf{\textit{fs}}};
\draw (5) -- (21.north) node [midway, right, shift={(.13,0)}, edge label] {\footnotesize \textbf{\textit{ff}}};
\node[above=2pt of 1, xshift=1.05cm] {\small \textbf{N: 1}};
\node[above=2pt of 2, xshift=0.5cm] {\small \textbf{N: 2}};
\node[above=2pt of 3] {\small \textbf{N: 3}};
\node[above=2pt of 4] {\small \textbf{N: 4}};
\node[above=2pt of 5] {\small \textbf{N: 5}};
\end{tikzpicture}
\vspace{-0.012\textwidth}
\caption{Decision tree for the \texttt{Relaxed Budget Case}.}
\label{fig:decision_tree_relaxed_budget}
\end{figure}
\vspace*{\fill}

\end{landscape}

\begin{figure}[H]
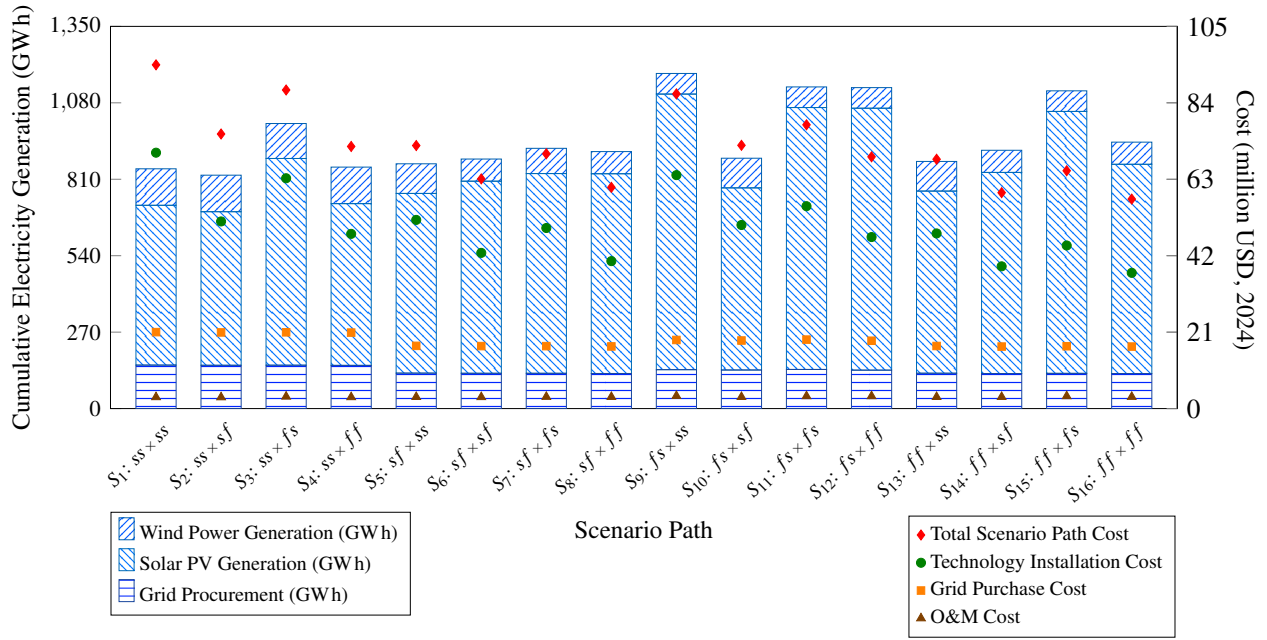

\centering

\caption{Cumulative electricity generation and grid electricity procurement, along with scenario path cost components, for each scenario path under the \texttt{Relaxed Budget Case}.}
\label{fig:scenario_paths_relaxed_budget}
\end{figure}

\vspace{-10pt}

\begin{figure}[H]
\centering
\begin{subfigure}[t]{0.46\textwidth}
\centering
%
\vspace{-5pt}
\caption{Annual electricity supply mix percentages and cost breakdowns for extreme scenario paths in the \texttt{Relaxed Budget Case}.}
\label{fig:extreme_scenario_paths_relaxed_budget}
\end{figure}

\clearpage
\begin{landscape}
  \vspace*{\fill}
    \begin{figure}[H]
    \centering
    %
    };
    \draw ([xshift=-1.05cm]1.south) -- ([xshift=1.5cm]2.north) node [midway, left, shift={(-.45,0)}, edge label] {\footnotesize \textbf{\textit{ss}}};
    \draw ([xshift=0.5cm]1.south) -- ([xshift=0.6cm]3.north) node [midway, left, shift={(-.05,0)}, edge label] {\footnotesize \textbf{\textit{sf}}};
    \draw ([xshift=1.55cm]1.south) -- ([xshift=-0.6cm]4.north) node [midway, right, shift={(.05,0)}, edge label] {\footnotesize \textbf{\textit{fs}}};
    \draw ([xshift=3.15cm]1.south) -- ([xshift=-1cm]5.north) node [midway, right, shift={(.45,0)}, edge label] {\footnotesize \textbf{\textit{ff}}};
    \draw ([xshift=-0.5cm]2.south) -- ([xshift=0.25cm]6.north) node [midway, left, shift={(-.13,0)}, edge label] {\footnotesize \textbf{\textit{ss}}};
    \draw ([xshift=0.2cm]2.south) -- (7.north) node [midway, left, shift={(.02,0)}, edge label] {\footnotesize \textbf{\textit{sf}}};
    \draw ([xshift=0.7cm]2.south) -- (8.north) node [midway, right, shift={(-.02,0)}, edge label] {\footnotesize \textbf{\textit{fs}}};
    \draw ([xshift=1.5cm]2.south) -- (9.north) node [midway, right, shift={(.13,0)}, edge label] {\footnotesize \textbf{\textit{ff}}};
    \draw (3) -- (10.north) node [midway, left, edge label] {\footnotesize \textbf{\textit{ss}}};
    \draw (3) -- (11.north) node [midway, left, edge label] {\footnotesize \textbf{\textit{sf}}};
    \draw (3) -- (12.north) node [midway, right, edge label] {\footnotesize \textbf{\textit{fs}}};
    \draw (3) -- (13.north) node [midway, right, edge label] {\footnotesize \textbf{\textit{ff}}};
    \draw (4) -- (14.north) node [midway, left, edge label] {\footnotesize \textbf{\textit{ss}}};
    \draw (4) -- (15.north) node [midway, left, edge label] {\footnotesize \textbf{\textit{sf}}};
    \draw (4) -- (16.north) node [midway, right, edge label] {\footnotesize \textbf{\textit{fs}}};
    \draw (4) -- (17.north) node [midway, right, edge label] {\footnotesize \textbf{\textit{ff}}};
    \draw (5) -- (18.north) node [midway, left, shift={(-.13,0)}, edge label] {\footnotesize \textbf{\textit{ss}}};
    \draw (5) -- (19.north) node [midway, left, shift={(.05,0)}, edge label] {\footnotesize \textbf{\textit{sf}}};
    \draw (5) -- (20.north) node [midway, right, shift={(-.03,0)}, edge label] {\footnotesize \textbf{\textit{fs}}};
    \draw (5) -- (21.north) node [midway, right, shift={(.13,0)}, edge label] {\footnotesize \textbf{\textit{ff}}};
    \node[above=2pt of 1, xshift=1.05cm] {\small \textbf{N: 1}};
    \node[above=2pt of 2, xshift=0.5cm] {\small \textbf{N: 2}};
    \node[above=2pt of 3] {\small \textbf{N: 3}};
    \node[above=2pt of 4] {\small \textbf{N: 4}};
    \node[above=2pt of 5] {\small \textbf{N: 5}};
    \end{tikzpicture}
    \vspace{-0.012\textwidth}
    \caption{Decision tree for the \texttt{Relaxed Emission Case}.}
    \label{fig:decision_tree_relaxed_emission}
    \end{figure}
  \vspace*{\fill}
\end{landscape}
\clearpage

\begin{figure}[H]
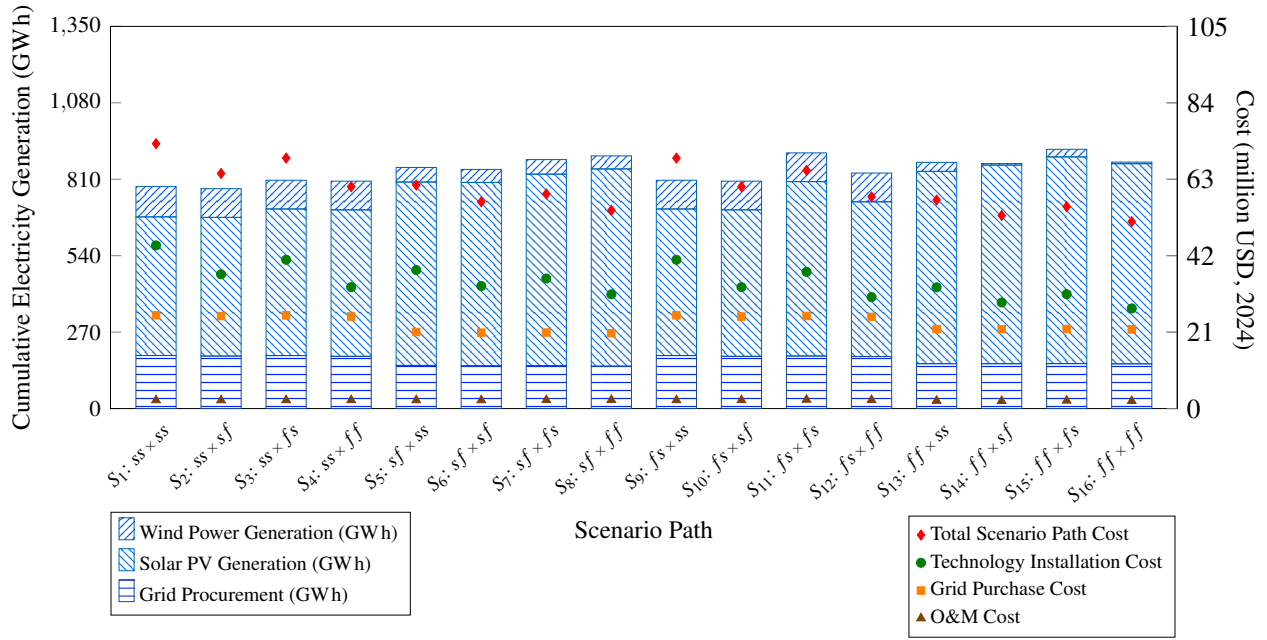

\centering

\caption{Cumulative electricity generation and grid electricity procurement, along with scenario path cost components, for each scenario path under the \texttt{Relaxed Emission Case}.}
\label{fig:scenario_paths_relaxed_emission}
\end{figure}

\vspace{-10pt}

\begin{figure}[H]
\centering
\begin{subfigure}[t]{0.46\textwidth}
\centering
%
\vspace{-5pt}
\caption{Annual electricity supply mix percentages and cost breakdowns for extreme scenario paths in the \texttt{Relaxed Emission Case}.}
\label{fig:extreme_scenario_paths_relaxed_emission}
\end{figure}

\clearpage
\begin{landscape}
  \vspace*{\fill}
    \begin{figure}[H]
    \centering
    %
    };
    \draw ([xshift=-1.05cm]1.south) -- ([xshift=1.5cm]2.north) node [midway, left, shift={(-.45,0)}, edge label] {\footnotesize \textbf{\textit{ss}}};
    \draw ([xshift=0.5cm]1.south) -- ([xshift=0.6cm]3.north) node [midway, left, shift={(-.05,0)}, edge label] {\footnotesize \textbf{\textit{sf}}};
    \draw ([xshift=1.55cm]1.south) -- ([xshift=-0.6cm]4.north) node [midway, right, shift={(.05,0)}, edge label] {\footnotesize \textbf{\textit{fs}}};
    \draw ([xshift=3.15cm]1.south) -- ([xshift=-1cm]5.north) node [midway, right, shift={(.45,0)}, edge label] {\footnotesize \textbf{\textit{ff}}};
    \draw ([xshift=-0.5cm]2.south) -- ([xshift=0.25cm]6.north) node [midway, left, shift={(-.13,0)}, edge label] {\footnotesize \textbf{\textit{ss}}};
    \draw ([xshift=0.2cm]2.south) -- (7.north) node [midway, left, shift={(.02,0)}, edge label] {\footnotesize \textbf{\textit{sf}}};
    \draw ([xshift=0.7cm]2.south) -- (8.north) node [midway, right, shift={(-.02,0)}, edge label] {\footnotesize \textbf{\textit{fs}}};
    \draw ([xshift=1.5cm]2.south) -- (9.north) node [midway, right, shift={(.13,0)}, edge label] {\footnotesize \textbf{\textit{ff}}};
    \draw (3) -- (10.north) node [midway, left, edge label] {\footnotesize \textbf{\textit{ss}}};
    \draw (3) -- (11.north) node [midway, left, edge label] {\footnotesize \textbf{\textit{sf}}};
    \draw (3) -- (12.north) node [midway, right, edge label] {\footnotesize \textbf{\textit{fs}}};
    \draw (3) -- (13.north) node [midway, right, edge label] {\footnotesize \textbf{\textit{ff}}};
    \draw (4) -- (14.north) node [midway, left, edge label] {\footnotesize \textbf{\textit{ss}}};
    \draw (4) -- (15.north) node [midway, left, edge label] {\footnotesize \textbf{\textit{sf}}};
    \draw (4) -- (16.north) node [midway, right, edge label] {\footnotesize \textbf{\textit{fs}}};
    \draw (4) -- (17.north) node [midway, right, edge label] {\footnotesize \textbf{\textit{ff}}};
    \draw (5) -- (18.north) node [midway, left, shift={(-.13,0)}, edge label] {\footnotesize \textbf{\textit{ss}}};
    \draw (5) -- (19.north) node [midway, left, shift={(.05,0)}, edge label] {\footnotesize \textbf{\textit{sf}}};
    \draw (5) -- (20.north) node [midway, right, shift={(-.03,0)}, edge label] {\footnotesize \textbf{\textit{fs}}};
    \draw (5) -- (21.north) node [midway, right, shift={(.13,0)}, edge label] {\footnotesize \textbf{\textit{ff}}};
    \node[above=2pt of 1, xshift=1.05cm] {\small \textbf{N: 1}};
    \node[above=2pt of 2, xshift=0.5cm] {\small \textbf{N: 2}};
    \node[above=2pt of 3] {\small \textbf{N: 3}};
    \node[above=2pt of 4] {\small \textbf{N: 4}};
    \node[above=2pt of 5] {\small \textbf{N: 5}};
    \end{tikzpicture}
    \vspace{-0.012\textwidth}
    \caption{Decision tree for the \texttt{Safety Margin Case}.}
    \label{fig:decision_tree_safety_margin}
    \end{figure}
  \vspace*{\fill}
\end{landscape}
\clearpage

\begin{figure}[H]
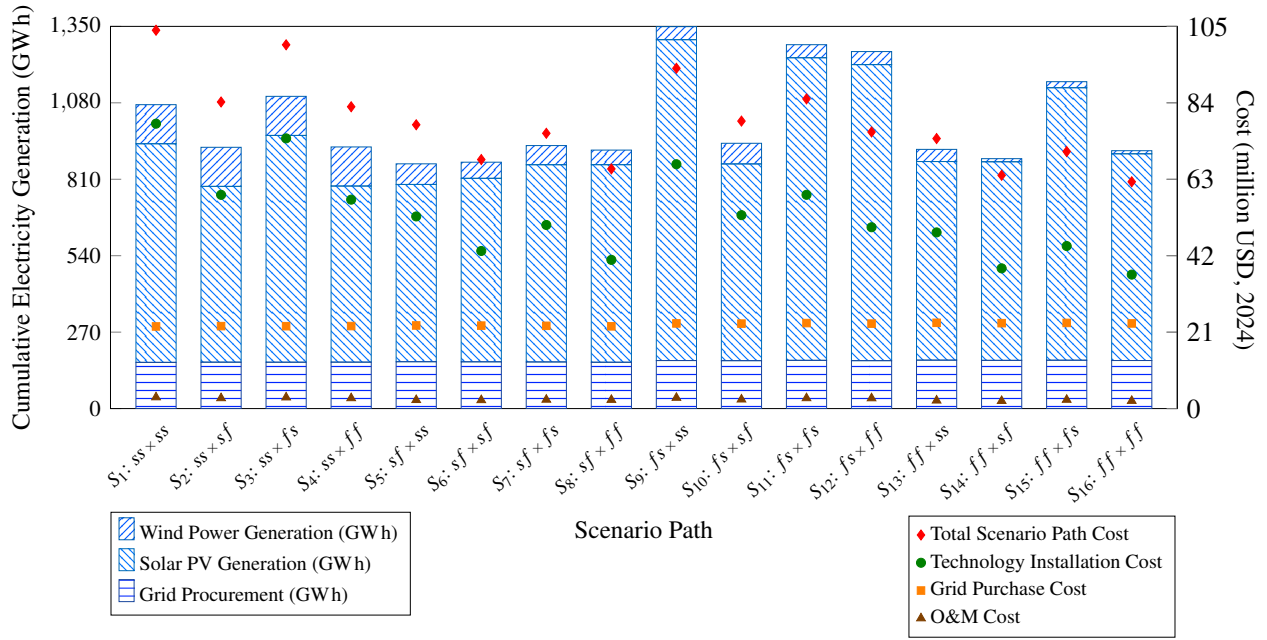

\centering

\caption{Cumulative electricity generation and grid electricity procurement, along with scenario path cost components, for each scenario path under the \texttt{Safety Margin Case}.}
\label{fig:scenario_paths_safety_margin}
\end{figure}

\vspace{-10pt}

\begin{figure}[H]
\centering
\begin{subfigure}[t]{0.46\textwidth}
\centering
%
\vspace{-5pt}
\caption{Annual electricity supply mix percentages and cost breakdowns for extreme scenario paths in the \texttt{Safety Margin Case}}
\label{fig:extreme_scenario_paths_safety_margin}
\end{figure}

\clearpage
\begin{landscape}
  \vspace*{\fill}
    \begin{figure}[H]
    \centering
    %
    };
    \draw ([xshift=-1.05cm]1.south) -- ([xshift=1.5cm]2.north) node [midway, left, shift={(-.45,0)}, edge label] {\footnotesize \textbf{\textit{ss}}};
    \draw ([xshift=0.5cm]1.south) -- ([xshift=0.6cm]3.north) node [midway, left, shift={(-.05,0)}, edge label] {\footnotesize \textbf{\textit{sf}}};
    \draw ([xshift=1.55cm]1.south) -- ([xshift=-0.6cm]4.north) node [midway, right, shift={(.05,0)}, edge label] {\footnotesize \textbf{\textit{fs}}};
    \draw ([xshift=3.15cm]1.south) -- ([xshift=-1cm]5.north) node [midway, right, shift={(.45,0)}, edge label] {\footnotesize \textbf{\textit{ff}}};
    \draw ([xshift=-0.5cm]2.south) -- ([xshift=0.25cm]6.north) node [midway, left, shift={(-.13,0)}, edge label] {\footnotesize \textbf{\textit{ss}}};
    \draw ([xshift=0.2cm]2.south) -- (7.north) node [midway, left, shift={(.02,0)}, edge label] {\footnotesize \textbf{\textit{sf}}};
    \draw ([xshift=0.7cm]2.south) -- (8.north) node [midway, right, shift={(-.02,0)}, edge label] {\footnotesize \textbf{\textit{fs}}};
    \draw ([xshift=1.5cm]2.south) -- (9.north) node [midway, right, shift={(.13,0)}, edge label] {\footnotesize \textbf{\textit{ff}}};
    \draw (3) -- (10.north) node [midway, left, edge label] {\footnotesize \textbf{\textit{ss}}};
    \draw (3) -- (11.north) node [midway, left, edge label] {\footnotesize \textbf{\textit{sf}}};
    \draw (3) -- (12.north) node [midway, right, edge label] {\footnotesize \textbf{\textit{fs}}};
    \draw (3) -- (13.north) node [midway, right, edge label] {\footnotesize \textbf{\textit{ff}}};
    \draw (4) -- (14.north) node [midway, left, edge label] {\footnotesize \textbf{\textit{ss}}};
    \draw (4) -- (15.north) node [midway, left, edge label] {\footnotesize \textbf{\textit{sf}}};
    \draw (4) -- (16.north) node [midway, right, edge label] {\footnotesize \textbf{\textit{fs}}};
    \draw (4) -- (17.north) node [midway, right, edge label] {\footnotesize \textbf{\textit{ff}}};
    \draw (5) -- (18.north) node [midway, left, shift={(-.13,0)}, edge label] {\footnotesize \textbf{\textit{ss}}};
    \draw (5) -- (19.north) node [midway, left, shift={(.05,0)}, edge label] {\footnotesize \textbf{\textit{sf}}};
    \draw (5) -- (20.north) node [midway, right, shift={(-.03,0)}, edge label] {\footnotesize \textbf{\textit{fs}}};
    \draw (5) -- (21.north) node [midway, right, shift={(.13,0)}, edge label] {\footnotesize \textbf{\textit{ff}}};
    \node[above=2pt of 1, xshift=1.05cm] {\small \textbf{N: 1}};
    \node[above=2pt of 2, xshift=0.5cm] {\small \textbf{N: 2}};
    \node[above=2pt of 3] {\small \textbf{N: 3}};
    \node[above=2pt of 4] {\small \textbf{N: 4}};
    \node[above=2pt of 5] {\small \textbf{N: 5}};
    \end{tikzpicture}
    \vspace{-0.012\textwidth}
    \caption{Decision tree for the \texttt{Low Turbine Price Case}.}
    \label{fig:decision_tree_low_turbine_price}
    \end{figure}
  \vspace*{\fill}
\end{landscape}
\clearpage

\begin{figure}[H]
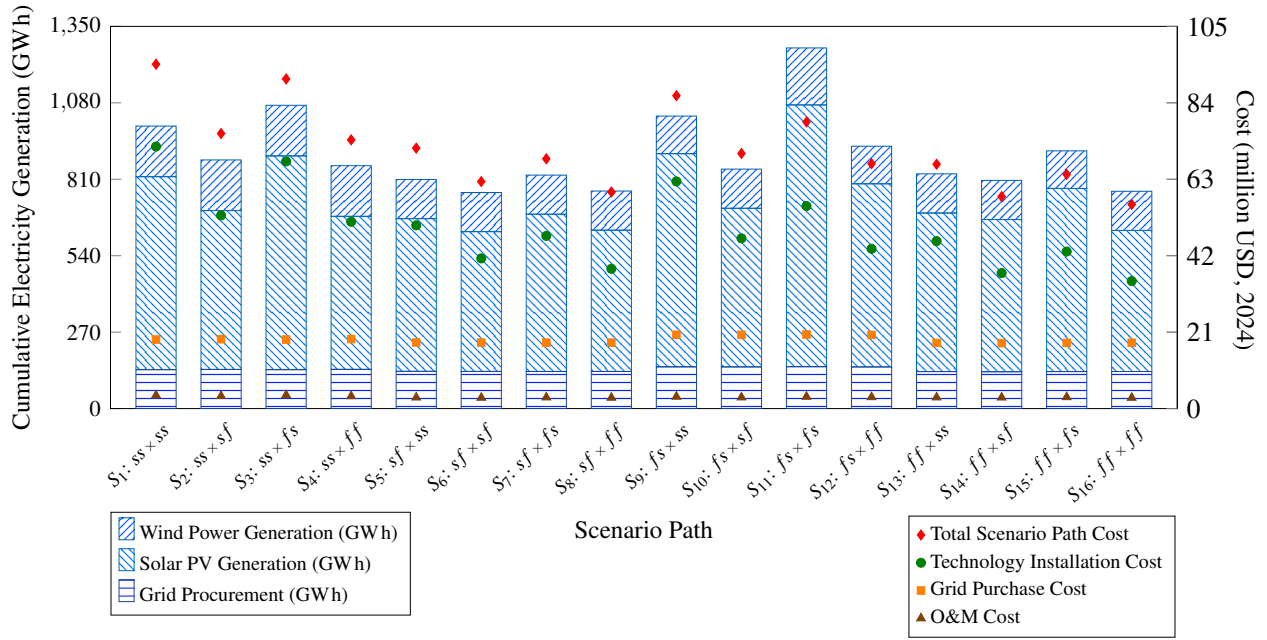

\centering

\caption{Cumulative electricity generation and grid electricity procurement, along with scenario path cost components, for each scenario path under the \texttt{Low Turbine Price Case}.}
\label{fig:scenario_paths_low_turbine_price}
\end{figure}

\vspace{-10pt}

\begin{figure}[H]
\centering
\begin{subfigure}[t]{0.46\textwidth}
\centering
%
\vspace{-5pt}
\caption{Annual electricity supply mix percentages and cost breakdowns for extreme scenario paths in the \texttt{Low Turbine Price Case}.}
\label{fig:extreme_scenario_paths_low_turbine_price}
\end{figure}

\end{document}